\documentclass[british,english,envcountsect, envcountsame, envcountreset, runningheads]{svjour3}
\usepackage[T1]{fontenc}
\usepackage[latin9]{inputenc}
\usepackage{color}
\usepackage{babel}
\usepackage{units}
\usepackage{url}
\usepackage{amsmath}
\usepackage{amssymb}
\usepackage{cancel}
\usepackage[unicode=true,pdfusetitle,
 bookmarks=true,bookmarksnumbered=false,bookmarksopen=false,
 breaklinks=true,pdfborder={0 0 0},pdfborderstyle={},backref=false,colorlinks=true]
 {hyperref}

\makeatletter
\numberwithin{equation}{section}
\newenvironment{svmultproof}{\begin{proof}}{\qed\end{proof}}

\usepackage{graphicx}

\makeatother

\addto\captionsbritish{}
\addto\captionsbritish{}
\addto\captionsbritish{}
\addto\captionsbritish{}
\addto\captionsbritish{}
\addto\captionsbritish{}
\addto\captionsbritish{}
\addto\captionsenglish{}
\addto\captionsenglish{}
\addto\captionsenglish{}
\addto\captionsenglish{}
\addto\captionsenglish{}
\addto\captionsenglish{}
\addto\captionsenglish{}

\begin{document}

\title{Model reduction of non-densely defined piecewise-smooth systems in
Banach spaces}

\titlerunning{Model reduction of PWS systems in Banach spaces}

\author{Robert Szalai}

\date{9/10/2018}

\institute{Department of Engineering Mathematics, University of Bristol, Merchant
Venturers Building, Woodland Road, BS8 1UB, UK, \email{r.szalai@bristol.ac.uk}}
\maketitle
\begin{abstract}
In this paper a model reduction technique is introduced for piecewise-smooth
(PWS) vector fields, whose trajectories fall into a Banach space,
but the domain of definition of the vector fields is a non-dense subset
of the Banach space. The vector fields depend on a parameter that
can assume different discrete values in two parts of the phase space
and a continuous family of values on the boundary that separates the
two parts of the phase space. In essence the parameter parametrizes
the possible vector fields on the boundary. The problem is to find
one or more values of the parameter so that the solution of the PWS
system on the boundary satisfies certain requirements. In this paper
we require continuous solutions. Motivated by the properties of applications,
we assume that when the parameter is forced to switch between the
two discrete values, trajectories become discontinuous. Discontinuous
trajectories exist in systems whose domain of definition is non-dense.
It is shown that under our assumptions the trajectories of such PWS
systems have unique forward-time continuation when the parameter of
the system switches. A finite-dimensional reduced order model is constructed,
which accounts for the discontinuous trajectories. It is shown that
this model retains uniqueness of solutions and other properties of
the original PWS system. The model reduction technique is illustrated
on a nonlinear bowed string model.

\end{abstract}

\section{Introduction}

The purpose of model reduction is to extract the essence of a complex
model, disregarding details that are irrelevant to a specific application.
Depending on the question asked from the model, different kinds of
model reduction are required. In many cases, only qualitative predictions
are needed, where low order analytically solvable models, such as
normal forms used in bifurcation theory \cite{KuznetsovBook}, are
useful. In other cases, the reduced order model has to be solved numerically
with a specified accuracy using constrained computational resources
\cite{benner2017model}. Similar to model reduction, any numerical
scheme that solves a continuum problem, such as finite elements, spectral
collocation or finite differences, turns an infinite-dimensional continuous-time
problem into a finite-dimensional problem. A numerical scheme, however
tends to emphasize quantitative accuracy, which might miss some qualitative
features, such as differentiability of solutions. In this paper we
focus on the qualitative properties of solutions of piecewise-smooth
(PWS) systems, with applications to numerical schemes and reduced
order models in mind.

For smooth systems there are rigorous ways to obtain reduced order
models. Center manifold reduction \cite{carr2012applications} about
an invariant set, such as an equilibrium or periodic orbit, captures
the slowest dynamics and can be used to study bifurcations, regardless
of the dimensionality of the system \cite{KuznetsovBook}. In multiple
time-scale systems \cite{kuehn2015multiple} attracting slow manifolds
that contain dynamics much slower than the rest of the system can
be used to obtain reduced order models.

This paper discusses model reduction for infinite-dimensional systems
that are piecewise smooth. The theory of PWS systems is summarised
in \cite{FilippovBook}, which contains the basic definitions and
results on existence of solutions in finite dimensions. There are
numerous applications of PWS systems, where discontinuities are essential
to the model or where rapid variations of the vector field over small
regions of the phase space naturally lead to discontinuous approximations.
Some applications in finite dimensions include neuron models with
resetting \cite{Coombes2012,IzhikevichNeuron2003}, DC-DC converters
\cite{diBernardoDCDC}, network dynamics, \cite{Danca20021813,DELELLIS20151},
friction oscillators \cite{Popp,SzalaiPolygons}, gene regulatory
networks \cite{GLASS1973103,MESTL1995291} and so on. We consider
the special case of differential equations that are discontinuous
along a codimension-one hypersurface of their phase space, called
the switching manifold. We assume that the phase space is a Banach
space and that the domain of definition of the differential equation
is not dense. 

In contrast to smooth systems, center manifolds or slow manifolds
that continue through switching manifolds do not exist for PWS systems.
In general, the dynamics of singularly perturbed PWS systems cannot
be reduced to an invariant manifold, because small scale instabilities
persist as the perturbation vanishes \cite{SieberSingular}. For a
special class of PWS systems, slow manifolds with similar properties
to smooth systems exist \cite{Fridman2002,cardin_da_silva_teixeira_2013,CARDIN20151166}.
It is also possible to find equivalents of invariant manifolds which
allow model reduction by considering the dynamics on the invariant
manifold. Invariant cones can be found in systems with equilibria
on the switching manifold \cite{WEISS20121895,WEISS201515}. Invariant
polygons may also appear when an unstable focus type periodic orbit
interacts with discontinuities of the vector field \cite{SzalaiPolygons},
which leads to periodic or chaotic dynamics \cite{SzalaiTongues}. 

In infinite dimensions, the theory of PWS systems is focused on sliding
mode control \cite{ORLOV1987} and PWS delay equations \cite{SieberRelay,Londono2012483}.
Sliding mode control applies a discontinuous control signal to a plant,
in order to restrict the system onto an engineered hypersurface with
a prescribed dynamics. The main objective of sliding mode control
is to establish conditions that guarantee the prescribed dynamics.
The results in this area concern systems that are densely defined
on reflexive Banach spaces \cite{Levaggi2002167,Levaggi2002508},
which suggests that these systems are similar to finite dimensional
PWS systems.

In this paper we relax the assumption of a dense domain of definition
and not surprisingly we find different dynamics to what has been studied
before. For this class of systems we are able to prove uniqueness
of solutions and also construct a reduced order model. One consequence
of the non-dense domain is the existence of discontinuous solutions,
which is just the inverse of the Hille-Yosida theorem \cite{Pazy}:
trajectories of a linear autonomous system (as described by a semigroup)
are strongly (or weakly) continuous if and only if the infinitesimal
generator is closed and densely defined. The relevant mathematics
describing our class of systems is the non-autonomous generalization
of integrated semigroups \cite{Neubrander1988,DaPrato1992}. To illustrate
that our class of systems are necessary to describe physical phenomena
we refer to \cite{McIntyre197993}. In \cite{McIntyre197993} the
authors have noticed that the measured impulse response function of
a string has a discontinuity in the velocity component, which is manifest
of the non-dense domain and that the initial condition is outside
of the closure of the domain. This is shown later in the paper for
the relevant mathematical model. Crucially, accounting for the discontinuity
of the impulse response explains the observed asymmetric hysteresis
of the stick-slip motion that causes `flattening' of notes when the
string is bowed in a certain way. The discontinuity of the impulse
response is exactly the property that allows us to carry out model
reduction.

The outline of the paper is as follows. We first carry out model reduction
on a simple linear example to illustrate each step of the process,
but without a rigorous justification of the steps. In section 3, we
review basic classes of PWS systems and highlight some cases where
solutions may be non-unique. Section 4 describes the model reduction
process in a general setting and shows that uniqueness of solutions
and some other properties carry over to the reduced order model. Section
5 describes a nonlinear example, the classical example of the bowed
string, which highlights the significance that nonlinearity plays
in the reduction process and uncovers some possibly surprising results
that were not known about friction oscillators.

\section{\label{sec:TrivialExample}The reduction procedure through an example}

To provide a straightforward template for the model reduction procedure
we take an idealized linear bowed string model with a single contact
point and systematically apply our abstract procedure without rigorous
justification. The list of steps is found at the start of section
\ref{sec:NLS}. Let us consider the equation of motion of a linear
bowed string
\begin{equation}
\ddot{u}(\xi,t)=u''(\xi,t),\;u(0,t)=u(1,t)=0,\,u'(\xi^{\star}-,t)-u'(\xi^{\star}+,t)=\lambda,\label{eq:Trivial String}
\end{equation}
where $u(\xi,t)$ is the scalar valued displacement of the string,
$t\in\mathbb{R}$ is time and $\xi\in[0,1]$ is the spatial coordinate
along the string; $\lambda\in[-1,1]$ is the force applied at the
contact point $\xi^{\star}$, prime denotes differentiation with respect
to $\xi$ and dot with respect to $t$; $u'(\xi^{\star}-,t)$, $u'(\xi^{\star}+,t)$
denote the left and the right derivative at $\xi^{\star}$, respectively.
The value of parameter $\lambda$ is given by 
\begin{equation}
\lambda=\begin{cases}
1 & h>0\\
-1 & h<0
\end{cases},\;\text{where}\;h=v_{0}-\dot{u}(\xi^{\star},t),\label{eq:Trivial SWcond}
\end{equation}
where $v_{0}$ is the speed of the bow and $\lambda$ represents the
friction force between the bow and the string. In general, $h$ is
a smooth scalar valued function of the state variables and is called
the \emph{switching function}. The equation $h=0$ implicitly defines
a surface in the phase space of (\ref{eq:Trivial String}) which is
called the \emph{switching manifold} \cite{diBernardoBook}. Note
that $\lambda$ is not defined for $h=0$ by equation (\ref{eq:Trivial SWcond}).
We assume that all values of $\lambda\in[-1,1]$ are possible when
$h=0$ and therefore the model is a differential inclusion \cite{Smirnov2002}.
Later on we will find a unique value for $\lambda$ using the condition
that the functions $u(\xi^{\star},\cdot)$ and $\dot{u}(\xi^{\star},\cdot)$,
that is, the solution of (\ref{eq:Trivial String}) and (\ref{eq:Trivial SWcond})
evaluated at $\xi=\xi^{\star}$, must be continuous in time.

We now consider the case when $\lambda$ is constant. For constant
$\lambda$ equation (\ref{eq:Trivial String}) has an equilibrium.
The equilibrium shape of the string for $\lambda=1$ is given by 
\begin{equation}
u_{0}(\xi)=(1-\xi^{\star})\xi-(\xi-\xi^{\star})H(\xi-\xi^{\star}),\label{eq:Trivial Equilibrium}
\end{equation}
where $H$ is the Heaviside function. Due to linearity, for a fixed
$\lambda$ the equilibrium is then $\lambda u_{0}(\xi)$. It is known
that free vibrations of a string, that is, the solutions of (\ref{eq:Trivial String})
with constant $\lambda$ can be written as 
\begin{equation}
u(\xi,t)=\lambda u_{0}(\xi)+\sum_{k=1}^{\infty}\sin k\pi\xi\left(a_{k}\sin k\pi t+b_{k}\cos k\pi t\right),\label{eq:Trivial AllModes}
\end{equation}
where $a_{k}$ and $b_{k}$ are determined from initial conditions
\cite[section 8.2]{rao2016mechanical}. We now consider solutions
for which $a_{k}=b_{k}=0$ for $k>1$ in (\ref{eq:Trivial AllModes}).
The remaining two parameters $a_{1}$, $b_{1}$ describe a set of
solutions that are restricted to a two dimensional invariant manifold,
which we denote by $\mathcal{M}$. For this set of solutions we denote
the displacement of the string at $\xi^{\star}$ by $y(t)=\lambda y^{\star}+\left(a_{1}\sin\pi t+b_{1}\cos\pi t\right)\sin\pi\xi^{\star}$,
where $y^{\star}$ is yet to be defined. The value $y(t)$ can be
used to recover the displacement of the whole string as
\begin{equation}
u_{\mathcal{M}}(y(t),\lambda;\xi)=\lambda u_{0}(\xi)+\left(y(t)-\lambda y^{\star}\right)\frac{\sin\pi\xi}{\sin\pi\xi^{\star}}.\label{eq:Trivial Immersion}
\end{equation}
Expression (\ref{eq:Trivial Immersion}) is the immersion of the invariant
manifold into the configuration space, but not into the full phase
space, which owing to the second order time derivative in (\ref{eq:Trivial String})
should also contain velocities. Note that any value of $y^{\star}$
gives the same manifold, $y^{\star}$ only influences the parametrization
of the manifold. To fix $y^{\star}$ we require that the manifold
does not move in the tangential direction of $\mathcal{M}$ when $\lambda$
changes. This means that the two partial derivatives $\frac{\partial}{\partial\lambda}u_{\mathcal{M}}(y,\lambda;\xi)$
and $\frac{\partial}{\partial y}u_{\mathcal{M}}(y,\lambda;\xi)$ must
be perpendicular, that is,
\begin{equation}
\int_{0}^{1}\frac{\partial}{\partial\lambda}u_{\mathcal{M}}(y,\lambda;\xi)\frac{\partial}{\partial y}u_{\mathcal{M}}(y,\lambda;\xi)\mathrm{d}\xi=0.\label{eq:Trivial Restriction}
\end{equation}
Solving equation (\ref{eq:Trivial Restriction}) for $y^{\star}$
we get 
\[
y^{\star}=\frac{2}{\pi^{2}}\sin^{2}\pi\xi^{\star}.
\]

Substituting the immersion (\ref{eq:Trivial Immersion}) of manifold
$\mathcal{M}$ into (\ref{eq:Trivial String}) and (\ref{eq:Trivial SWcond})
while assuming that $\lambda$ is constant, we get
\begin{equation}
\left.\begin{array}{l}
\ddot{y}+\pi^{2}(y-\lambda y^{\star})=0\\
\lambda=\begin{cases}
1 & h>0\\
-1 & h<0
\end{cases}
\end{array}\right\} ,\label{eq:Trivial 1DOF}
\end{equation}
where $h=v_{0}-\dot{y}$. Equation (\ref{eq:Trivial 1DOF}) has the
form of a common PWS system, which is widely used as a reduced order
model of (\ref{eq:Trivial String}). However, the assumption that
$\lambda$ is constant does not hold when $h=0$, therefore we consider
(\ref{eq:Trivial 1DOF}) a skeleton of a more accurate description
and call (\ref{eq:Trivial 1DOF}) the \emph{skeleton model}.

The skeleton model (\ref{eq:Trivial 1DOF}) is a typical friction
oscillator and therefore can be solved using techniques known from
mechanics. At \emph{stick}, when $h=0$, $\lambda$ is not explicitly
defined by the skeleton model (\ref{eq:Trivial 1DOF}), instead we
need to argue the following. If $h=0$ on some interval of time, then
$\dot{y}=v_{0}$ on this interval, and consequently, $\ddot{y}=0$.
Substituting the conclusion of this argument into the first line of
(\ref{eq:Trivial 1DOF}) yields 
\begin{equation}
\lambda=\frac{y}{y^{\star}}.\label{eq:Trivial Closure}
\end{equation}
We only allow $\lambda\in[-1,1]$, therefore if the result of (\ref{eq:Trivial Closure})
is outside of the interval $[-1,1]$, $\lambda$ simply swaps from
$1$ to $-1$ or vice versa. The argument made to find (\ref{eq:Trivial Closure})
is equivalent to \emph{Filippov's closure}, which is summarized in
section \ref{sec:FilUtkin}.

\begin{figure}
\begin{centering}
\includegraphics[width=0.7\linewidth]{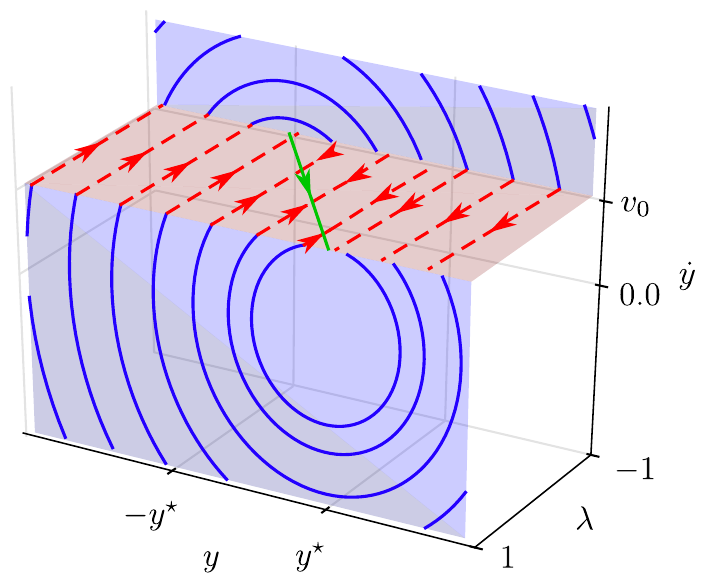}
\par\end{centering}
\caption{\label{fig:Trivial Fil Phase}Phase portrait of the skeleton model
(\ref{eq:Trivial 1DOF}). Solid lines represent solutions where $\lambda$
is continuous and dashed lines represent solutions where $\lambda$
switches between two values. The horizontal plane with $\dot{y}=v_{0}$
contains sticking solutions. There is only one sticking solution represented
by the solid line, which is also marked with an arrowhead denoting
the direction of time.}
\end{figure}

The phase portrait of the skeleton model (\ref{eq:Trivial 1DOF})
can be seen in figure \ref{fig:Trivial Fil Phase}. We focus on the
dynamics at stick, which occurs on the switching manifold, highlighted
by the horizontal red shaded plane. The dashed red lines correspond
to discontinuities in $\lambda$. The solid green line on the horizontal
red shaded plane is the stick solution, where the friction force $\lambda$
grows with a constant rate with respect to $t$ and $y$ until it
reaches the limit $\lambda=\pm1$ and slip ensues.

By assuming constant $\lambda$, we made an error when the relative
velocity between the bow and the string becomes zero, i.e., $\dot{y}=v_{0}$,
because contrary to the assumption $\lambda$ is not constant and
can even jump between $-1$ and $1$ or from $\pm1$ to the value
given by (\ref{eq:Trivial Closure}). The desire to correct for this
error is the subject of the paper, because this is the source of the
qualitative discrepancy between solutions of equations (\ref{eq:Trivial String}),
(\ref{eq:Trivial SWcond}) and the skeleton model (\ref{eq:Trivial 1DOF}).
To account for the error made, the solution of equations (\ref{eq:Trivial String})
and (\ref{eq:Trivial SWcond}) is now written as 
\begin{equation}
u(\xi,t)=u_{\mathcal{M}}(y(t),\lambda(t);\xi)+w(\xi,t),\label{eq:Trivial ExactSolution}
\end{equation}
where $w(\xi,t)$ is a correction term. Assuming that $y(t)$ satisfies
(\ref{eq:Trivial 1DOF}) and substituting (\ref{eq:Trivial ExactSolution})
into (\ref{eq:Trivial String}) without assuming constant $\lambda$
we get the governing equation of the correction term
\begin{equation}
\ddot{w}(\xi,t)=w^{\prime\prime}(\xi,t)-\ddot{\lambda}\left(u_{0}(\xi)-y^{\star}\frac{\sin\pi\xi}{\sin\pi\xi^{\star}}\right),\label{eq:Trivial Pert Wave}
\end{equation}
which describes how the dynamics depart from the invariant manifold
when $\lambda$ varies. Starting from the invariant manifold $\mathcal{M}$,
we have initial conditions
\begin{equation}
w(\xi,0)=\dot{w}(\xi,0)=0\;\forall\xi\in[0,1].\label{eq:Trivial Pert Wave IC}
\end{equation}
The choice of $y^{\star}$ dictated by equation (\ref{eq:Trivial Restriction})
guarantees that the correction term as described by (\ref{eq:Trivial Pert Wave})
does not include vibrations with the first natural frequency of the
string. (In abstract terms, $w$ is restricted to the invariant normal
bundle of $\mathcal{M}$.) Consequently, the skeleton model (\ref{eq:Trivial 1DOF})
does not require any correction even when $\lambda$ varies. Instead,
the switching function (\ref{eq:Trivial SWcond}) needs to be revised
by substituting the corrected solution (\ref{eq:Trivial ExactSolution}).
If we use (\ref{eq:Trivial ExactSolution}) in equation (\ref{eq:Trivial SWcond})
we can write the switching function as
\begin{equation}
h=v_{0}-\dot{y}-\dot{\lambda}\left((1-\xi^{\star})\xi^{\star}-y^{\star}\right)-\dot{w}(\xi^{\star},t).\label{eq:Trivial SWcond Pert}
\end{equation}
The difficulty of evaluating equation (\ref{eq:Trivial SWcond Pert})
lies with solving (\ref{eq:Trivial Pert Wave}), which we carry out
in appendix \ref{sec:Trivial Pert Solution} in detail.

\begin{figure}
\begin{centering}
\includegraphics[width=1\linewidth]{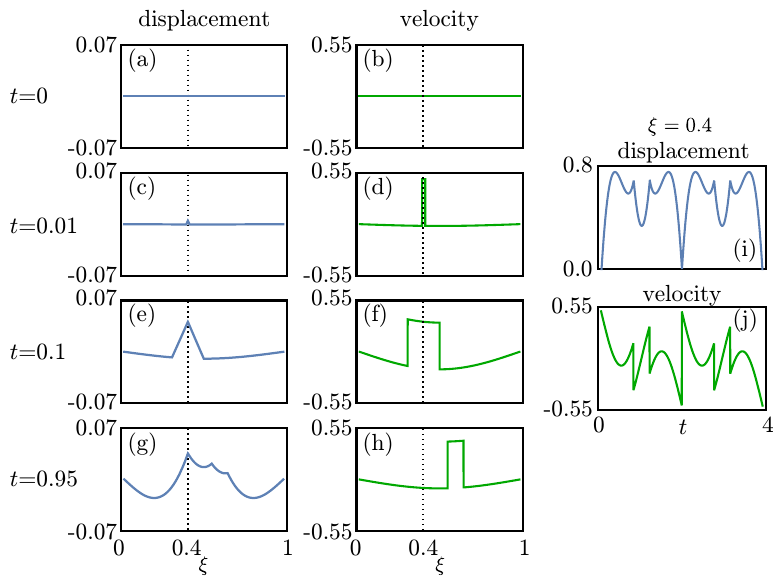}
\par\end{centering}
\caption{\label{fig:Intro String Solution}The solution of equation (\ref{eq:Trivial Pert Wave})
with zero initial conditions and with forcing $\lambda(t)=H(t)$ .
The velocity of the string develops discontinuities both in space
and time.}
\end{figure}

We now illustrate how a discontinuity of $\lambda$ leads to a jump
in the velocity $\dot{w}(\xi^{\star},t)$. For this we assume that
$\lambda(t)=H(t)$ in equation (\ref{eq:Trivial Pert Wave}). The
solution of equation (\ref{eq:Trivial Pert Wave}) can be seen in
figure \ref{fig:Intro String Solution}, with initial conditions (\ref{eq:Trivial Pert Wave IC}).
By comparing figures \ref{fig:Intro String Solution}(b) and \ref{fig:Intro String Solution}(d),
it can be seen that the velocity $\dot{w}(\xi^{\star},t)$ at $t=0$
has a discontinuity, whose gap is proportional to the jump in $\lambda$.
The time history of this velocity in figure \ref{fig:Intro String Solution}(j)
has further discontinuities. Discontinuities for $t>0$ are due to
reflections at the boundaries and they are specific to this example
that lacks damping. Typically, wave dispersion or damping that is
present in other mechanical systems would destroy discontinuities
for $t>0$, but not at $t=0$. At $t=0$, we have 
\begin{equation}
\lim_{t\downarrow0}\dot{w}(\xi^{\star},t)-\dot{w}(\xi^{\star},0)=\frac{1}{2},\label{eq:Trivial Discontinuity Gap}
\end{equation}
which we call the \emph{normal discontinuity gap}. Due to the linearity
of equation (\ref{eq:Trivial Pert Wave}), any jump in $\lambda$
is translated into a discontinuity of the velocity $\dot{w}(\xi^{\star},t)$.
This velocity jump also appears in the switching function (\ref{eq:Trivial SWcond Pert}),
which makes a qualitative difference between the dynamics of the original
model and the skeleton model (\ref{eq:Trivial 1DOF}) on the switching
manifold as we show below.

After solving equation (\ref{eq:Trivial Pert Wave}) on the interval
$0\le t<\min\left(2\xi^{\star},2-2\xi^{\star}\right)$, before any
discontinuity is reflected back to $\xi=\xi^{\star}$, we find that
the switching function (\ref{eq:Trivial SWcond Pert}) in equation
(\ref{eq:Trivial 1DOF}) becomes
\begin{equation}
h=v_{0}-\dot{y}-\frac{1}{2}\lambda-\kappa+\frac{1}{2}\lambda(0),\label{eq:Trivial SWcond RED}
\end{equation}
where $\kappa$ is a variable satisfying the initial value problem
\begin{equation}
\ddot{\kappa}=\pi^{2}\left(y^{\star}\dot{\lambda}-\kappa\right),\;\kappa(0)=0,\,\dot{\kappa}(0)=0.\label{eq:Trivial Kappa Dyn}
\end{equation}
Note that by using $h$ as defined by (\ref{eq:Trivial SWcond RED})
in equation (\ref{eq:Trivial 1DOF}) we get an exact representation
of the dynamics of the original problem (\ref{eq:Trivial String})
and (\ref{eq:Trivial SWcond}) on the time interval $0\le t<\min\left(2\xi^{\star},2-2\xi^{\star}\right)$.
The valid time interval can be extended to any length by considering
the full solution of (\ref{eq:Trivial Pert Wave}) derived in appendix
\ref{sec:Trivial Pert Solution}.

The switching function (\ref{eq:Trivial SWcond RED}) depends on $\lambda$,
because of the presence of the normal discontinuity gap (\ref{eq:Trivial Discontinuity Gap}).
Therefore, when $h=0$, $\lambda$ can be solved for, without any
closure rule, such as Filippov's or Utkin's (see Section \ref{sec:FilUtkin}).
In our case, solving the equation $h=0$ for $\lambda$ yields
\begin{equation}
\lambda=\lambda(0)+2\left(v_{0}-\dot{y}-\kappa\right).\label{eq:Trivial LambdaSol}
\end{equation}

A non-zero normal discontinuity gap, as calculated in equation (\ref{eq:Trivial Discontinuity Gap}),
turns  the skeleton model (\ref{eq:Trivial 1DOF}) at stick into
an index-1 differential algebraic equation. When gathering all dynamic
equations at stick we get 
\begin{equation}
\left.\begin{array}{c}
\ddot{y}+\pi^{2}(y-\lambda y^{\star})=0\\
\lambda=\lambda(0)+2\left(v_{0}-\dot{y}-\kappa\right)\\
\ddot{\kappa}=\pi^{2}\left(y^{\star}\dot{\lambda}-\kappa\right)
\end{array}\right\} .\label{eq:Trivial DAE}
\end{equation}
By definition, an index-1 differential algebraic equation can be turned
into an ordinary differential equation by differentiating the algebraic
equation (\ref{eq:Trivial LambdaSol}) once, which for equation (\ref{eq:Trivial DAE})
of the bowed string at stick gives 
\begin{equation}
\left.\begin{array}{c}
\ddot{y}+\pi^{2}(y-\lambda y^{\star})=0\\
\dot{\lambda}=2\pi^{2}(y-\lambda y^{\star})-2\dot{\kappa}\\
\ddot{\kappa}+\pi^{2}\left(2y^{\star}\dot{\kappa}+\kappa\right)=2\pi^{4}y^{\star}(y-\lambda y^{\star})
\end{array}\right\} .\label{eq:Trivial DAE Traf}
\end{equation}
Note that the differentiation of (\ref{eq:Trivial LambdaSol}) also
transforms the stick constraint $h=0$ into $\dot{h}=0$, therefore
(\ref{eq:Trivial DAE Traf}) is valid for initial conditions that
satisfy $h=0$ with $\lambda\in[-1,1]$.

We can now put together the whole model with the correction into a
single system
\begin{equation}
\left.\begin{array}{l}
\ddot{y}+\pi^{2}(y-\lambda y^{\star})=0\\
\ddot{\kappa}=\pi^{2}\left(y^{\star}\dot{\lambda}-\kappa\right)\\
\lambda=\begin{cases}
1 & h+\lambda>1\\
\lambda(0)+2\left(v_{0}-\dot{y}-\kappa\right) & -1<h+\lambda<1\\
-1 & h+\lambda<-1
\end{cases}
\end{array}\right\} .\label{eq:Trivial Reduced Model}
\end{equation}
We call equation (\ref{eq:Trivial Reduced Model}) the \emph{reduced
order model} of the initial problem (\ref{eq:Trivial String}) and
(\ref{eq:Trivial SWcond}), because it exactly reproduces the dynamics
for initial conditions in $\mathcal{M}$ and for the time interval
$0\le t<\min\left(2\xi^{\star},2-2\xi^{\star}\right)$. Appendix \ref{sec:Trivial Pert Solution}
shows that the valid time interval can be extended to any length by
including delayed values of $\lambda$ in the switching function.
It is noteworthy that instead of the two phase space regions defined
by (\ref{eq:Trivial SWcond}), the reduced order model (\ref{eq:Trivial Reduced Model})
has three regions, where the dynamics is defined. The additional phase
space region corresponds to the stick phase of motion, which has its
own regular dynamics. This dynamics follows from the assumption that
the velocity $y_{2}$ is continuous in time and there is a normal
discontinuity gap, i.e., the correction term $w$ is discontinuous
at $t=0$, when $\lambda=H(t)$.

\begin{figure}
\begin{centering}
\includegraphics[width=0.7\linewidth]{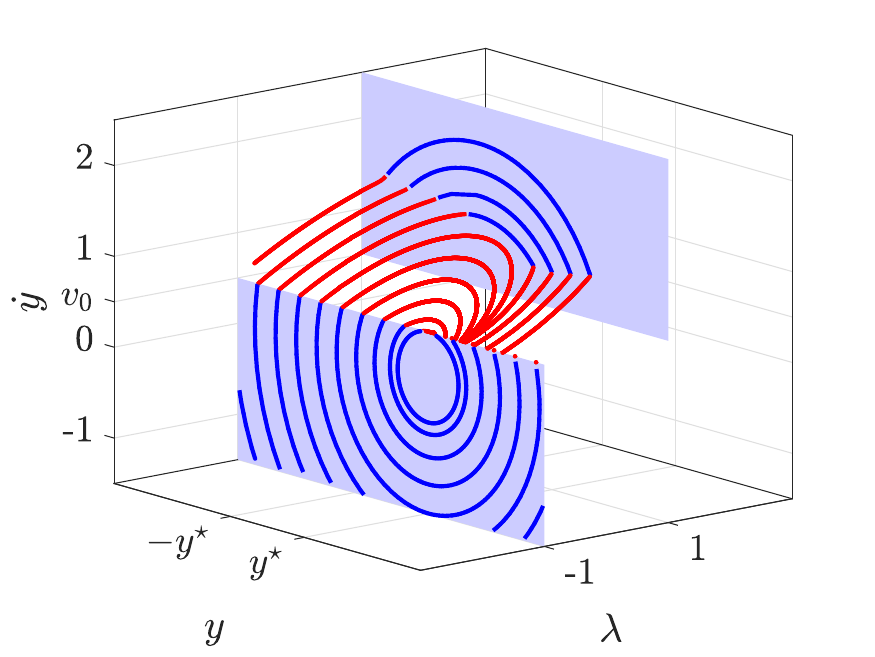}
\par\end{centering}
\caption{\label{fig:Trivial FullSim}Phase portrait of the full model (\ref{eq:Trivial String})
and (\ref{eq:Trivial SWcond}). The simulation was carried out using
an exact representation of the switching function $h$, a further
extension of (\ref{eq:Trivial SWcond RED}) to $t\ge0$, which is
derived in appendix \ref{sec:Trivial Pert Solution}.}
\end{figure}

The phase portrait of the reduced order model (\ref{eq:Trivial Reduced Model})
can be seen in figure \ref{fig:Trivial FullSim}. In the simulation
we have used the result of appendix \ref{sec:Trivial Pert Solution}
to extend the valid time interval to an appropriate length. In comparison
to the skeleton model (\ref{eq:Trivial 1DOF}) shown in figure \ref{fig:Trivial Fil Phase}
the dynamics at stick becomes more complicated. The dynamics in slip
are the same, because $\lambda$ is constant and decoupled from the
rest of the variables due to the choice of immersion (\ref{eq:Trivial Restriction}).
The stick dynamics is now described by the differential equation (\ref{eq:Trivial DAE Traf})
and therefore there is no discontinuity of $\lambda$. Due to the
higher dimensional dynamics that arise from the inclusion of $\kappa$
as a dynamic variable and delayed values of $\lambda$, the dynamics
depicted in figure \ref{fig:Trivial FullSim} is only a projection.
Regardless of the differences, the phase portrait in figure \ref{fig:Trivial FullSim}
appears as a smoothed version of the same dynamics in figure \ref{fig:Trivial Fil Phase},
even though no smoothing or regularization was applied. Furthermore,
to solve the reduced order model (\ref{eq:Trivial Reduced Model})
we did not need an arbitrary closure, such as Filippov's to define
the dynamics at stick, instead the solution followed straight from
the initial problem (\ref{eq:Trivial String}) and (\ref{eq:Trivial SWcond}).

In section \ref{sec:Continuum} we explore a generalization of equations
(\ref{eq:Trivial String}) and (\ref{eq:Trivial SWcond}). We consider
models whose solutions may be norm-discontinuous as illustrated by
the linear string model. Before we embark on the general theory we
recall basic definitions and properties of PWS models.

\section{\label{sec:FilUtkin}Finite dimensional PWS models}

In this section we summarize two commonly used closures of PWS systems.
As it turns out, these common PWS systems are special forms of the
skeleton model to be defined in section \ref{sec:skeleton}. An introduction
to the state-of-the-art can be found in \cite{GlendinngJeffrey2017},
however the book of Filippov \cite{FilippovBook} contains the most
general definitions of PWS systems. Below, we review the cases defined
in \cite[Chapter 2, \S{}4]{FilippovBook}, which are used most commonly
in applications. We avoid cases where the vector field is a set-valued
function \cite[Chapter 2, \S{}5,\S{}6]{FilippovBook}. We also limit
the description to the bi-modal case, where the discontinuity occurs
along a single implicitly defined manifold in the phase space. 
\begin{note}
In addition to various notation for derivatives, in what follows $D$
is also used to denote the Frechet derivative of a function; a subscript
of $D$, such as $D_{k}$ denotes the partial derivative with respect
to the $k$th argument of a function and a superscript such as $D_{k}^{j}$
denotes the $j$th derivative with respect to the $k$th argument.
\end{note}
Let us consider the vector field
\begin{align}
\dot{\boldsymbol{y}} & =\boldsymbol{f}(\boldsymbol{y},\lambda),\label{eq:LowDim}\\
\lambda & =\begin{cases}
1 & \mbox{for}\;h(\boldsymbol{y})>0,\\
-1 & \mbox{for}\;h(\boldsymbol{y})<0,
\end{cases}\label{eq:SwCond}
\end{align}
where either 
\begin{align}
\boldsymbol{f} & \in C^{p}(G\times\{-1,1\},\mathbb{R}^{n})\;\text{or}\label{eq:Review-PointDomain}\\
\boldsymbol{f} & \in C^{p}(G\times[-1,1],\mathbb{R}^{n}),\label{eq:Review-IntervalDomain}
\end{align}
$G$ is a compact and connected subset of $\mathbb{R}^{n}$ and $n,p\in\mathbb{N}^{+}$.
The function $h\in C^{p}(G,\mathbb{R})$ is called the switching function
and its zero set defines the switching manifold 
\begin{equation}
\Sigma=\left\{ \boldsymbol{y}\in G:h(\boldsymbol{y})=0\right\} .\label{eq:SWsurf}
\end{equation}
A solution $\boldsymbol{y}:I\to\mathbb{R}^{n}$ of equations (\ref{eq:LowDim})
and (\ref{eq:SwCond}) is defined on a closed interval of non-zero
length $I\subset\mathbb{R}$. There is no information in equations
(\ref{eq:LowDim}) and (\ref{eq:SwCond}) that helps to deduce a value
for $\lambda$ on $\Sigma$. To explore all possibilities (\ref{eq:LowDim})
can be turned into a differential inclusion
\begin{equation}
\dot{\boldsymbol{y}}\in\overline{\mathrm{co}}\boldsymbol{f}(\boldsymbol{y},[-1,1])\;\text{or}\;\dot{\boldsymbol{y}}\in\overline{\mathrm{co}}\boldsymbol{f}(\boldsymbol{y},\{-1,1\}),\;\boldsymbol{y}\in\Sigma,\label{eq:Review-Inclusion}
\end{equation}
where $\overline{\mathrm{co}}$ denotes the closure of the convex
hull of a set. The existence of solutions of (\ref{eq:Review-Inclusion})
is investigated in \cite[Chapter 2, \S{}7]{FilippovBook}. In this
section we review different definitions of $\lambda$ on $\Sigma$.

We have assumed two possibilities, (\ref{eq:Review-PointDomain})
or (\ref{eq:Review-IntervalDomain}), for the domain of definition
of $\boldsymbol{f}$. The case of (\ref{eq:Review-PointDomain}) is
the minimum necessary to make equations (\ref{eq:LowDim}) and (\ref{eq:SwCond})
consistent. In many applications, such as mechanics, the larger domain
of definition (\ref{eq:Review-IntervalDomain}) is naturally given,
which is useful to define the solutions of (\ref{eq:LowDim}) and
(\ref{eq:SwCond}) on $\Sigma$ as we show later in this section.

The system (\ref{eq:LowDim}) and (\ref{eq:SwCond}) has a unique
solution on an interval of non-zero length for initial condition $\boldsymbol{y}(0)\in G$
if $h(\boldsymbol{y}(0))\neq0$, because $\boldsymbol{f}$ is a smooth
vector field \cite{CoddingtonLevinson}. However for $h(\boldsymbol{y}(0))=0$
the vector field is not defined and one needs to reason how trajectories
continue once they reach $\Sigma$.

The simplest case of a trajectory interacting with $\Sigma$ is when
the trajectory approaches $\Sigma$ transversely from one side and
continues on the other side; this is called crossing. The exact condition
for crossing is that the value of $h$ must change monotonically through
$h=0$ with a non-zero speed along the trajectory at $h=0$. If the
common point of the trajectory with $\Sigma$ is denoted by $\boldsymbol{y}^{\star}$,
the speeds at which $h$ increases are $Dh(\boldsymbol{y}^{\star})\boldsymbol{f}(\boldsymbol{y}^{\star},\pm1)$.
$h$ changes monotonically through $h=0$ with a non-zero speed if
and only if
\begin{equation}
\left(Dh(\boldsymbol{y}^{\star})\boldsymbol{f}(\boldsymbol{y}^{\star},1)\right)\left(Dh(\boldsymbol{y}^{\star})\boldsymbol{f}(\boldsymbol{y}^{\star},-1)\right)>0.\label{eq:Review-CrossingCondition}
\end{equation}
In case of (\ref{eq:Review-CrossingCondition}), there is no need
to define the dynamics on $\Sigma$ because the value of $\lambda$
simply switches between $\pm1$. In our argument we have used $\boldsymbol{f}$
for $\lambda=\pm1$ only, therefore, to resolve crossing, it is sufficient
to assume (\ref{eq:Review-PointDomain}). The subset of $\Sigma$,
where (\ref{eq:Review-CrossingCondition}) holds is called the \emph{crossing
region} and denoted by $\Sigma_{cr}$.

Now we discuss the case when
\begin{equation}
\left(Dh(\boldsymbol{y}^{\star})\boldsymbol{f}(\boldsymbol{y}^{\star},1)\right)\left(Dh(\boldsymbol{y}^{\star})\boldsymbol{f}(\boldsymbol{y}^{\star},-1)\right)<0.\label{eq:Review-SlidingCondition}
\end{equation}
When (\ref{eq:Review-SlidingCondition}) holds, $\Sigma$ is an attractor
in either forward or backward time. This means that trajectories cannot
immediately escape $\Sigma$ once they are on $\Sigma$. The subset
of $\Sigma$, where (\ref{eq:Review-SlidingCondition}) holds and
attracts solutions in forward time is called the \emph{sliding region}
and denoted by $\Sigma_{sl}$. The repelling subset of $\Sigma$,
where (\ref{eq:Review-SlidingCondition}) holds, is called the escaping
region and denoted by $\Sigma_{esc}$.

In case of (\ref{eq:Review-SlidingCondition}), equations (\ref{eq:LowDim})
and (\ref{eq:SwCond}) are not sufficient to define a solution and
an assumption is required that specifies how a trajectory continues
on $\Sigma$. In this paper, we call such an assumption a \emph{closure},
because it completes (\ref{eq:LowDim}) and (\ref{eq:SwCond}). In
the following, we discuss two commonly used closures. The first closure
is attributed to Filippov \cite[Chapter 2, \S{}4, 2.a)]{FilippovBook},
the second closure is due to Utkin \cite{utkin1992sliding}, and also
explored in Filippov's book \cite[Chapter 2, \S{}4, 2.b)]{FilippovBook}.
We note that there is no common terminology in the literature for
closures, various closures have their own name. For example, Filippov's
closure is commonly called Filippov's method and Utkin's closure is
called the equivalent control method, due to its origin in control
theory. There are many possibilities to define a closure, for example,
Filippov explores systems where the closure is not explicitly defined,
but only certain constraints are placed on it \cite[Chapter 4]{FilippovBook}.

\begin{figure}
\begin{centering}
\includegraphics[width=1\linewidth]{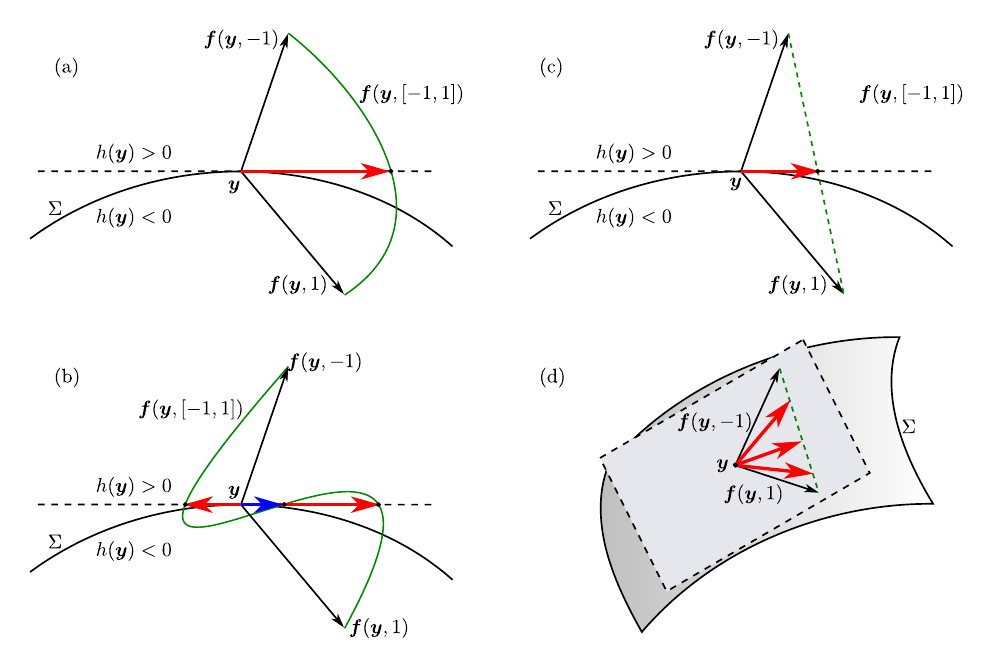}
\par\end{centering}
\caption{\label{fig:FilUtkin}Illustration of Utkin's and Filippov's closure.
(a) The thick red arrow shows the vector field on the switching manifold,
which is the only tangential vector to $\Sigma$ chosen from the family
of vectors $\boldsymbol{f}(\boldsymbol{y},[-1,1])$ as defined by
Utkin's closure. In this example there is a unique element that is
tangent to $\Sigma$. (b) An example, where the family of vectors
$\boldsymbol{f}(\boldsymbol{y},[-1,1])$ has three elements that are
tangent to $\Sigma$, hence no unique solution can be found. (c) Filippov's
convex combination of the vectors $\boldsymbol{f}(\boldsymbol{y},-1)$
and $\boldsymbol{f}(\boldsymbol{y},1)$ is illustrated by the green
dashed line. The only tangential vector in this convex set is the
thick red arrow. (d) Both vectors $\boldsymbol{f}(\boldsymbol{y},-1)$
and $\boldsymbol{f}(\boldsymbol{y},1)$ are tangential to $\Sigma$,
hence all vectors in their convex hull are equally plausible, hence
there is no unique solution.}
\end{figure}

\subsection{\label{subsec:Filippovs-closure}Filippov's closure}

Filippov's closure defines a vector field on $\Sigma$, when condition
(\ref{eq:Review-SlidingCondition}) holds, by interpolating between
the vector fields $\boldsymbol{f}(\boldsymbol{y},\pm1)$, such that
$\Sigma$ becomes an invariant manifold of the new vector field. The
interpolation is carried out as follows. We define a new vector field
\begin{equation}
\dot{\boldsymbol{y}}=\boldsymbol{r}(\boldsymbol{y})+\boldsymbol{b}(\boldsymbol{y})\lambda,\;\lambda\in[-1,1]\label{eq:LowDimFil}
\end{equation}
where
\[
\boldsymbol{r}(\boldsymbol{y})=\frac{\boldsymbol{f}(\boldsymbol{y},1)+\boldsymbol{f}(\boldsymbol{y},-1)}{2},\;\boldsymbol{b}(\boldsymbol{y})=\frac{\boldsymbol{f}(\boldsymbol{y},1)-\boldsymbol{f}(\boldsymbol{y},-1)}{2}.
\]
For $\boldsymbol{y}\notin\Sigma$ equation (\ref{eq:SwCond}) still
defines $\lambda$ and equation (\ref{eq:LowDimFil}) is identical
to (\ref{eq:LowDim}) for $\lambda=\pm1$. On $\Sigma$ and when (\ref{eq:Review-SlidingCondition})
holds we calculate $\lambda$ from 

\begin{equation}
Dh(\boldsymbol{y})\cdot\left(\boldsymbol{r}(\boldsymbol{y})+\boldsymbol{b}(\boldsymbol{y})\lambda\right)=0,\label{eq:Review-Fil-Invariance}
\end{equation}
which stipulates that the vector field (\ref{eq:LowDimFil}) is tangential
to $\Sigma$. The solution of (\ref{eq:Review-Fil-Invariance}) is
\begin{equation}
\lambda=-\frac{Dh(\boldsymbol{y})\cdot\boldsymbol{r}(\boldsymbol{y})}{Dh(\boldsymbol{y})\cdot\boldsymbol{b}(\boldsymbol{y})}.\label{eq:FilDef}
\end{equation}

\begin{definition}
Assume that (\ref{eq:Review-SlidingCondition}) holds. We call the
vector field (\ref{eq:LowDimFil}), where $\lambda$ is given by (\ref{eq:FilDef}),
Fillipov's closure.
\end{definition}
It can be shown that $\lambda\in(-1,1)$ when condition (\ref{eq:Review-SlidingCondition})
holds \cite{FilippovBook}. Fillipov's closure is illustrated in figure
\ref{fig:FilUtkin}(c), which shows that the vector field given by
(\ref{eq:LowDimFil}) and (\ref{eq:FilDef}) is chosen from all convex
combinations of $\boldsymbol{f}(\boldsymbol{y},\pm1)$ so that $\boldsymbol{r}(\boldsymbol{y})+\boldsymbol{b}(\boldsymbol{y})\lambda$
is tangential to $\Sigma$. A trajectory at its first point of contact
with $\Sigma$ is continuous, but not continuously differentiable,
because $\lambda$ becomes discontinuous due to (\ref{eq:FilDef}).

When neither (\ref{eq:Review-CrossingCondition}) nor (\ref{eq:Review-SlidingCondition})
holds for $\boldsymbol{y}^{\star}\in\Sigma$, we have 
\begin{equation}
\left(Dh(\boldsymbol{y}^{\star})\boldsymbol{f}(\boldsymbol{y}^{\star},1)\right)\left(Dh(\boldsymbol{y}^{\star})\boldsymbol{f}(\boldsymbol{y}^{\star},-1)\right)=0.\label{eq:Review-TangencyCondition}
\end{equation}
Equation (\ref{eq:Review-TangencyCondition}) means that one or both
of the vector fields $\boldsymbol{f}(\boldsymbol{y},\pm1)$ is tangential
to $\Sigma$, which we call a \emph{tangency}. The boundaries of crossing,
sliding and escaping regions are formed by tangencies, which generally
occur as codimension-one surfaces of $\Sigma$. 

Trajectories may not have unique continuation when they are tangent
to $\Sigma$. When both vector fields $\boldsymbol{f}(\boldsymbol{y},\pm1)$
are tangential to $\Sigma$ at a point, $\lambda$ is not uniquely
defined by (\ref{eq:FilDef}) as both the numerator and denominator
of (\ref{eq:FilDef}) vanish. Consequently, the forward-time solution
of (\ref{eq:LowDim}), (\ref{eq:SwCond}) and (\ref{eq:FilDef}) is
not unique. This case is illustrated in figure \ref{fig:FilUtkin}(d),
which shows a set of possible directions that a solution can follow.
A particular case of this double tangency is the Teixeira singularity,
where an open set of initial conditions generate trajectories that
go through the double tangency. The Teixeira singularity \cite{Teixeira1982}
was studied extensively \cite{ColomboJeffrey,FilippovBook,KristiansenHogan2015,SzalaiNonDet}
in various contexts.

\subsection{\label{subsec:Utkins-closure}Utkin's closure}

In this section we assume that the domain of definition of $\boldsymbol{f}$
is given by (\ref{eq:Review-IntervalDomain}). In case of (\ref{eq:Review-SlidingCondition})
we similarly construct the vector field on $\Sigma$, such that $\Sigma$
becomes invariant under the vector field. The invariance of $\Sigma$
is expressed as
\begin{equation}
Dh(\boldsymbol{y})\cdot\boldsymbol{f}(\boldsymbol{y},\lambda)=0.\label{eq:UtkinInvCond}
\end{equation}
Equation (\ref{eq:UtkinInvCond}) has at least one solution for some
$\lambda\in[-1,1]$, because (\ref{eq:Review-IntervalDomain}) implies
that $Dh(\boldsymbol{y})\cdot\boldsymbol{f}(\boldsymbol{y},\pm1)$
has different signs and due to Bolzano's theorem there must be a root. 
\begin{definition}
Assume that (\ref{eq:Review-SlidingCondition}) holds. We call the
vector field (\ref{eq:LowDimFil}), where $\lambda$ is given by the
solutions of equation (\ref{eq:UtkinInvCond}), Utkin's closure.
\end{definition}
The root of (\ref{eq:UtkinInvCond}) may not be unique, which renders
the solution of (\ref{eq:LowDim}) and (\ref{eq:SwCond}) non-unique.
We also note that (\ref{eq:UtkinInvCond}) can have a solution even
when (\ref{eq:Review-CrossingCondition}) holds in the crossing region.

A simple case of Utkin's closure is illustrated in figure \ref{fig:FilUtkin}(a).
The green curve connecting $\boldsymbol{f}(\boldsymbol{y},1)$ to
$\boldsymbol{f}(\boldsymbol{y},-1)$ represents the possible values
of the vector field on $\Sigma$. There is one intersection of this
family of vectors with the tangent plane of $\Sigma$, represented
by the thick red arrow, which satisfies equation (\ref{eq:UtkinInvCond}).
Figure \ref{fig:FilUtkin}(b) shows that there can be multiple intersections
of $\boldsymbol{f}(\boldsymbol{y},[-1,1])$ with the tangent plane
of $\Sigma$, that then yields multiple solutions. Note that in the
case of figure \ref{fig:FilUtkin}(b), Filippov's closure yields a
unique solution. The contrary, when Utkin's closure predicts a unique
solution and Filippov's closure predicts a family of solutions, is
also possible. For example, when the convex hull represented by the
green dashed line in figure \ref{fig:FilUtkin}(d) is deformed slightly,
the possible number of solutions can be reduced to three. Out of these
three solutions there is only one with $\lambda\in(-1,1)$.

\section{\label{sec:Continuum}Model reduction}

We start with a general continuum model, in the form of

\begin{equation}
\left.\begin{array}{rl}
\dot{\boldsymbol{x}} & =\boldsymbol{F}(\boldsymbol{x},\lambda)\\
\lambda & =\begin{cases}
1 & \mbox{for}\;h(\boldsymbol{x})>0,\\
-1 & \mbox{for}\;h(\boldsymbol{x})<0,
\end{cases}
\end{array}\right\} ,\label{eq:GenEq}
\end{equation}
where $\boldsymbol{x}$ is a function of time $t\in[s,\infty)$, that
is $\boldsymbol{x}:[s,\infty)\to\boldsymbol{X}$ with an initial condition
$\boldsymbol{x}(s)=\boldsymbol{x}_{0}$ and $\boldsymbol{X}$ is an
appropriately defined Banach space. The domain of definition of $\boldsymbol{F}(\cdot,\lambda)$,
for a fixed $\lambda$, is denoted by $\boldsymbol{\mathcal{D}}_{\lambda}(\boldsymbol{F})\subset\boldsymbol{X}$
so that the full domain of definition is $\boldsymbol{\mathcal{D}}(\boldsymbol{F})=\{(\boldsymbol{\mathcal{D}}_{\lambda}(\boldsymbol{F}),\lambda):\lambda\in[-1,1]\}$
and $\boldsymbol{F}:\boldsymbol{\mathcal{D}}(\boldsymbol{F})\to\boldsymbol{X}$.
The switching function $h$ is defined on $\overline{\cup_{\lambda}\boldsymbol{\mathcal{D}}_{\lambda}(\boldsymbol{F})}$
and has values in $\mathbb{R}$. When $h=0$, the most general definition
of the dynamics is $\dot{\boldsymbol{x}}\in\overline{\mathrm{co}}\boldsymbol{F}(\boldsymbol{x},[-1,1])$.
We also require that trajectories are continuous, even when $h=0$.
The smoothness of $\boldsymbol{F}$ and $h$ is not assumed globally,
instead we assume the smoothness of an invariant manifold of $\boldsymbol{F}$
and related quantities in the next section.
\begin{remark}
The notation of equation (\ref{eq:GenEq}) facilitates that $\lambda$
is an unknown, which needs to be found when $h(\boldsymbol{x})=0$.
Therefore $\lambda$ may not be a function of $\boldsymbol{x}$, but
it may become part of the phase space. The solution for $\lambda$,
when $h(\boldsymbol{x})=0$ is defined in section \ref{subsec:NormalDynamics}.
This is a similar setting to section \ref{sec:FilUtkin}, except that
the phase space is now infinite dimensional and therefore a different
kind of solution is required for $\lambda$.
\end{remark}

\subsection{\label{subsec:The-invariant-manifold}The invariant manifold}

For PWS systems, such as equation (\ref{eq:GenEq}), differentiable
invariant manifolds that extend over switching boundaries do not exist.
This fact makes model reduction more complicated than for smooth models.
It is however possible to find invariant manifolds for constant $\lambda$
of the PWS system (\ref{eq:GenEq}). Our approach is therefore to
first consider the smooth system
\begin{equation}
\left.\begin{array}{rl}
\dot{\boldsymbol{x}} & =\boldsymbol{F}(\boldsymbol{x},\lambda)\\
\dot{\lambda} & =0
\end{array}\right\} .\label{eq:VF const lambda}
\end{equation}
We make the following initial assumptions and definitions:
\begin{description}
\item [{\textbf{(A1)\label{A-InvariantManifold}}}] Existence of an invariant
manifold. We assume that there exists a function $\boldsymbol{W}\in C^{p}(\mathbb{R}^{n}\times[-1,1],\boldsymbol{X})$,
$p\ge2$ and a vector field $\boldsymbol{f}\in C^{p}(G\times[-1,1],\mathbb{R}^{n})$,
which satisfies the invariance condition 
\begin{equation}
\boldsymbol{F}(\boldsymbol{W}(\boldsymbol{y},\lambda),\lambda)=D_{1}\boldsymbol{W}(\boldsymbol{y},\lambda)\boldsymbol{f}(\boldsymbol{y},\lambda),\label{eq:ASS-Invariance}
\end{equation}
where $G$ is a compact and connected subset of $\mathbb{R}^{n}$.
The invariant manifold is given by 
\begin{equation}
\mathcal{M}_{\lambda}=\left\{ \boldsymbol{W}(\boldsymbol{y},\lambda):\boldsymbol{y}\in G,\lambda\in[-1,1]\right\} \label{eq:AS-manifold-def}
\end{equation}
and the dynamics of (\ref{eq:VF const lambda}) on $\mathcal{M}_{\lambda}$
is described by 
\begin{equation}
\dot{\boldsymbol{y}}=\boldsymbol{f}(\boldsymbol{y},\lambda).\label{eq:ASS-VectorFieldOnManifold}
\end{equation}
$\boldsymbol{W}$ is called the immersion of $\mathcal{M}_{\lambda}$.
\item [{\textbf{(A2)\label{A-DomainOfDefinition}}}] We assume that for
every $\lambda\in[-1,1]$, $\boldsymbol{F}(\cdot,\lambda)$ is Frechet
differentiable on $\mathcal{M}_{\lambda}$. This derivative is denoted
by
\[
\boldsymbol{A}_{1}(\boldsymbol{y},\lambda)=D_{1}\boldsymbol{F}(\boldsymbol{W}(\boldsymbol{y},\lambda),\lambda).
\]
We also assume that the domain of definition of $\boldsymbol{A}_{1}$,
i.e., $\boldsymbol{\mathcal{D}}(\boldsymbol{A}_{1}(\boldsymbol{y},\lambda))=\left\{ \boldsymbol{x}\in\boldsymbol{X}:\boldsymbol{A}_{1}(\boldsymbol{y},\lambda)\boldsymbol{x}\in\boldsymbol{X}\right\} $,
is independent of $\boldsymbol{y}$ and $\lambda$, and we define
$\boldsymbol{Z}=\overline{\boldsymbol{\mathcal{D}}(\boldsymbol{A}_{1}(\boldsymbol{y},\lambda))}$.
(In general, $\boldsymbol{\mathcal{D}}(\boldsymbol{A}_{1}(\boldsymbol{y},\lambda))\neq\boldsymbol{\mathcal{D}}_{\lambda}(\boldsymbol{F})$.)
\item [{\textbf{(A3)\label{A-EvolutionOperator}}}] Unique continuous solutions.
We assume that the abstract Cauchy problem
\begin{equation}
\left.\begin{array}{rl}
\dot{\boldsymbol{y}} & =\boldsymbol{f}(\boldsymbol{y},\lambda)\\
\dot{\boldsymbol{z}} & =\boldsymbol{A}_{1}(\boldsymbol{y},\lambda)\boldsymbol{z}-D_{2}\boldsymbol{W}(\boldsymbol{y},\lambda)\dot{\lambda}
\end{array}\right\} \label{eq:ASS-VariationalEquation}
\end{equation}
with initial conditions $\boldsymbol{y}(s)\in G$, $\boldsymbol{z}(s)\in\boldsymbol{Z}$,
$s\in\mathbb{R}$ and with $\lambda\in C^{1}([s,\infty),\mathbb{R})$
has a unique solution $\left(\boldsymbol{y},\boldsymbol{z}\right)\in C([s,\infty),G\times\boldsymbol{Z})$,
even though we only have $D_{2}\boldsymbol{W}(\boldsymbol{y}(t),\lambda(t))\in\boldsymbol{X}$.
We also assume that the $\boldsymbol{Z}$ component of the solution
can be written as 
\begin{equation}
\boldsymbol{z}(t)=\boldsymbol{U}(t,s)\boldsymbol{z}(s)-\int_{s}^{t}\boldsymbol{K}(t,\tau)\dot{\lambda}(\tau)\mathrm{d}\tau,\label{eq:ASS-ConvolutionRepresentation}
\end{equation}
where $\boldsymbol{K}$ is bounded for $\tau\ge s$ and continuous
in both variables for $\tau>s$. The underlying conditions of existence
of unique solutions can be found in \cite{DaPrato1992}. For discussion
see remarks \ref{rem:ASS-DefOfSolution} and \ref{rem:ASS-evolutionOp}.
\item [{\textbf{(A4)\label{A-NormalHyperbolicity}}}] $\mathcal{M}_{\lambda}$
is attracting and normally hyperbolic. We assume that there exist
two families of projections $\Pi^{c}(\boldsymbol{y},\lambda)$ and
$\Pi^{s}(\boldsymbol{y},\lambda)$, strongly continuous in $\boldsymbol{y}$
and $\lambda$ such that
\begin{align}
\Pi^{c}(\boldsymbol{y},\lambda)+\Pi^{s}(\boldsymbol{y},\lambda) & =\boldsymbol{I},\nonumber \\
\Pi^{c}(\boldsymbol{y},\lambda)D_{1}\boldsymbol{W}(\boldsymbol{y},\lambda) & =D_{1}\boldsymbol{W}(\boldsymbol{y},\lambda),\label{eq:ASS-TangentBundle}\\
\boldsymbol{U}(t,s)\Pi^{s}(\boldsymbol{y}(s),\lambda(s))\boldsymbol{z} & =\Pi^{s}(\boldsymbol{y}(t),\lambda(t))\boldsymbol{U}(t,s)\boldsymbol{z},\;\forall\boldsymbol{z}\in\boldsymbol{Z},\;t\ge s.\label{eq:ASS-StableBundle}
\end{align}
Consider the non-autonomous ordinary differential equation $\dot{\boldsymbol{\eta}}=D_{1}\boldsymbol{f}(\boldsymbol{y},\lambda)\boldsymbol{\eta}$,
whose solutions with initial condition $\boldsymbol{\eta}_{0}$ at
$t=s$ are denoted by $\boldsymbol{\eta}(t,s,\boldsymbol{\eta}_{0})$.
We assume that there exist real numbers $\sigma_{s}<-\sigma_{c}$,
$M_{c}>0$ and $M_{s}>0$ such that
\begin{align*}
\forall(t-s)\in\mathbb{R},\,\boldsymbol{\eta}_{0}\in\mathbb{R}^{n} & :\left\Vert \boldsymbol{\eta}(t,s,\boldsymbol{\eta}_{0})\right\Vert \le M_{c}\left\Vert \boldsymbol{\eta}_{0}\right\Vert \mathrm{e}^{\sigma_{c}\left|t-s\right|},\\
\forall s\le t,\,\Pi^{s}(\boldsymbol{y}(s),\lambda(s))\boldsymbol{z}=\boldsymbol{z} & :\left\Vert \boldsymbol{U}(t,s)\boldsymbol{z}\right\Vert \le M_{s}\left\Vert \boldsymbol{z}\right\Vert \mathrm{e}^{\sigma_{s}(t-s)}.
\end{align*}
\item [{\textbf{(A5)\label{A-PertTrajectories}}}] We assume that for $t\ge s$
there exists $0<M<\infty$ and $\sigma<0$ such that
\begin{equation}
\left\Vert \boldsymbol{K}(t,s)\right\Vert \le M\mathrm{e}^{\sigma(t-s)}.\label{eq:MR-dichotomy}
\end{equation}
\end{description}
\begin{remark}
For systems with an equilibrium it is natural to consider spectral
submanifolds \cite{Haller2016}, that are the smoothest invariant
manifolds tangent to an invariant linear subspace of the variational
problem about the equilibrium. The uniqueness and existence of such
manifolds is established in \cite{CabreLlave2003}. In order to be
meaningful, these manifolds need to contain the slowest dynamics within
the system to capture long-time behavior. This requirement is outlined
in points R1 and R2 of \cite{HallerExact2016}.
\end{remark}
\begin{remark}
\label{rem:ASS-DefOfSolution}We do not fully specify the definition
of a solution of (\ref{eq:GenEq}) and (\ref{eq:ASS-VariationalEquation})
apart from the solution being continuous. The results of this paper
only depend on the form of the solution as given by (\ref{eq:ASS-ConvolutionRepresentation})
and not how it is obtained. However it might be helpful to think of
\emph{F-solutions} of (\ref{eq:ASS-VariationalEquation}) as defined
by \cite{DaPrato1992}: $\boldsymbol{z}$$\left(\boldsymbol{y},\boldsymbol{z}\right)\in C([s,\infty),G\times\boldsymbol{Z})$
is an F-solution of (\ref{eq:ASS-VariationalEquation}) if there exists
a sequence $\boldsymbol{z}_{k}\in C^{1}([s,\infty),\boldsymbol{Z})\cap C^{1}([s,\infty),\boldsymbol{Z})$
such that
\[
\lim_{k\to\infty}\left[\left\Vert \boldsymbol{z}_{k}(s)-\boldsymbol{z}_{0}\right\Vert +\left\Vert \boldsymbol{z}-\boldsymbol{z}_{k}\right\Vert _{\infty}+\left\Vert \dot{\boldsymbol{z}}_{k}-\boldsymbol{A}_{1}(\boldsymbol{y},\lambda)\boldsymbol{z}_{k}+D_{2}\boldsymbol{W}(\boldsymbol{y},\lambda)\dot{\lambda}\right\Vert _{\infty}\right]=0,
\]
where $\left\Vert \boldsymbol{z}\right\Vert _{\infty}=\sup_{t\in[s,\infty)}\left\Vert \boldsymbol{z}(t)\right\Vert $
and $\boldsymbol{y}$ satisfies $\dot{\boldsymbol{y}}=\boldsymbol{f}(\boldsymbol{y},\lambda)$.
\end{remark}
\begin{remark}
\label{rem:ASS-evolutionOp}The existence of unique F-solutions of
equation (\ref{eq:ASS-VariationalEquation}) is established in \cite{DaPrato1992}
in theorem 5.1. However for many examples, e.g., elastodynamics \cite{GurtinSternbergElastoUnique1961,Martins1987407}
and delay equations \cite{Diekmann1995}, existence and uniqueness
results are already known and it is not necessary to check the conditions
listed in \cite{DaPrato1992}. The existence and regularity of a convolution
kernel for non-autonomous problems is not discussed in the literature.
However, the autonomous problem is discussed in \cite{THIEMEA1990416,Thieme2008},
which implies that the convolution integral is
\begin{equation}
\int_{s}^{t}\boldsymbol{K}(t,\tau)\dot{\lambda}(\tau)\mathrm{d}\tau=\lim_{\mu\to\infty}\int_{s}^{t}\boldsymbol{U}(t,\tau)\mu\left(\mu-\boldsymbol{A}_{1}(\boldsymbol{y}(\tau),\lambda(\tau))\right)^{-1}D_{2}\boldsymbol{W}(\boldsymbol{y}(\tau),\lambda(\tau))\dot{\lambda}(\tau)\mathrm{d}\tau,\label{eq:REM-Convolution}
\end{equation}
because $\left(\mu-\boldsymbol{A}_{1}(\boldsymbol{y}(\tau),\lambda(\tau))\right)^{-1}:\boldsymbol{X}\to\boldsymbol{\mathcal{D}}$.
Reference \cite{DaPrato1992} uses a similar technique to approximate
the unique solution. The kernel $\boldsymbol{K}$, however has two
parameters and therefore its smoothness properties are not trivial
even if we know that the convolution (\ref{eq:REM-Convolution}) is
continuous in $t$. We have therefore assumed the continuity of $\boldsymbol{K}$
for $t>s$, which allows for a discontinuity at $t=s$ due to $D_{2}\boldsymbol{W}(\boldsymbol{y}(\tau),\lambda(\tau))\notin\boldsymbol{Z}$.
\end{remark}
\begin{remark}
\label{rem:ASS-Persistence}The uniqueness or persistence of $\mathcal{M}_{\lambda}$
are not addressed by the assumptions. For persistence of $\mathcal{M}_{\lambda}$
under a perturbation, additional smoothness conditions on the solutions
of (\ref{eq:VF const lambda}) have to hold, which can be found in
\cite{Bates1998}.
\end{remark}
\begin{remark}
The condition (\ref{eq:MR-dichotomy}) implies that the convolution
in (\ref{eq:ASS-ConvolutionRepresentation}) remains bounded when
$\dot{\lambda}$ is bounded. This will be useful later when the reduced
order model is constructed.
\end{remark}
\begin{figure}
\begin{centering}
\includegraphics[width=1\linewidth]{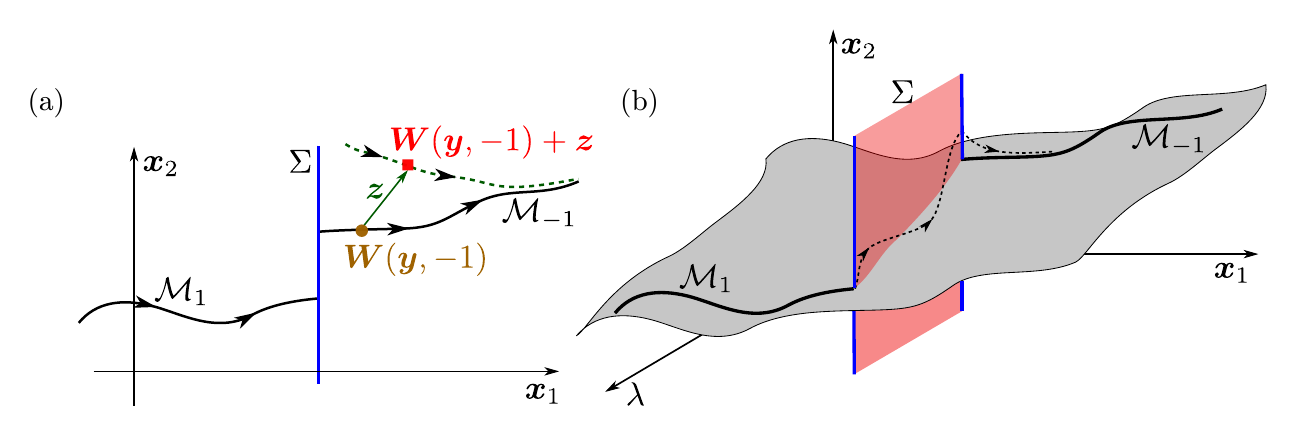}
\par\end{centering}
\caption{\label{fig:GeomSketch}(a) The manifolds $\mathcal{M}_{\pm1}$ are
the images of $\boldsymbol{W}(\cdot,\pm1)$. The two sections $\mathcal{M}_{-1}$
and $\mathcal{M}_{1}$ of the invariant manifold $\mathcal{M}_{\lambda}$
are not joined up along $\Sigma$. This implies that trajectories
restricted to $\mathcal{M}_{\lambda}$ cannot cross $\Sigma$ instantaneously
without being discontinuous. The difference between a trajectory on
$\mathcal{M}_{\lambda}$ and a trajectory about $\mathcal{M}_{\lambda}$
is represented by $\boldsymbol{z}$. The dashed line illustrates that
trajectories near $\mathcal{M}_{\lambda}$ are asymptotic to $\mathcal{M}_{\lambda}$
in forward time\@. (b) The dynamics of (\ref{eq:GenEq}) is considered
in the extended space $\boldsymbol{X}\times[-1,1]$ and in the neighborhood
of $\mathcal{M}_{\lambda}$. For a trajectory that enters $\Sigma$,
the value of $\lambda$ becomes time dependent and for such trajectories
$\mathcal{M}_{\lambda}$ is not invariant. The dotted line shows a
hypothetical trajectory leaving $\mathcal{M}_{\lambda}$.}
\end{figure}

Figure \ref{fig:GeomSketch} shows the invariant manifold $\mathcal{M}_{\lambda}$
and its intersection with the switching manifold $\Sigma$. The Banach
space $\boldsymbol{X}$ is represented by two coordinates $\boldsymbol{x}_{1}$
and $\boldsymbol{x}_{2}$. The two parts of the invariant manifold
($\mathcal{M}_{-1}$ and $\mathcal{M}_{1}$) do not join up in figure
\ref{fig:GeomSketch}(a). If trajectories cross $\Sigma$ instantaneously,
they are discontinuous. When discontinuity is not allowed, the crossing
cannot be instantaneous. In certain cases, however, the disconnected
nature of $\mathcal{M}_{\lambda}$ may be overlooked. For example,
when the switching function solely depends on the parameter $\boldsymbol{y}$
of the immersion $\boldsymbol{W}$, i.e., $\boldsymbol{y}=\boldsymbol{x}_{1}$
in figure \ref{fig:GeomSketch}. The case when the dynamics is restricted
to $\mathcal{M}_{\lambda}$ is discussed in section \ref{sec:skeleton}.

Figure \ref{fig:GeomSketch}(b) shows the extended phase space and
how solutions of (\ref{eq:GenEq}) behave about $\mathcal{M}_{\lambda}$,
when instantaneous crossing is not allowed. In figure \ref{fig:GeomSketch}(b)
$\mathcal{M}_{\lambda}$ is a connected manifold. When a solution
of (\ref{eq:GenEq}) arrives at $\Sigma$, the value of $\lambda$
must change, so that a trajectory can enter $\Sigma$. $\mathcal{M}_{\lambda}$
is only invariant for constant $\lambda$ and therefore a trajectory
(denoted by dotted lines) will not continue on $\mathcal{M}_{\lambda}$,
while also in $\Sigma$. Once a trajectory has left $\Sigma$ it will
be attracted to $\mathcal{M}_{\lambda}$ as per assumption \nameref{A-NormalHyperbolicity}.

In the following sections we discuss how the departure of a trajectory
from $\mathcal{M}_{\lambda}$ can be captured and whether or not capturing
this dynamics makes a qualitative difference in the predictions of
the model. In section \ref{sec:TrivialExample} we have already seen
that including a correction that captures the departure from $\mathcal{M}_{\lambda}$
makes a difference and trajectories can no longer cross $\Sigma$
instantaneously.

\subsection{\label{sec:skeleton}The skeleton model}

Having assumed the existence of an invariant manifold $\mathcal{M}_{\lambda}$,
it is natural to consider the dynamics on $\mathcal{M}_{\lambda}$
in the presence of switching. This can be done by substituting the
immersion $\boldsymbol{W}$ into the full problem (\ref{eq:GenEq})
and disregarding that $\lambda$ may not be constant on $\Sigma$.
We start with the switching function 
\begin{equation}
h_{0}(\boldsymbol{y},\lambda)=h(\boldsymbol{W}(\boldsymbol{y},\lambda)).\label{eq:RED-h0}
\end{equation}
In contrast to section \ref{sec:FilUtkin}, the switching function
(\ref{eq:RED-h0}) depends on $\lambda$ and therefore the closures
described in section \ref{sec:FilUtkin} may not apply. Using the
vector field (\ref{eq:ASS-VectorFieldOnManifold}) on $\mathcal{M}_{\lambda}$
and (\ref{eq:RED-h0}) we obtain
\begin{equation}
\left.\begin{array}{rl}
\dot{\boldsymbol{y}} & =\boldsymbol{f}(\boldsymbol{y},\lambda)\\
\lambda & =\begin{cases}
1 & \mbox{for}\;h_{0}(\boldsymbol{y},\lambda)>0\\
-1 & \mbox{for}\;h_{0}(\boldsymbol{y},\lambda)<0
\end{cases}
\end{array}\right\} ,\label{eq:RED-skeleton}
\end{equation}
where $\boldsymbol{y}(t)\in G$ for all $t\in[s,s+\Delta)$, $\Delta>0$. 
\begin{definition}
\label{def:RED-skeleton}Equation (\ref{eq:RED-skeleton}) is called
the \emph{skeleton model} of (\ref{eq:GenEq}) on the invariant manifold
$\mathcal{M}_{\lambda}$.
\end{definition}
Definition \ref{def:RED-skeleton} alludes to what follows next. We
will use equation (\ref{eq:RED-skeleton}) to build upon and not consider
it as an end result. Equation (\ref{eq:RED-skeleton}) is inaccurate
when $\lambda$ varies and that causes solutions to become non-unique,
even if they were unique in the full problem (\ref{eq:GenEq}). Nevertheless,
we highlight some properties of the skeleton model that carry over
to the reduced order model.

We note that already in the skeleton model the switching function
$h_{0}$ can become dependent on $\lambda$. This means that the dynamics
when $h_{0}=0$ may be defined as an index-1 differential algebraic
equation. To describe such dynamics, in the introductory example in
equation (\ref{eq:Trivial Reduced Model}) we needed to separate the
switching manifold into two components. Here, we formalize this splitting
and define two new switching manifolds
\begin{equation}
\Sigma_{0}^{\pm}=\left\{ \boldsymbol{y}\in G:h_{0}(\boldsymbol{y},\pm1)=0\right\} .\label{eq:RED-SwitchManif-PM}
\end{equation}
In the extended state space $\left(\boldsymbol{y},\lambda\right)\in G\times[-1,1]$,
$\Sigma_{0}^{\pm}$ is the boundary of the $n$-dimensional manifold
\begin{equation}
\Sigma_{0}=\left\{ \left(\boldsymbol{y},\lambda\right)\in G\times[-1,1]:h_{0}(\boldsymbol{y},\lambda)=0\right\} .\label{eq:RED-SwitchManif-Full}
\end{equation}
$\Sigma_{0}^{\pm}$ cannot intersect each other in the extended state
space. 

When $h_{0}(\boldsymbol{y},\lambda)\neq0$ the trajectories are described
by the vector field
\begin{equation}
\left.\begin{array}{rl}
\dot{\boldsymbol{y}} & =\boldsymbol{f}(\boldsymbol{y},\lambda)\\
\dot{\lambda} & =0
\end{array}\right\} .\label{eq:RED-SmoothEq}
\end{equation}
Otherwise, we have an index-1 differential algebraic equation 
\begin{equation}
\left.\begin{array}{rl}
\dot{\boldsymbol{y}} & =\boldsymbol{f}(\boldsymbol{y},\lambda)\\
0 & =h_{0}(\boldsymbol{y},\lambda)
\end{array}\right\} .\label{eq:RED-DAE}
\end{equation}
A unique solution to (\ref{eq:RED-DAE}) is guaranteed by the Implicit
Function Theorem if $D_{2}h_{0}(\boldsymbol{y},\lambda)\neq0$, so
that there is a unique $C^{p}$ smooth function $\lambda(\boldsymbol{y})$
satisfying $h_{0}(\boldsymbol{y},\lambda(\boldsymbol{y}))=0$. This
also implies that $\lambda(t)=\lambda(\boldsymbol{y}(t))$ is continuous,
hence there is no discontinuity of $\lambda$ when a trajectory reaches
$\Sigma_{0}^{\pm}$ transversely. Trajectories must spend nonzero
time on $\Sigma_{0}$ in order to keep $\lambda$ continuous. This
short argument highlights a major difference between PWS models described
in section \ref{sec:FilUtkin}, where we have $D_{2}h(\boldsymbol{y},\lambda)=0$
and where the Implicit Function Theorem does not apply. Models in
section \ref{sec:FilUtkin} are special cases of the skeleton model.

We can also write the index-1 differential algebraic equation (\ref{eq:RED-DAE})
in a differential form by differentiating the constraint $h_{0}(\boldsymbol{y},\lambda)=0$,
that is, 
\begin{equation}
\frac{\mathrm{d}}{\mathrm{d}t}h_{0}(\boldsymbol{y},\lambda)=D_{1}h_{0}(\boldsymbol{y},\lambda)\boldsymbol{f}(\boldsymbol{y},\lambda)+D_{2}h_{0}(\boldsymbol{y},\lambda)\dot{\lambda}=0.\label{eq:RED-hDerivative}
\end{equation}
As discussed, whether solutions are well defined, depends on the term
\begin{equation}
D_{2}h_{0}(\boldsymbol{y},\lambda)=Dh(\boldsymbol{W}(\boldsymbol{y},\lambda))\cdot D_{2}\boldsymbol{W}(\boldsymbol{y},\lambda).\label{eq:RED-dMinus}
\end{equation}
If (\ref{eq:RED-dMinus}) is non-zero, equation (\ref{eq:RED-hDerivative})
can be solved for $\dot{\lambda},$ which yields the differential
form of (\ref{eq:RED-DAE}) for $\left(\boldsymbol{y},\lambda\right)\in\Sigma_{0}$,
that is, 
\begin{equation}
\left.\begin{array}{rl}
\dot{\boldsymbol{y}} & =\boldsymbol{f}(\boldsymbol{y},\lambda)\\
\dot{\lambda} & =-\dfrac{D_{1}h_{0}(\boldsymbol{y},\lambda)\boldsymbol{f}(\boldsymbol{y},\lambda)}{D_{2}h_{0}(\boldsymbol{y},\lambda)}
\end{array}\right\} .\label{eq:RED-SwitchingEq}
\end{equation}
The continuous concatenation of solutions of equations (\ref{eq:RED-SmoothEq})
and (\ref{eq:RED-SwitchingEq}) gives the full solution of equation
(\ref{eq:RED-skeleton}). This concatenation is a PWS problem, where
$\Sigma_{0}^{\pm}$ are now separating the phase space into three
regions. The following theorem looks at the case when there is no
need to define the dynamics on $\Sigma_{0}^{\pm}$.
\begin{theorem}
\label{thm:ExtUnique}Consider a point $(\boldsymbol{y}^{\star},\lambda^{\star})\in\Sigma_{0}^{\pm}$
and assume that 
\begin{equation}
D_{2}h_{0}(\boldsymbol{y}^{\star},\lambda^{\star})<0.\label{eq:UniquenessCondition}
\end{equation}
Further assume a solution $(\boldsymbol{y}(t),\lambda(t))$, for $t\in I=(-\delta,0]$
(or $t\in I=[0,\delta)$) of either equation (\ref{eq:RED-SmoothEq})
or equation (\ref{eq:RED-SwitchingEq}) that reaches $(\boldsymbol{y}^{\star},\lambda^{\star})$
at $t=0$. The corresponding trajectory is defined as $\mathcal{T}=\left\{ (\boldsymbol{y}(t),\lambda(t)):t\in I\right\} $.
Trajectory $\mathcal{T}$ has a unique continuation for $t>0$ (or
$t<0$) sufficiently small as a solution of the skeleton model (\ref{eq:RED-skeleton})
if one of the following conditions holds:
\begin{enumerate}
\item $\mathcal{T}$ is not tangent to $\Sigma_{0}^{\pm}$, i.e., $D_{1}h_{0}(\boldsymbol{y}^{\star},\lambda^{\star})\boldsymbol{f}(\boldsymbol{y}^{\star},\lambda^{\star})\neq0$
\item $\mathcal{T}$ is tangent to $\Sigma_{0}^{\pm}$ and the order of
the tangency is less than the smoothness order ($C^{p}$) of $h_{0}$.
In other words, there exists $0<\ell\le p$ such that
\begin{equation}
\frac{\mathrm{d}^{\ell}}{\mathrm{d}t^{\ell}}h_{0}(\boldsymbol{y}(t),\lambda^{\star})\vert_{t=0}\neq0.\label{eq:ExtUni-TangencyCondition}
\end{equation}
\end{enumerate}
\end{theorem}
\begin{svmultproof}
The proof can be found in appendix \ref{sec:Proof-Skeleton-Unique}.
\end{svmultproof}

\begin{remark}
Theorem \ref{thm:ExtUnique} excludes the case $D_{2}h_{0}(\boldsymbol{y},\lambda)>0$.
For $D_{2}h_{0}(\boldsymbol{y},\lambda)>0$, transverse trajectories
(case 1 of theorem \ref{thm:ExtUnique}) cannot cross $\Sigma_{0}^{\pm}$.
Tangential trajectories with even $\ell$ may have multiple continuation,
which is the case of the Teixeira singularity \cite{ColomboJeffrey}.
Tangential trajectories with odd $\ell$ cannot cross $\Sigma_{0}^{\pm}$,
similar to transverse trajectories. To investigate the case of $D_{2}h_{0}(\boldsymbol{y},\lambda)>0$
in detail, a definition of how trajectories move along $\Sigma_{0}^{\pm}$
(with $\lambda=\pm1$) is also required, which falls outside of the
scope of this paper.
\end{remark}

\subsection{\label{subsec:NormalDynamics}Dynamics about manifold $\mathcal{M}_{\lambda}$
due to switching}

This section describes a correction to the skeleton model (\ref{eq:RED-skeleton})
that resolves the dynamics in the neighborhood of $\mathcal{M}_{\lambda}$
up to linear order. The correction is necessary, because the \textbf{$\dot{\lambda}=0$}
assumption does not hold: equation (\ref{eq:RED-SwitchingEq}) states
that $\lambda$ varies on $\Sigma_{0}$. The correction that is introduced
here captures trajectories that depart from $\mathcal{M}_{\lambda}$
when $h=0$ (see dashed line in figure \ref{fig:GeomSketch}(b)).

Let us suppose that 
\begin{align}
\boldsymbol{x} & =\boldsymbol{W}(\boldsymbol{y},\lambda)+\boldsymbol{z},\label{eq:zDefinition}
\end{align}
where $\boldsymbol{z}$ represents the difference between the trajectories
of the full model (\ref{eq:GenEq}) and the skeleton model (\ref{eq:RED-skeleton}).
This set-up is illustrated in figure \ref{fig:GeomSketch}(a). To
derive an equation for $\boldsymbol{z}$, we substitute (\ref{eq:zDefinition})
into (\ref{eq:GenEq}) while taking into account that $\lambda$ is
a function of time. This substitution yields
\begin{equation}
\dot{\boldsymbol{x}}=D_{1}\boldsymbol{W}(\boldsymbol{y},\lambda)\dot{\boldsymbol{y}}+D_{2}\boldsymbol{W}(\boldsymbol{y},\lambda)\dot{\lambda}+\dot{\boldsymbol{z}}=\boldsymbol{F}(\boldsymbol{W}(\boldsymbol{y},\lambda)+\boldsymbol{z},\lambda).\label{eq:zSubtitution}
\end{equation}
We assume that $\boldsymbol{z}$ is a small deviation from $\mathcal{M}_{\lambda}$
and Taylor expand $\boldsymbol{F}(\boldsymbol{W}(\boldsymbol{y},\lambda)+\boldsymbol{z},\lambda)$
in $\boldsymbol{z}$ about $\boldsymbol{z}=\boldsymbol{0}$, that
is,
\begin{equation}
\boldsymbol{F}(\boldsymbol{W}(\boldsymbol{y},\lambda)+\boldsymbol{z},\lambda)=\boldsymbol{F}(\boldsymbol{W}(\boldsymbol{y},\lambda),\lambda)+D_{1}\boldsymbol{F}(\boldsymbol{W}(\boldsymbol{y},\lambda),\lambda)\boldsymbol{z}+\mathcal{O}\left(\left\Vert \boldsymbol{z}\right\Vert ^{2}\right).\label{eq:zLinearExpansion}
\end{equation}
The expansion (\ref{eq:zLinearExpansion}), when substituted into
(\ref{eq:zSubtitution}) yields
\begin{equation}
D_{1}\boldsymbol{W}(\boldsymbol{y},\lambda)\dot{\boldsymbol{y}}+D_{2}\boldsymbol{W}(\boldsymbol{y},\lambda)\dot{\lambda}+\dot{\boldsymbol{z}}=\boldsymbol{F}(\boldsymbol{W}(\boldsymbol{y},\lambda),\lambda)+D_{1}\boldsymbol{F}(\boldsymbol{W}(\boldsymbol{y},\lambda),\lambda)\boldsymbol{z}+\mathcal{O}\left(\left\Vert \boldsymbol{z}\right\Vert ^{2}\right).\label{eq:zLinearSubstitution}
\end{equation}
We now use the invariance equation (\ref{eq:ASS-Invariance}) and
the dynamics on $\mathcal{M}_{\lambda}$ as given by (\ref{eq:ASS-VectorFieldOnManifold})
and notice that two terms cancel in (\ref{eq:zLinearSubstitution}),
so that we get
\begin{equation}
D_{2}\boldsymbol{W}(\boldsymbol{y},\lambda)\dot{\lambda}+\dot{\boldsymbol{z}}=D_{1}\boldsymbol{F}(\boldsymbol{W}(\boldsymbol{y},\lambda),\lambda)\boldsymbol{z}+\mathcal{O}\left(\left\Vert \boldsymbol{z}\right\Vert ^{2}\right).\label{eq:zFinalLinearExp}
\end{equation}
Combining the skeleton model (\ref{eq:RED-skeleton}) with (\ref{eq:zFinalLinearExp})
yields the corrected model

\begin{equation}
\left.\begin{array}{rl}
\dot{\boldsymbol{y}} & =\boldsymbol{f}(\boldsymbol{y},\lambda)\\
\dot{\boldsymbol{z}} & =\boldsymbol{A}_{1}(\boldsymbol{y},\lambda)\boldsymbol{z}-D_{2}\boldsymbol{W}(\boldsymbol{y},\lambda)\dot{\lambda}\\
\lambda & =\begin{cases}
1 & \mbox{for}\;h(\boldsymbol{W}(\boldsymbol{y},\lambda)+\boldsymbol{z})>0\\
-1 & \mbox{for}\;h(\boldsymbol{W}(\boldsymbol{y},\lambda)+\boldsymbol{z})<0
\end{cases}
\end{array}\right\} ,\label{eq:zLinearEq}
\end{equation}
where $\boldsymbol{A}_{1}(\boldsymbol{y},\lambda)=D_{1}\boldsymbol{F}(\boldsymbol{W}(\boldsymbol{y},\lambda),\lambda)$
is defined in assumption \nameref{A-DomainOfDefinition}. A unique
solution of (\ref{eq:zLinearEq}) is assumed in \nameref{A-EvolutionOperator}
with a continuously differentiable $\lambda$. In this paper we do
not investigate whether the corrected model (\ref{eq:zLinearEq})
is a faithful representation of the fully nonlinear system (\ref{eq:GenEq});
for some discussion, see remark \ref{rem:linear-correction}. 

We define the switching manifolds as
\[
\Sigma=\left\{ \left(\boldsymbol{y},\boldsymbol{z},\lambda\right)\in G\times\boldsymbol{Z}\times[-1,1]:h(\boldsymbol{W}(\boldsymbol{y},\lambda)+\boldsymbol{z})=0\right\} 
\]
and
\[
\Sigma^{\pm}=\left\{ \left(\boldsymbol{y},\boldsymbol{z},\pm1\right)\in G\times\boldsymbol{Z}\times[-1,1]:h(\boldsymbol{W}(\boldsymbol{y},\pm1)+\boldsymbol{z})=0\right\} .
\]
When a trajectory is restricted to $\Sigma$, the solution must satisfy
\begin{equation}
h(\boldsymbol{W}(\boldsymbol{y},\lambda)+\boldsymbol{z})=0.\label{eq:INF-constraint}
\end{equation}
Similar to the skeleton model we evaluate how $h$ changes in time
and restrict this change to zero on $\Sigma$ to find an equation
for $\lambda$ (cf. equation (\ref{eq:RED-hDerivative})). To evaluate
equation (\ref{eq:INF-constraint}), we use the following lemma.
\begin{lemma}
\label{lem:INF-discontinuity-gap}Assume \nameref{A-EvolutionOperator}
and that $\lambda$ is continuously differentiable and $\boldsymbol{y}$,
$\boldsymbol{z}$ satisfy the differential equations
\begin{align*}
\dot{\boldsymbol{y}} & =\boldsymbol{f}(\boldsymbol{y},\lambda)\\
\dot{\boldsymbol{z}} & =\boldsymbol{A}_{1}(\boldsymbol{y},\lambda)\boldsymbol{z}-D_{2}\boldsymbol{W}(\boldsymbol{y},\lambda)\dot{\lambda}
\end{align*}
on the interval $t\in[s,s+\epsilon)$, $\epsilon>0$ with an initial
condition $\boldsymbol{y}(s)\in G$, $\boldsymbol{z}(s)\in\boldsymbol{\mathcal{D}}$.
Then the right-side derivative of $h$ as a function of time is calculated
as
\begin{equation}
\frac{\mathrm{d}}{\mathrm{d}t^{+}}h(\boldsymbol{W}(\boldsymbol{y},\lambda)+\boldsymbol{z})=Dh(\boldsymbol{W}(\boldsymbol{y},\lambda)+\boldsymbol{z})\cdot D_{1}\boldsymbol{W}(\boldsymbol{y},\lambda)\boldsymbol{f}(\boldsymbol{y},\lambda)-d^{\pm}(\boldsymbol{y},\boldsymbol{z},\lambda)\dot{\lambda}+\boldsymbol{A}_{1}(\boldsymbol{y},\lambda)\boldsymbol{z},\label{eq:INF-hDerivative}
\end{equation}
where
\begin{equation}
d^{\pm}(\boldsymbol{y},\boldsymbol{z},\lambda)=\lim_{\delta\downarrow0}Dh(\boldsymbol{W}(\boldsymbol{y},\lambda)+\boldsymbol{z})\cdot\left(\boldsymbol{K}(t+\delta,t)-D_{2}\boldsymbol{W}(\boldsymbol{y},\lambda)\right).\label{eq:dPlusMinusDefinition}
\end{equation}
\end{lemma}
\begin{svmultproof}
The proof can be found in appendix \ref{sec:UniquenessProof}.
\end{svmultproof}

\begin{remark}
The quantity $d^{\pm}(\boldsymbol{y},\boldsymbol{z},\lambda)$ in
(\ref{eq:dPlusMinusDefinition}) measures the discontinuity of the
convolution kernel $\boldsymbol{K}$ at $t=s$. A discontinuous $\boldsymbol{K}$
is possible, because $D_{2}\boldsymbol{W}(\boldsymbol{y},\lambda)\in\boldsymbol{X}\backslash\boldsymbol{Z}$,
and the continuity assumption \nameref{A-EvolutionOperator} does
not apply at $t=s$. Such a discontinuity allowed us to find a differential
equation for $\lambda$ in section \ref{sec:TrivialExample}.
\end{remark}
\begin{definition}
\label{def:INF-NormalDiscontinuityGap}We call the quantity $d^{\pm}(\boldsymbol{y},\boldsymbol{z},\lambda)$
in equation (\ref{eq:dPlusMinusDefinition}) the \emph{normal discontinuity
gap}.
\end{definition}
We also define two other quantities that will be useful later. These
are
\begin{align}
d^{-}(\boldsymbol{y},\boldsymbol{z},\lambda) & =Dh(\boldsymbol{W}(\boldsymbol{y},\lambda)+\boldsymbol{z})\cdot D_{2}\boldsymbol{W}(\boldsymbol{y},\lambda),\label{eq:INF-dMinus}\\
d^{+}(\boldsymbol{y},\boldsymbol{z},\lambda) & =\lim_{\delta\downarrow0}Dh(\boldsymbol{W}(\boldsymbol{y},\lambda)+\boldsymbol{z})\cdot\boldsymbol{K}(t+\delta,t),\label{eq:INF-dPlus}
\end{align}
and therefore we have the identity $d^{\pm}(\boldsymbol{y},\boldsymbol{z},\lambda)=d^{+}(\boldsymbol{y},\boldsymbol{z},\lambda)-d^{-}(\boldsymbol{y},\boldsymbol{z},\lambda)$. 

We now find the governing equation of the dynamics on $\Sigma$. We
solve equation $\frac{\mathrm{d}}{\mathrm{d}t^{+}}h=0$, where $\frac{\mathrm{d}}{\mathrm{d}t^{+}}h$
is given by (\ref{eq:INF-hDerivative}) for $\dot{\lambda}$, which
yields
\begin{equation}
\left.\begin{array}{rl}
\dot{\boldsymbol{y}} & =\boldsymbol{f}(\boldsymbol{y},\lambda)\\
\dot{\boldsymbol{z}} & =\boldsymbol{A}_{1}(\boldsymbol{y},\lambda)\boldsymbol{z}-D_{2}\boldsymbol{W}(\boldsymbol{y},\lambda)\dot{\lambda}\\
\dot{\lambda} & ={\displaystyle \frac{1}{d^{\pm}(\boldsymbol{y},\boldsymbol{z},\lambda)}Dh(\boldsymbol{W}(\boldsymbol{y},\lambda)+\boldsymbol{z})\cdot\left(D_{1}\boldsymbol{W}(\boldsymbol{y},\lambda)\boldsymbol{f}(\boldsymbol{y},\lambda)+\boldsymbol{A}_{1}(\boldsymbol{y},\lambda)\boldsymbol{z}\right)}
\end{array}\right\} .\label{eq:INF-SlidingDynamics}
\end{equation}
The trajectories of equation (\ref{eq:INF-SlidingDynamics}) are concatenated
with trajectories of
\begin{equation}
\left.\begin{array}{rl}
\dot{\boldsymbol{y}} & =\boldsymbol{f}(\boldsymbol{y},\lambda)\\
\dot{\boldsymbol{z}} & =\boldsymbol{A}_{1}(\boldsymbol{y},\lambda)\boldsymbol{z}\\
\dot{\lambda} & =0
\end{array}\right\} \label{eq:INF-NormalDynamics}
\end{equation}
along the boundaries $\Sigma^{\pm}$ and form the trajectories of
the corrected model (\ref{eq:zLinearEq}). The following theorem provides
a sufficient condition for a unique continuation of trajectories through
$\Sigma^{\pm}$.
\begin{theorem}
\label{thm:Cont-Uniqueness}Assume \nameref{A-InvariantManifold}-\nameref{A-PertTrajectories}.
A trajectory $\mathcal{T}$ of either (\ref{eq:INF-SlidingDynamics})
or (\ref{eq:INF-NormalDynamics}) with an end point $\left(\boldsymbol{y},\boldsymbol{z},\lambda\right)\in\Sigma^{\pm}$
at $t=s$ has a unique continuation for $t>s$ with $t-s$ sufficiently
small, as a solution of the corrected model (\ref{eq:zLinearEq}),
if the following conditions hold:
\begin{enumerate}
\item 
\begin{equation}
d^{\pm}(\boldsymbol{y},\boldsymbol{z},\lambda)>0,\label{eq:ContUniquenessCondition}
\end{equation}
\item $Dh(\boldsymbol{W}(\boldsymbol{y},\lambda)+\boldsymbol{z})\cdot\boldsymbol{U}(t,s)\boldsymbol{z}$
is continuously differentiable with respect to $t$ for $t\ge s$
and
\item one of the vector fields, (\ref{eq:INF-SlidingDynamics}) or (\ref{eq:INF-NormalDynamics})
is not tangent to $\Sigma^{\pm}$, that is, 
\begin{equation}
Dh(\boldsymbol{W}(\boldsymbol{y},\lambda)+\boldsymbol{z})\cdot\left(D_{1}\boldsymbol{W}(\boldsymbol{y},\lambda)\boldsymbol{f}(\boldsymbol{y},\lambda)+\boldsymbol{A}_{1}(\boldsymbol{y},\lambda)\boldsymbol{z}\right)\neq0.\label{eq:INF-tangency}
\end{equation}
\end{enumerate}
\end{theorem}
\begin{svmultproof}
The proof of theorem \ref{thm:Cont-Uniqueness} can be found in appendix
\ref{sec:UniquenessProof}.
\end{svmultproof}

\begin{remark}
\label{rem:linear-correction}The linear correction about the invariant
manifold is carried out here without an assessment whether trajectories
of the corrected model (\ref{eq:zLinearEq}) and the full model (\ref{eq:GenEq})
are qualitatively the same. If $\left\Vert \boldsymbol{z}\right\Vert \ll1$
the linear correction is accurate. Because on $\text{\ensuremath{\mathcal{M}}}_{\lambda}$
we have $\boldsymbol{z}=\boldsymbol{0}$, when a trajectory enters
$\Sigma$, the rate of change of $\boldsymbol{z}$ is determined by
$\dot{\lambda}$. The magnitude of $\dot{\lambda}$ depends on the
$\boldsymbol{f}$ and $d^{\pm}$. Smaller $d^{\pm}$ makes $\lambda$
faster. The value of $d^{\pm}$ is not necessarily a small parameter
and therefore the deviation from $\text{\ensuremath{\mathcal{M}}}_{\lambda}$
can stay small. For the linear string $d^{\pm}=\frac{1}{2}$. In the
literature of regularized PWS systems \cite{JeffreyTautology,KristiansenHogan2015},
to stay close to the skeleton model, fast $\lambda$ is assumed.
\end{remark}
\begin{remark}
If $d^{\pm}(\boldsymbol{y},\boldsymbol{z},\lambda)=0$, the dynamics
about $\mathcal{M}_{\lambda}$ as captured by variable $\boldsymbol{z}$
can only have a second order effect on $h$ due to the nonlinearity
of $h$. Therefore (\ref{eq:INF-hDerivative}) is independent of $\dot{\lambda}$
and $\frac{\mathrm{d}}{\mathrm{d}t^{+}}h=0$ cannot be solved for
$\dot{\lambda}$. When $d^{\pm}(\boldsymbol{y},\boldsymbol{z},\lambda)=0$
the corrected model (\ref{eq:zLinearEq}) needs a closure, such as
Filippov's or Utkin's. $d^{\pm}(\boldsymbol{y},\boldsymbol{z},\lambda)=0$
occurs when $\boldsymbol{U}$ is strongly continuous on the whole
of $\boldsymbol{X}$, i.e., $\boldsymbol{Z}=\boldsymbol{X}$. This
case for linear systems is explored in \cite{ORLOV1995473,Levaggi2002508,Levaggi2002167}.
\end{remark}
\begin{remark}
The transversality condition (\ref{eq:INF-tangency}) is the equivalent
of case 1 of theorem \ref{thm:ExtUnique}. The equivalent of case
2 of theorem \ref{thm:ExtUnique} is not proven here, but a similar
argument can be made while carefully accounting for the infinite dimensional
nature of the problem.
\end{remark}
\begin{remark}
It is possible to consider a nonlinear correction, so that (\ref{eq:zDefinition})
becomes exact. Let us define the nonlinear term
\[
\boldsymbol{N}(\boldsymbol{y},\lambda;\boldsymbol{z})=\boldsymbol{F}(\boldsymbol{W}(\boldsymbol{y},\lambda)+\boldsymbol{z},\lambda)-\boldsymbol{F}(\boldsymbol{W}(\boldsymbol{y},\lambda),\lambda)-\boldsymbol{A}_{1}(\boldsymbol{y},\lambda)\boldsymbol{z}
\]
without discussing the constraints on $\boldsymbol{N}$. The equation
of the exact correction can be written as
\begin{equation}
\left.\begin{array}{rl}
\dot{\boldsymbol{y}} & =\boldsymbol{f}(\boldsymbol{y},\lambda)\\
\dot{\boldsymbol{z}} & =\boldsymbol{A}_{1}(\boldsymbol{y},\lambda)\boldsymbol{z}+\boldsymbol{N}(\boldsymbol{y},\lambda;\boldsymbol{z})-D_{2}\boldsymbol{W}(\boldsymbol{y},\lambda)\dot{\lambda}
\end{array}\right\} .\label{eq:INF-ExactCorrection}
\end{equation}
Equation (\ref{eq:INF-ExactCorrection}) is a semi-linear abstract
Cauchy problem, which is frequently analyzed in the mathematical literature.
The solutions of (\ref{eq:INF-ExactCorrection}) are formally obtained
from the integral equation
\begin{equation}
\boldsymbol{z}(t)=\boldsymbol{U}(t,s)\boldsymbol{z}(s)+\int_{s}^{t}\boldsymbol{U}(t,\tau)\left(\boldsymbol{N}(\boldsymbol{y}(\tau),\lambda(\tau);\boldsymbol{z}(\tau))-D_{2}\boldsymbol{W}(\boldsymbol{y}(\tau),\lambda(\tau))\dot{\lambda}(\tau)\right)\mathrm{d}\tau.\label{eq:INF-ExactSolution}
\end{equation}
In general, existence and uniqueness of solutions of (\ref{eq:INF-ExactSolution})
is established using a contraction mapping argument. However under
our assumptions the convolution is not justified because $\boldsymbol{U}(t,s)$
is only defined on $\boldsymbol{Z}$, but 
\begin{equation}
\boldsymbol{N}(\boldsymbol{y},\lambda;\boldsymbol{z})-D_{2}\boldsymbol{W}(\boldsymbol{y},\lambda)\dot{\lambda}\notin\boldsymbol{Z}.\label{eq:INF-NotInDomain}
\end{equation}
Regardless of (\ref{eq:INF-NotInDomain}), the autonomous case \cite{Thieme2008,magal2009}
has unique solutions under appropriate conditions. The author is confident
that a similar argument can be made to establish unique solutions
(\ref{eq:INF-ExactSolution}) although that might require that the
nonlinearity $\boldsymbol{N}(\boldsymbol{y},\lambda;\cdot):\boldsymbol{\mathcal{D}}\to\boldsymbol{X}$
be bounded.
\end{remark}

\subsection{Time-scale separation}

We already have some indication that switching has a great influence
on the normal dynamics. For example, ignoring the normal dynamics
as in the skeleton model (\ref{eq:RED-skeleton}) leads to a different
uniqueness condition than for the corrected model (\ref{eq:zLinearEq}).
In this section we restrict the analysis to the simplest case where
there is a separation of time scales. We assume a parameter $0\le\varepsilon\le1$
and denote the dependence on $\varepsilon$ by a subscript, that is
$\boldsymbol{F}_{\varepsilon}$. Here, the $\varepsilon=0$ limit
is represented by the skeleton model (\ref{eq:RED-skeleton}) and
$\varepsilon=1$ refers to the corrected model (\ref{eq:zLinearEq}).
Naturally, the immersion $\boldsymbol{W}_{\varepsilon}(\boldsymbol{y},\lambda)$
of the invariant manifold also depends on $\varepsilon$, which implicitly
assumes that $\text{\ensuremath{\mathcal{M}}}_{\lambda}$ persists
for $0\le\varepsilon\le1$. Whenever we write $\boldsymbol{F}_{0}$
or $\boldsymbol{W}_{0}$ we mean the $\varepsilon=0$ limit. 

Let us define the scaled Frechet derivative as
\[
\boldsymbol{A}_{\varepsilon}(\boldsymbol{y},\lambda)=\varepsilon D_{1}\boldsymbol{F}_{\varepsilon}(\boldsymbol{W}_{\varepsilon}(\boldsymbol{y},\lambda),\lambda).
\]
With this notation the corrected model (\ref{eq:zLinearEq}) becomes
\begin{equation}
\left.\begin{array}{rl}
\dot{\boldsymbol{y}} & =\boldsymbol{f}_{\varepsilon}(\boldsymbol{y},\lambda)\\
\dot{\boldsymbol{z}} & =\varepsilon^{-1}\boldsymbol{A}_{\varepsilon}(\boldsymbol{y},\lambda)\boldsymbol{z}-D_{2}\boldsymbol{W}_{\varepsilon}(\boldsymbol{y},\lambda)\dot{\lambda}\\
\lambda & =\begin{cases}
1 & \mbox{for}\;h(\boldsymbol{W}_{\varepsilon}(\boldsymbol{y},\lambda)+\boldsymbol{z})>0\\
-1 & \mbox{for}\;h(\boldsymbol{W}_{\varepsilon}(\boldsymbol{y},\lambda)+\boldsymbol{z})<0
\end{cases}
\end{array}\right\} .\label{eq:TimeScale-FullLinEq}
\end{equation}
Changing the time-scales by introducing $t=\varepsilon\theta$ we
get 
\begin{equation}
\left.\begin{array}{rl}
\mathring{\boldsymbol{y}} & =\varepsilon\boldsymbol{f}_{\varepsilon}(\boldsymbol{y},\lambda)\\
\mathring{\boldsymbol{z}} & =\boldsymbol{A}_{\varepsilon}(\boldsymbol{y},\lambda)\boldsymbol{z}-D_{2}\boldsymbol{W}_{\varepsilon}(\boldsymbol{y},\lambda)\mathring{\lambda}
\end{array}\right\} ,\label{eq:TimeScale-Rescaled}
\end{equation}
where $\mathring{\;}$ stand for $\nicefrac{\mathrm{d}}{\mathrm{d}\theta}$.
When setting $\varepsilon=0$ we arrive at the layer system 
\begin{equation}
\left.\begin{array}{rl}
\mathring{\boldsymbol{y}} & =\boldsymbol{0}\\
\mathring{\boldsymbol{z}} & =\boldsymbol{A}_{0}(\boldsymbol{y},\lambda)\boldsymbol{z}-D_{2}\boldsymbol{W}_{0}(\boldsymbol{y},\lambda)\mathring{\lambda}
\end{array}\right\} ,\label{eq:TimeScale-LayerSys}
\end{equation}
which stipulates that variable $\boldsymbol{y}$ is constant along
trajectories. We assume the following:
\begin{description}
\item [{\textbf{($\overline{\textrm{\textbf{A3}}}$)\label{A3BAR-UniqueSolutions}}}] Assumptions
\nameref{A-EvolutionOperator} holds when (\ref{eq:ASS-VariationalEquation})
is replaced by (\ref{eq:TimeScale-Rescaled}) for all $\varepsilon\in[0,1]$.
The unique solution of (\ref{eq:TimeScale-Rescaled}) can be written
as 
\[
\boldsymbol{z}(t)=\boldsymbol{U}_{\varepsilon}(t,s)\boldsymbol{z}(s)-\int_{s}^{t}\boldsymbol{K}_{\varepsilon}(t,\tau)\dot{\lambda}(\tau)\mathrm{d}\tau.
\]
\item [{\textbf{($\overline{\textrm{\textbf{A4}}}$)\label{A4BAR-Hyperbolicity}}}] Assumption
\nameref{A-NormalHyperbolicity} holds when (\ref{eq:ASS-VariationalEquation})
is replaced by (\ref{eq:TimeScale-Rescaled}) and $\sigma_{s}<-\varepsilon\sigma_{c}$.
\item [{\textbf{($\overline{\textrm{\textbf{A5}}}$)\label{A5BAR-NormalBundle}}}] The
perturbation $D_{2}\boldsymbol{W}_{0}(\boldsymbol{y},\lambda)$ acts
in the invariant normal bundle of $\mathcal{M}_{\lambda}$, that is,\textbf{
\begin{equation}
\Pi^{c}(\boldsymbol{y},\lambda)D_{2}\boldsymbol{W}_{0}(\boldsymbol{y},\lambda)=\boldsymbol{0}.\label{eq:Manifold Q condition}
\end{equation}
}
\end{description}
\begin{remark}
\label{rem:A5BAR-consequence}As a consequence of \nameref{A3BAR-UniqueSolutions}
and \nameref{A4BAR-Hyperbolicity}, $\boldsymbol{A}_{0}(\boldsymbol{y},\lambda)$
has an $n$ dimensional kernel spanned by $D_{1}\boldsymbol{W}_{0}(\boldsymbol{y},\lambda)$,
and $\Pi^{c}(\boldsymbol{y},\lambda)\boldsymbol{A}_{0}(\boldsymbol{y},\lambda)=\boldsymbol{0}$.
Because of \nameref{A5BAR-NormalBundle} and for $t\ge s$ we also
have 
\[
\left\Vert \boldsymbol{K}_{\varepsilon}(t,s)\right\Vert \le M\mathrm{e}^{\sigma_{s}(t-s)}.
\]
\end{remark}
We investigate the non-smooth dynamics for $\varepsilon=0$. The case
of constant $\lambda$ is trivial, because we have assumed that $\text{\ensuremath{\mathcal{M}}}_{\lambda}$
is attracting for $0\le\varepsilon\le1$. Next we consider the dynamics
in $\Sigma$, which is described by 

\begin{equation}
\left.\begin{array}{rl}
\mathring{\boldsymbol{y}} & =\boldsymbol{0}\\
\mathring{\boldsymbol{z}} & =\boldsymbol{A}_{0}(\boldsymbol{y},\lambda)\boldsymbol{z}-D_{2}\boldsymbol{W}_{0}(\boldsymbol{y},\lambda)\mathring{\lambda}\\
0 & =h(\boldsymbol{W}_{0}(\boldsymbol{y},\lambda)+\boldsymbol{z})
\end{array}\right\} .\label{eq:TimeScale-Stick}
\end{equation}
Any point in $\text{\ensuremath{\mathcal{M}}}_{\lambda}\cap\Sigma$,
i.e., $\boldsymbol{y}\in G$, $\boldsymbol{z}=\boldsymbol{0}$, $\lambda\in[-1,1]$
is an equilibrium of (\ref{eq:TimeScale-Stick}); therefore $\text{\ensuremath{\mathcal{M}}}_{\lambda}$
is invariant under all the dynamics for $\varepsilon=0$. It is however
not obvious whether $\text{\ensuremath{\mathcal{M}}}_{\lambda}\cap\Sigma$
is attracting for $\varepsilon=0$, which is addressed by the next
theorem.
\begin{theorem}
\label{thm:Normal-hyperbolicity}Assume \nameref{A-InvariantManifold},\nameref{A-DomainOfDefinition},\nameref{A3BAR-UniqueSolutions},\nameref{A4BAR-Hyperbolicity}
and $d^{-}(\boldsymbol{y},\boldsymbol{0},\lambda)\neq0$. Let $\text{\ensuremath{\mathcal{M}}}_{\lambda}^{crit}$
be a compact set, such that 
\begin{equation}
\text{\ensuremath{\mathcal{M}}}_{\lambda}^{crit}\subset\left\{ \left(\boldsymbol{y},\lambda\right)\in\text{\ensuremath{\mathcal{M}}}_{\lambda}\cap\Sigma\,:\,\sup\left\{ \text{real part of roots of }\Delta(s)\right\} <0\right\} ,\label{eq:TimeScale-CritManifDef}
\end{equation}
where
\begin{equation}
\Delta(s)=sDh(\boldsymbol{W}_{0}(\boldsymbol{y},\lambda))\cdot\left(s-\boldsymbol{A}_{0}(\boldsymbol{y},\lambda)\right)^{-1}D_{2}\boldsymbol{W}_{0}(\boldsymbol{y},\lambda)-d^{-}(\boldsymbol{y}_{0},\boldsymbol{0},\lambda_{0}).\label{eq:TimeScale-CharFunction}
\end{equation}
Then $\text{\ensuremath{\mathcal{M}}}_{\lambda}^{crit}$ is a normally
hyperbolic and attracting critical manifold of equation (\ref{eq:TimeScale-FullLinEq})
for $\varepsilon=0$.
\end{theorem}
\begin{svmultproof}
In order to calculate whether the critical manifold is attracting,
we linearize equation (\ref{eq:TimeScale-Stick}) by using $\lambda=\lambda_{0}+\alpha$
as a perturbation
\begin{align}
\mathring{\boldsymbol{y}} & =\boldsymbol{0},\label{eq:TimeScale-Layer-A}\\
\mathring{\boldsymbol{z}} & =\boldsymbol{A}_{0}(\boldsymbol{y}_{0},\lambda_{0})\boldsymbol{z}-D_{2}\boldsymbol{W}_{0}(\boldsymbol{y}_{0},\lambda_{0})\mathring{\alpha},\label{eq:TimeScale-Layer-B}\\
0 & =Dh(\boldsymbol{W}_{0}(\boldsymbol{y}_{0},\lambda_{0}))\cdot\left(D_{2}\boldsymbol{W}_{0}(\boldsymbol{y}_{0},\lambda_{0})\alpha+\boldsymbol{z}\right).\label{eq:TimeScale-Layer-C}
\end{align}
The initial conditions $\alpha(0)$ and $\boldsymbol{z}(0)$ are linked
through equation (\ref{eq:TimeScale-Layer-C}), such that 
\[
\alpha(0)=-\frac{Dh(\boldsymbol{W}_{0}(\boldsymbol{y}_{0},\lambda_{0}))\cdot\boldsymbol{z}(0)}{d^{-}(\boldsymbol{y}_{0},\boldsymbol{0},\lambda_{0})}.
\]
It is sufficient to show that $\alpha$ decays, because by assumptions
\nameref{A3BAR-UniqueSolutions} \nameref{A4BAR-Hyperbolicity} and
without forcing, the $\boldsymbol{z}$ component decays to a constant;
if the initial condition satisfies $\Pi^{s}(\boldsymbol{y},\lambda)\boldsymbol{z}(0)=\boldsymbol{0}$,
$\boldsymbol{z}$ decays to zero. Applying the Laplace transform to
(\ref{eq:TimeScale-Layer-B}) we find that 
\begin{equation}
\boldsymbol{z}(s)=\left(s-\boldsymbol{A}_{0}(\boldsymbol{y}_{0},\lambda_{0})\right)^{-1}\left(\boldsymbol{z}(0)-D_{2}\boldsymbol{W}_{0}(\boldsymbol{y}_{0},\lambda_{0})\left(s\alpha(s)-\alpha(0)\right)\right),\label{eq:TimeScale-zLaplace}
\end{equation}
where $s$ is the Laplace parameter. By substituting (\ref{eq:TimeScale-zLaplace})
into (\ref{eq:TimeScale-Layer-C}) we find
\begin{multline*}
Dh\left(\boldsymbol{W}_{0}(\boldsymbol{y}_{0},\lambda_{0})\right)\cdot\biggl(D_{2}\boldsymbol{W}_{0}(\boldsymbol{y}_{0},\lambda_{0})\alpha(s)\\
+\left(s-\boldsymbol{A}_{0}(\boldsymbol{y}_{0},\lambda_{0})\right)^{-1}\left(\boldsymbol{z}(0)-D_{2}\boldsymbol{W}_{0}(\boldsymbol{y}_{0},\lambda_{0})\left(s\alpha(s)-\alpha(0)\right)\right)\biggr)=0,
\end{multline*}
which can be rearranged into
\begin{equation}
\alpha(s)=\frac{Dh\left(\boldsymbol{W}_{0}(\boldsymbol{y}_{0},\lambda_{0})\right)\cdot\left(s-\boldsymbol{A}_{0}(\boldsymbol{y}_{0},\lambda_{0})\right)^{-1}\left(D_{2}\boldsymbol{W}(\boldsymbol{y}_{0},\lambda_{0})\alpha(0)+\boldsymbol{z}(0)\right)}{sDh\left(\boldsymbol{W}_{0}(\boldsymbol{y}_{0},\lambda_{0})\right)\cdot\left(s-\boldsymbol{A}_{0}(\boldsymbol{y}_{0},\lambda_{0})\right)^{-1}D_{2}\boldsymbol{W}(\boldsymbol{y}_{0},\lambda_{0})-d^{-}(\boldsymbol{y}_{0},\boldsymbol{0},\lambda_{0})}.\label{eq:TimeScale-alpha(s)}
\end{equation}
The asymptotic properties of $\alpha(t)$ are determined by the poles
of (\ref{eq:TimeScale-alpha(s)}). The poles of the numerator are
already given by the spectrum of $\boldsymbol{A}_{0}(\boldsymbol{y}_{0},\lambda_{0})$,
which is assumed to be in the left half of the complex plane because
$\text{\ensuremath{\mathcal{M}}}_{\lambda}$ is attracting for constant
$\lambda$. Therefore only the roots of the denominator can cause
instability, hence the condition that $\Delta(s)$ has roots in the
left half of the complex plane is sufficient.
\end{svmultproof}

\begin{remark}
The proof can be extended to calculate the initial and final values
of $\alpha$. According to the Laplace Final Value Theorem we have
$\lim_{t\to\infty}\alpha(t)=\lim_{s\to0}s\alpha(s)$. We observe that
\[
\lim_{s\to0}sDh\left(\boldsymbol{W}_{0}(\boldsymbol{y}_{0},\lambda_{0})\right)\cdot\left(s-\boldsymbol{A}_{0}(\boldsymbol{y}_{0},\lambda_{0})\right)^{-1}D_{2}\boldsymbol{W}(\boldsymbol{y}_{0},\lambda_{0})=\lim_{t\to\infty}\boldsymbol{K}_{0}(t,s)=0
\]
by assumption \nameref{A5BAR-NormalBundle} and remark \ref{rem:A5BAR-consequence}.
Therefore we have 
\[
\lim_{t\to\infty}\alpha(t)=-\frac{\lim_{t\to\infty}\mathrm{e}^{\boldsymbol{A}_{0}(\boldsymbol{y}_{0},\lambda_{0})t}\boldsymbol{z}(0)}{d^{-}(\boldsymbol{y}_{0},\boldsymbol{0},\lambda_{0})}
\]
and if $\Pi^{c}(\boldsymbol{y},\lambda)\boldsymbol{z}(0)=\boldsymbol{0}$
we also have $\lim_{t\to\infty}\alpha(t)=0$. Applying the Laplace
Initial Value Theorem to equation (\ref{eq:TimeScale-alpha(s)}) yields
\[
\lim_{t\downarrow0}\alpha(t)=\lim_{s\to\infty}s\alpha(s)=\frac{Dh\left(\boldsymbol{W}_{0}(\boldsymbol{y}_{0},\lambda_{0})\right)\cdot\boldsymbol{z}(0)+d^{+}(\boldsymbol{y},\boldsymbol{z},\lambda)\alpha(0)}{d^{\pm}(\boldsymbol{y},\boldsymbol{z},\lambda)},
\]
where $d^{+}(\boldsymbol{y},\boldsymbol{z},\lambda)$ is given by
(\ref{eq:INF-dPlus}). We also have $Dh\left(\boldsymbol{W}_{0}(\boldsymbol{y}_{0},\lambda_{0})\right)\cdot\boldsymbol{z}(0)=-d^{-}(\boldsymbol{y}_{0},\boldsymbol{0},\lambda_{0})\alpha(0)$
according to (\ref{eq:TimeScale-Layer-C}) and therefore $\lim_{t\downarrow0}\alpha(t)=\alpha(0)$,
which makes $\alpha$ continuous at $t=s$.
\end{remark}
\begin{remark}
Similar to remark \ref{rem:ASS-Persistence}, normal hyperbolicity
does not imply the persistence of $\mathcal{M}_{\lambda}^{crit}$
under variations in $\varepsilon$. The theorem of Bates, Lu and Zeng
\cite{Bates1998} suggests that the evolution operator $\boldsymbol{U}$
needs to be differentiable (among other conditions) for $\mathcal{M}_{\lambda}^{crit}$
to persist for small $\varepsilon>0$. Note that the nonlinear string
example in section \ref{sec:NLS} generates such a differentiable
$\boldsymbol{U}$ on $\boldsymbol{Z}$.
\end{remark}
\begin{remark}
When both regions of $\mathcal{M}_{\lambda}$, that is $\mathcal{M}_{\lambda}\cap\Sigma$
and $\mathcal{M}_{\lambda}\backslash\left(\mathcal{M}_{\lambda}\cap\Sigma\right)$,
persist for $\varepsilon>0$, they most likely become discontinuous
at the boundaries $\Sigma^{\pm}$, hence as a whole, $\mathcal{M}_{\lambda}$
does not persist. This is because the vector fields are discontinuous.
Therefore for $\varepsilon>0$, trajectories that followed one part
of $\mathcal{M}_{\lambda}$ must jump to the other part of $\mathcal{M}_{\lambda}$,
which induces fast transients that we are unable to characterize under
general settings.
\end{remark}

\subsection{\label{subsec:QualitativeModel}Qualitative approximation of normal
dynamics and the reduced order model}

A key difference between the skeleton model (\ref{eq:RED-skeleton})
and the corrected model (\ref{eq:zLinearEq}) is that they have unique
solutions under different conditions. This difference is caused by
the fact that the skeleton model does not take into account the normal
discontinuity gap $d^{\pm}$. To rectify the omission of $d^{\pm}$,
the skeleton model is extended by a scalar variable, which represents
the dynamics of the convolution kernel $\boldsymbol{K}$ in equation
(\ref{eq:ASS-ConvolutionRepresentation}). We call this extension
the reduced order model. It is then shown that the reduced order model
reproduces uniqueness of solutions and the existence of a critical
manifold under equivalent conditions to those of theorems \ref{thm:Cont-Uniqueness}
and \ref{thm:Normal-hyperbolicity}.

To simplify the ensuing analysis we assume that
\begin{description}
\item [{\textbf{(A6)\label{(A6)-h-Linearity}}}] $h(\boldsymbol{x})$ is
linear, therefore $h(\boldsymbol{x})=h(\boldsymbol{0})+Dh\cdot\boldsymbol{x}$,
where $Dh$ is a constant linear functional. 
\end{description}
Assumption \nameref{(A6)-h-Linearity} allows us to derive a scalar
representation of $\boldsymbol{z}(t)$ without worrying about a varying
$Dh(\boldsymbol{x})$. The switching between parts of the state space
depends on 
\begin{equation}
h(\boldsymbol{x})=h(\boldsymbol{0})+Dh\cdot\left(\boldsymbol{W}(\boldsymbol{y},\lambda)+\boldsymbol{z}\right).\label{eq:Conc-h}
\end{equation}
In what follows we approximate the scalar valued $Dh\cdot\boldsymbol{z}$
in (\ref{eq:Conc-h}) by a convolution integral. Combining equations
(\ref{eq:ASS-ConvolutionRepresentation}), (\ref{eq:REM-Convolution})
and $\boldsymbol{z}(0)=\boldsymbol{0}$ yields
\[
Dh\cdot\boldsymbol{z}(t)=-\int_{0}^{t}Dh\cdot\boldsymbol{U}(t,\vartheta)\lim_{\mu\to\infty}\mu\left(\mu-\varepsilon^{-1}\boldsymbol{A}_{\varepsilon}(\boldsymbol{y}(\vartheta),\lambda(\vartheta))\right)^{-1}D_{2}\boldsymbol{W}(\boldsymbol{y}(\vartheta),\lambda(\vartheta))\dot{\lambda}(\vartheta)\mathrm{d}\vartheta.
\]
In order to proceed either \nameref{A-PertTrajectories} or time-scale
separation with \nameref{A5BAR-NormalBundle} can be assumed. When
\nameref{A-PertTrajectories} is assumed we are restricted to use
$\varepsilon=1$ and if \nameref{A5BAR-NormalBundle} is assumed we
set $\sigma=\sigma_{s}$. Now we can approximate that
\[
Dh\cdot\boldsymbol{z}(t)\approx\int_{0}^{t}\mathrm{e}^{\varepsilon^{-1}\sigma\left(t-\vartheta\right)}\lim_{\mu\to\infty}\mu Dh\cdot\left(\mu-\varepsilon^{-1}\boldsymbol{A}_{\varepsilon}(\boldsymbol{y}(\vartheta),\lambda(\vartheta))\right)^{-1}D_{2}\boldsymbol{W}(\boldsymbol{y}(\vartheta),\lambda(\vartheta))\dot{\lambda}(\vartheta)\mathrm{d}\vartheta,
\]
which neglects trajectories in the normal bundle of $\mathcal{M}_{\lambda}$
that are decaying with exponents smaller than $\varepsilon^{-1}\sigma$.
Note that 
\[
\lim_{\mu\to\infty}\mu Dh\cdot\left(\mu-\varepsilon^{-1}\boldsymbol{A}_{\varepsilon}(\boldsymbol{y}(s),\lambda(s))\right)^{-1}D_{2}\boldsymbol{W}(\boldsymbol{y}(s),\lambda(s))=d^{+}(\boldsymbol{y}(s),\boldsymbol{0},\lambda(s)).
\]
By simply defining $d^{+}(\boldsymbol{y},\lambda)=d^{+}(\boldsymbol{y},\boldsymbol{0},\lambda)$
we get
\begin{equation}
Dh\cdot\boldsymbol{z}(t)\approx-\int_{0}^{t}d^{+}(\boldsymbol{y}\left(\vartheta\right),\lambda\left(\vartheta\right))\mathrm{e}^{\varepsilon^{-1}\sigma\left(t-\vartheta\right)}\dot{\lambda}(\vartheta)\mathrm{d}\vartheta.\label{eq:Conc-z-Approx-Last}
\end{equation}
After defining $\kappa=Dh\cdot\boldsymbol{z}(t)$, we find that the
approximation (\ref{eq:Conc-z-Approx-Last}) satisfies the differential
equation 
\begin{equation}
\dot{\kappa}=\varepsilon^{-1}\sigma\kappa-d^{+}(\boldsymbol{y},\lambda)\dot{\lambda}\label{eq:Conc-kappa-DE}
\end{equation}
 with initial condition $\kappa(0)=0$. The switching function (\ref{eq:Conc-h})
using the new variable $\kappa$ becomes
\begin{equation}
h(\boldsymbol{x})\approx h_{\varepsilon}(\boldsymbol{y},\kappa,\lambda)=h(\boldsymbol{0})+Dh\cdot\boldsymbol{W}(\boldsymbol{y},\lambda)+\kappa.\label{eq:Conc-h-epsilon}
\end{equation}
We can also re-define the switching manifolds
\begin{align*}
\Sigma_{\varepsilon} & =\left\{ \left(\boldsymbol{y},\kappa,\lambda\right)\in G\times\mathbb{R}\times[-1,1]:h_{\varepsilon}(\boldsymbol{y},\kappa,\lambda)=0\right\} ,\\
\Sigma_{\varepsilon}^{\pm} & =\left\{ \left(\boldsymbol{y},\kappa,\pm1\right)\in G\times\mathbb{R}\times[-1,1]:h_{\varepsilon}(\boldsymbol{y},\kappa,\pm1)=0\right\} .
\end{align*}
With this notation, the skeleton model extended with the approximate
normal dynamics becomes
\begin{equation}
\left.\begin{array}{rl}
\dot{\boldsymbol{y}} & =\boldsymbol{f}(\boldsymbol{y},\lambda)\\
\dot{\kappa} & =\varepsilon^{-1}\sigma\kappa-d^{+}(\boldsymbol{y},\lambda)\dot{\lambda}\\
\lambda & =\begin{cases}
1 & \mbox{for}\;h_{\varepsilon}(\boldsymbol{y},\kappa,\lambda)>0\\
-1 & \mbox{for}\;h_{\varepsilon}(\boldsymbol{y},\kappa,\lambda)<0
\end{cases}
\end{array}\right\} .\label{eq:Conc-FullSystem}
\end{equation}

\begin{definition}
We call equation (\ref{eq:Conc-FullSystem}) the reduced order model
of (\ref{eq:GenEq}).
\end{definition}
When $h_{\varepsilon}(\boldsymbol{y},\kappa,\lambda)\neq0$ the dynamics
of $\kappa$ is decoupled from the rest of the variables and $\kappa$
exponentially vanishes, because $\sigma<0$ - due to assumption \nameref{A-PertTrajectories}
or \nameref{A5BAR-NormalBundle}. When $h_{\varepsilon}(\boldsymbol{y},\kappa,\lambda)=0$,
we apply the same technique as in section \ref{sec:skeleton} to find
a differential equation for $\lambda$. We express that
\begin{equation}
\frac{\mathrm{d}}{\mathrm{d}t}h_{\varepsilon}(\boldsymbol{y},\kappa,\lambda)=Dh\cdot D_{1}\boldsymbol{W}(\boldsymbol{y},\lambda)\boldsymbol{f}(\boldsymbol{y},\lambda)+d^{-}(\boldsymbol{y},\lambda)\dot{\lambda}+\dot{\kappa}=0\label{eq:Conc-h-deri}
\end{equation}
in $\Sigma_{\varepsilon}$, where $d^{-}(\boldsymbol{y},\lambda)=d^{-}(\boldsymbol{y},\boldsymbol{0},\lambda)$
is defined by (\ref{eq:INF-dMinus}). We drop the $(\boldsymbol{y},\lambda)$
arguments and solve (\ref{eq:Conc-kappa-DE}) and (\ref{eq:Conc-h-deri})
for $\dot{\kappa}$ and $\dot{\lambda}$ to arrive at

\begin{equation}
\left.\begin{array}{rl}
\dot{\boldsymbol{y}} & =\boldsymbol{f}\\
\dot{\kappa} & =-{\displaystyle \frac{d^{+}Dh\cdot D_{1}\boldsymbol{W}\,\boldsymbol{f}+\varepsilon^{-1}d^{-}\sigma\kappa}{d^{\pm}}}\\
\dot{\lambda} & ={\displaystyle \frac{Dh\cdot D_{1}\boldsymbol{W}\,\boldsymbol{f}+\varepsilon^{-1}\sigma\kappa}{d^{\pm}}}
\end{array}\right\} ,\label{eq:Conc-Sigma-dynamics}
\end{equation}
which governs the dynamics on $\Sigma_{\varepsilon}$. 

We can now check that the reduced order model (\ref{eq:Conc-FullSystem})
has the same key properties as the corrected model (\ref{eq:zLinearEq}).
In what follows we outline the equivalents of theorems \ref{thm:Cont-Uniqueness}
and \ref{thm:Normal-hyperbolicity} for the reduced order model (\ref{eq:Conc-FullSystem}).
\begin{proposition}
\label{prop:Conc-Uniqueness}A trajectory $\mathcal{T}$ of the reduced
order model (\ref{eq:Conc-FullSystem}) with an end point at $(\boldsymbol{y}^{\star},\kappa^{\star},\lambda^{\star})\in\Sigma_{\varepsilon}^{\pm}$
has a unique continuation through $(\boldsymbol{y}^{\star},\kappa^{\star},\lambda^{\star})$
if 
\end{proposition}
\begin{enumerate}
\item $d^{\pm}(\boldsymbol{y}^{\star},0,\lambda^{\star})>0$ as defined
by equation (\ref{eq:dPlusMinusDefinition}) and
\item when trajectory $\mathcal{T}$ is not tangent to $\Sigma^{\pm}$ or
trajectory $\mathcal{T}$ is tangent to $\Sigma_{\varepsilon}^{\pm}$
and the of order of the tangency is not greater than the smoothness
order ($C^{p}$) of $h_{\varepsilon}$, that is, there exists $0<\ell\le p$
such that
\[
\frac{\mathrm{d}^{\ell}}{\mathrm{d}t^{\ell}}h_{\varepsilon}(\boldsymbol{y}(t),\kappa(t),\lambda^{\star})\vert_{\boldsymbol{y}=\boldsymbol{y}^{\star},\kappa=\kappa^{\star}}\neq0.
\]
\end{enumerate}
\begin{svmultproof}
The proof is the same as for theorem \ref{thm:ExtUnique} if we replace
$\left(\boldsymbol{y},\kappa\right)\to\boldsymbol{y}$ and $d^{\pm}\to-D_{2}h_{0}(\boldsymbol{y},\lambda)$.
\end{svmultproof}

\begin{proposition}
\label{prop:Conc-Hyperbolicity}Let 
\begin{equation}
\text{\ensuremath{\mathcal{M}}}_{\lambda}^{crit}\subset\left\{ (\boldsymbol{y},\kappa,\lambda)\in\Sigma_{\varepsilon}\,:\,\frac{\sigma d^{-}(\boldsymbol{y},\lambda)}{d^{\pm}(\boldsymbol{y},\lambda)}>0,\,\kappa=0\right\} \label{eq:Conc-CritManifDefinition}
\end{equation}
be a compact set for $\varepsilon=0$. Then $\text{\ensuremath{\mathcal{M}}}_{\lambda}^{crit}$
is an attracting critical manifold of equation (\ref{eq:Conc-Sigma-dynamics})
which persists for a sufficiently small $\varepsilon>0$. The dynamics
on the critical manifold is governed by the skeleton model (\ref{eq:RED-SwitchingEq}).
\end{proposition}
\begin{svmultproof}
First, time is rescaled by $t=\varepsilon\theta$ in equation (\ref{eq:Conc-Sigma-dynamics})
which yields 
\[
\left.\begin{array}{rl}
\mathring{\boldsymbol{y}} & =\varepsilon\boldsymbol{f}\\
\mathring{\kappa} & =-{\displaystyle \frac{\varepsilon d^{+}Dh\cdot D_{1}\boldsymbol{W}\,\boldsymbol{f}+d^{-}\sigma\kappa}{d^{\pm}}}\\
\mathring{\alpha} & ={\displaystyle \frac{\varepsilon Dh\cdot D_{1}\boldsymbol{W}\,\boldsymbol{f}+\sigma\kappa}{d^{\pm}}}
\end{array}\right\} ,
\]
where $\mathring{\;}$ stands for $\nicefrac{\mathrm{d}}{\mathrm{d}\theta}$.
Setting $\varepsilon\to0$ yields 
\begin{equation}
\mathring{\boldsymbol{y}}=\boldsymbol{0},\;\mathring{\kappa}=-\frac{d^{-}\sigma}{d^{\pm}}\kappa,\;\mathring{\lambda}=-\frac{\sigma}{d^{\pm}}\kappa.\label{eq:Conc-Layer-Eqn}
\end{equation}
Assuming initial conditions $\kappa(0)=\kappa_{0}$ and $\lambda(0)=\lambda_{0}$
of (\ref{eq:Conc-Layer-Eqn}) at $t=0$ we get $\lim_{t\to\infty}\kappa(t)=0$
and $\lim_{t\to\infty}\lambda(t)=\lambda_{0}-\kappa_{0}/d^{-}$, if
$\frac{d^{-}\sigma}{d^{\pm}}>0$. This means that the critical manifold
is attracting and normally hyperbolic. Therefore $\text{\ensuremath{\mathcal{M}}}_{\lambda}^{crit}$
persists for a sufficiently small $\varepsilon>0$, according to Fenichel
\cite{Fenichel}. 

While $\kappa=0$ on the critical manifold, $\lim_{t\to\infty}\lim_{\varepsilon\to0}\varepsilon^{-1}\kappa(t)$
may not be zero, that is, the limits $t\to\infty$ and $\varepsilon\to0$
do not commute. After introducing $\varepsilon\gamma=\kappa$, we
can write that 
\[
\left.\begin{array}{rl}
\dot{\boldsymbol{y}} & =\boldsymbol{f}\\
\varepsilon\dot{\gamma} & ={\displaystyle \frac{d^{+}Dh\cdot D_{1}\boldsymbol{W}\,\boldsymbol{f}-d^{-}\sigma\gamma}{d^{\pm}}}\\
\dot{\lambda} & ={\displaystyle \frac{Dh\cdot D_{1}\boldsymbol{W}\,\boldsymbol{f}+\sigma\gamma}{d^{\pm}}}
\end{array}\right\} .
\]
Setting $\varepsilon\to0$ and some algebraic manipulation yields
\[
\dot{\lambda}=-\frac{Dh\cdot D_{1}\boldsymbol{W}\,\boldsymbol{f}}{d^{-}},
\]
which is the same equation as (\ref{eq:RED-SwitchingEq}) of the skeleton
model.
\end{svmultproof}

\begin{remark}
We know that $\sigma<0$, because of assumption \nameref{A-PertTrajectories}
or \nameref{A5BAR-NormalBundle}. If proposition \ref{prop:Conc-Uniqueness}
also holds, $\text{\ensuremath{\mathcal{M}}}_{\lambda}^{crit}$ is
attracting when $d^{-}<0$. This is the same condition under which
solutions of the skeleton model (\ref{eq:RED-skeleton}) are unique
due to theorem \ref{thm:ExtUnique}.
\end{remark}
Next we investigate in what sense the reduced order model (\ref{eq:Conc-FullSystem})
is similar to the corrected model (\ref{eq:TimeScale-FullLinEq})
with time-scale separation. It turns out that on $\Sigma_{\varepsilon}$
the critical manifold is likely to be attracting or repelling under
the same conditions. The precise statement is in the following proposition.
\begin{proposition}
\label{prop:same-stability}Assume that $d^{\pm}(\boldsymbol{y},\boldsymbol{0},\lambda)>0$
and $Dh\cdot\boldsymbol{A}_{0}^{-1}(\boldsymbol{y},\lambda)D_{2}\boldsymbol{W}(\boldsymbol{y},\lambda)>0$
along a smooth curve $\gamma=\left\{ \left(\boldsymbol{y}(\alpha),\lambda(\alpha)\right)\in\Sigma_{0}:\alpha\in(-\delta,\delta)\right\} $
with $\left(\boldsymbol{y},\lambda\right)\in C^{1}\left((-\delta,\delta),\Sigma_{0}\right)$
and $\delta>0$. For $\varepsilon=0$, the stability of equilibria
along $\gamma$ changes through a zero root (saddle-node bifurcation)
at the same value(s) of $\alpha\in(-\delta,\delta)$ for both systems
(\ref{eq:Conc-Sigma-dynamics}) and (\ref{eq:TimeScale-Stick}).
\end{proposition}
\begin{svmultproof}
Because of the assumption $d^{\pm}(\boldsymbol{y},\boldsymbol{0},\lambda)>0$,
the stability of an equilibrium of (\ref{eq:Conc-Layer-Eqn}) purely
depends on $d^{-}$, i.e., the equilibrium is attracting when $d^{-}<0$.
On the other hand, substituting $s=0$ into $\Delta(s)$ as given
by (\ref{eq:TimeScale-CharFunction}), we note that $\Delta(0)=d^{-}(\boldsymbol{y},\lambda)$.
This means that we have a zero root of $\Delta(s)$ when $d^{-}=0$.
Next we show that this zero root of $\Delta(s)$ becomes of the same
sign as $d^{-}$ as $\boldsymbol{y}(\alpha),\lambda(\alpha)$ changes
along $\gamma$. Let us now assume that at $\alpha=0$ we have $d^{-}(\boldsymbol{y}(0),\lambda(0))=0$
and denote the root of $\Delta$ that smoothly depends on $\alpha$
by $s:(-\delta,\delta)\to\mathbb{R}$ and for which $s(0)=0$. We
denote the derivative with respect to $\alpha$ by $^{\prime}$ and
calculate the derivative of $s$ from the definition (\ref{eq:TimeScale-CharFunction}),
that is, 
\[
s^{\prime}(0)=\left(Dh\cdot\boldsymbol{A}_{0}^{-1}(\boldsymbol{y},\lambda)D_{2}\boldsymbol{W}(\boldsymbol{y},\lambda)\right)^{-1}\left(D_{\text{1}}d^{-}(\boldsymbol{y},\lambda)\boldsymbol{y}^{\prime}+D_{\text{2}}d^{-}(\boldsymbol{y},\lambda)\lambda^{\prime}\right),
\]
where we omitted that $\boldsymbol{y}$, $\lambda$ are evaluated
at $\alpha=0$. We also calculate the derivative $d^{-\prime}(\boldsymbol{y},\lambda)=D_{\text{1}}d^{-}(\boldsymbol{y},\lambda)\boldsymbol{y}^{\prime}+D_{\text{2}}d^{-}(\boldsymbol{y},\lambda)\lambda^{\prime}$
and notice that the derivative $s^{\prime}$ and $d^{-\prime}$ have
the same sign when $Dh\cdot\boldsymbol{A}_{0}^{-1}(\boldsymbol{y},\lambda)D_{2}\boldsymbol{W}(\boldsymbol{y},\lambda)>0$,
which proves the proposition.
\end{svmultproof}

\begin{remark}
In the next section, for the example of the nonlinear string, $d^{-}$
is a small parameter, which measures how well the equilibrium shape
of the string is approximated by a truncated Fourier series. The error
gets smaller with increasing number of terms in the truncated series,
therefore $d^{-}$ also gets smaller. Without damping or nonlinearity
$d^{-}$ entirely vanishes, as was the case in \cite{SzalaiMZfriction}.
When both $d^{-}$ and $\varepsilon$ vanish, we arrive at a system
that is subject to Utkin's closure in section \ref{subsec:Utkins-closure}.
If $d^{-}$ vanishes, but we have $d^{+}>0$ then for $\varepsilon>0$
the trajectories are still unique, but there is no critical manifold
in $\Sigma$ that is being perturbed. 
\end{remark}
\begin{remark}
A more rigorous analysis would inspect the dynamics in the perturbed
vector bundle corresponding to the near zero root of (\ref{eq:TimeScale-CharFunction})
for $\varepsilon=1$. If this dynamics has a Lyapunov exponent $\sigma_{0}$
such that $\sigma_{s}<-\left|\sigma_{0}\right|$ as in assumption
\nameref{A-NormalHyperbolicity}, then this perturbed vector bundle
could be attached to $\mathcal{M}_{\lambda}$, which would become
a normally hyperbolic invariant manifold of the corrected model (\ref{eq:zLinearEq})
in $\Sigma$.
\end{remark}

\section{\label{sec:NLS}A bowed nonlinear string model reduced to single
degree-of-freedom }

In this section we illustrate the theory through a non-trivial example.
In this example, the invariant manifold $\mathcal{M}_{\lambda}$ is
a linear subspace about an equilibrium that depends nonlinearly on
the switching parameter $\lambda$. The dynamics within the invariant
manifold given by $\boldsymbol{f}(\boldsymbol{y},\lambda)$ and the
switching function $h_{0}(\boldsymbol{y},\lambda)$ are also nonlinear,
which yields neither a Filippov nor an Utkin type model, but the skeleton
model described in section \ref{sec:skeleton}. In addition to the
nonlinearity we also include damping to make the invariant manifold
attracting. 

\begin{figure}
\begin{centering}
\includegraphics[width=0.7\linewidth]{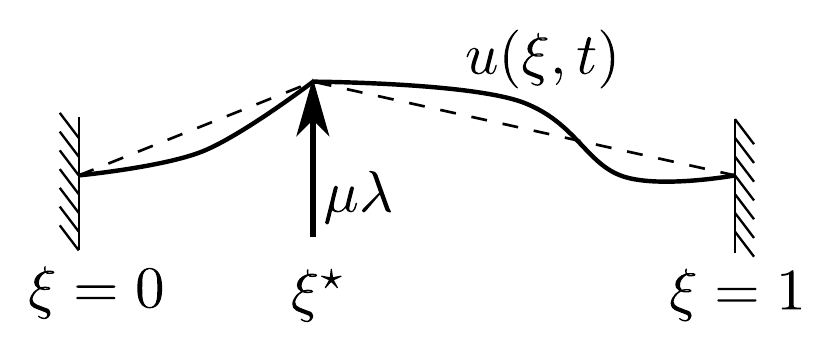}
\par\end{centering}
\caption{\label{fig:NLS-Schematic}Schematic of the nonlinear bowed string
model. The continuous line represents the deformation of the string
under vibration, the dashed line represents the equilibrium shape
of the string. $\mu\lambda$ represents the friction force between
the bow and the string, which acts at the contact point $\xi=\xi^{\star}$.}
\end{figure}

We consider a nonlinear string with both ends rigidly held as illustrated
in figure \ref{fig:NLS-Schematic}. The string has no resistance to
bending, any motion that occurs is due to the tension within the string.
Whenever lateral deformation occurs, the string becomes stretched,
which in turn causes an increase in tension and makes the model nonlinear.
The tension is uniform along the length of the string. We denote the
lateral deformation of the string by $u(\xi,t)$, where $\xi\in[0,1]$
represents the distance along the string and $t$ represents time.
Moreover, we assume that this deformation occurs within a fixed plane
so that $u$ is a scalar valued function. We also ignore any gravitational
effect. We use primes to denote differentiation with respect to $\xi$
and dots to denote differentiation in time. The dimensionless equation
under our simplifying assumptions is

\begin{equation}
\ddot{u}=Tu'',\quad T=1+\Gamma\left(\int_{0}^{1}\sqrt{1+u'^{2}}d\xi-1\right),\label{eq:NLS FullSqrt}
\end{equation}
where $T$ is the tension within the string and $\Gamma$ controls
the nonlinearity of the string. The boundary conditions are $u(0,t)=u(1,t)=0$
and 
\begin{equation}
u^{\prime}(\xi^{\star}-,t)-u^{\prime}(\xi^{\star}+,t)=\mu\lambda/T,\label{eq:NLS-innerBC}
\end{equation}
where $\mu$ is a friction coefficient, $\xi^{\star}$ is the position
of the contact point with the bow and $u^{\prime}(\xi^{\star}-,t)$,
$u^{\prime}(\xi^{\star}+,t)$ represent the left and right derivative
of $u$ with respect to $\xi$ at $\xi^{\star}$, respectively. The
boundary condition (\ref{eq:NLS-innerBC}) reflects the equilibrium
of forces at the contact point. The slope of the string together with
the tension forms a force vector on both sides of the contact point.
Since the string at the contact point is not smooth, the two force
vectors do not cancel and therefore to reach equilibrium an external
force is necessary, supplied by the friction force $\mu\lambda$.
The switching parameter $\lambda$ decides the direction of the friction
force and therefore changes sign as the relative velocity $h=v_{0}-\dot{u}(\xi^{\star},t)$
between the bow and the string reverses, that is,
\[
\lambda=\begin{cases}
1 & h>0\\
-1 & h<0
\end{cases}.
\]
To further simplify equation (\ref{eq:NLS FullSqrt}) we use second
order Taylor expansion, that is, $\sqrt{1+u'^{2}}\approx1+\frac{1}{2}u'^{2}$,
which gives us the equation
\begin{equation}
\ddot{u}=\left(1+\frac{\Gamma}{2}\int_{0}^{1}u'^{2}\mathrm{d}\xi\right)u'',\label{eq:NLS-u-undamp}
\end{equation}
with boundary conditions
\[
u(0,t)=u(1,t)=0,\;\left(1+\frac{\Gamma}{2}\int_{0}^{1}u'^{2}\mathrm{d}\xi\right)\left(u^{\prime}(\xi^{\star}-,t)-u^{\prime}(\xi^{\star}+,t)\right)=\mu\lambda.
\]
We require that $u(\cdot,t)\in\mathrm{Lip}\left([0,1],\text{\ensuremath{\mathbb{R}}}\right)$,
i.e., $u(\cdot,t)$ is Lipschitz continuous, which allows a finite
contact force on the string. We define the operator 
\[
\left(\mathfrak{D}^{2}v\right)(\xi)=-v^{\prime\prime}(\xi),\;v\in\left\{ v,v^{\prime\prime}\in\mathrm{Lip}\left([0,1],\text{\ensuremath{\mathbb{R}}}\right):v(0)=v(1)=0\right\} .
\]
The square root of $\mathfrak{D}^{2}$, can be represented on the
series $u=\sum a_{k}\sin k\pi\xi$ by $\mathfrak{D}u=\sum k\pi a_{k}\sin k\pi\xi$.
Note that $\mathfrak{D}$ is not producing the first order derivative.
To represent all boundary conditions, we define a restricted $\mathfrak{D}^{2}$
as 
\[
\overline{\mathfrak{D}}^{2}u=\mathfrak{D}^{2}u,\;u^{\prime}(\xi^{\star}-)-u^{\prime}(\xi^{\star}+)=\mu\lambda/T.
\]
We also introduce damping with a constant damping ratio $\beta\in[0,1)$
for all vibration modes which transforms equation (\ref{eq:NLS-u-undamp})
into 
\begin{equation}
\ddot{u}=-\left(1+\frac{\Gamma}{2}\int_{0}^{1}u'^{2}\mathrm{d}\xi\right)\overline{\mathfrak{D}}^{2}u-2\beta\mathfrak{D}\dot{u}.\label{eq:StringPDEDamp}
\end{equation}
Let us define $\boldsymbol{x}_{1}=u(\cdot),\,\boldsymbol{x}_{2}=\dot{u}(\cdot)$
and $\boldsymbol{x}=(\boldsymbol{x}_{1},\boldsymbol{x}_{2})^{T}$,
hence we can write the system (\ref{eq:NLS-u-undamp}) as the infinite
dimensional dynamical system
\begin{equation}
\dot{\boldsymbol{x}}=\boldsymbol{F}(\boldsymbol{x},\lambda)=\begin{pmatrix}\begin{array}{l}
\boldsymbol{x}_{2}\\
-\left(1+\frac{\Gamma}{2}\int_{0}^{1}\boldsymbol{x}_{1}^{\prime2}\mathrm{d}\xi\right)\overline{\mathfrak{D}}^{2}\boldsymbol{x}_{1}-2\beta\mathfrak{D}\boldsymbol{x}_{2}
\end{array}\end{pmatrix}\label{eq:NLS fullVF}
\end{equation}
and the switching function is 
\begin{equation}
h(\boldsymbol{x})=v_{0}-\boldsymbol{x}_{2}(\xi^{\star}).\label{eq:NLS h Function}
\end{equation}

In order to represent solutions that were encountered in section \ref{sec:TrivialExample},
we chose 
\[
\boldsymbol{X}=\mathrm{Lip}([0,1],\mathbb{R})\times L^{\infty}([0,1],\mathbb{R})
\]
for the phase space of (\ref{eq:NLS fullVF}), where $L^{\infty}$
stands for the space of bounded functions. The domain of definition
is 
\[
\boldsymbol{\mathcal{D}}(\boldsymbol{F})=\left\{ \left(\boldsymbol{x}_{1},\boldsymbol{x}_{2}\right)\in\boldsymbol{X}:\overline{\mathfrak{D}}^{2}\boldsymbol{x}_{1},\mathfrak{D}\boldsymbol{x}_{2}\in L^{\infty}([0,1],\mathbb{R}),\boldsymbol{x}_{2}\in\mathrm{Lip}([0,1],\mathbb{R})\right\} .
\]

In what follows we carry out a number of steps to arrive at the reduced
order model. These steps are applicable to systems where the invariant
manifold is a spectral submanifold of an equilibrium. The steps are 
\begin{enumerate}
\item Calculate the equilibrium of (\ref{eq:NLS fullVF}) as a function
of $\lambda$, which is denoted by $\boldsymbol{x}^{\star}$.
\item Find the smoothest two-dimensional spectral submanifold \cite{Haller2016}
$\mathcal{M}_{\lambda}$ about $\boldsymbol{x}^{\star}$, corresponding
to the pair of complex conjugate eigenvalues with the least negative
real part. The immersion of the manifold is denoted by $\boldsymbol{W}:\text{\ensuremath{\mathbb{R}}}^{2}\times[-1,1]\to\boldsymbol{X}$.
Assume a function $\boldsymbol{y}^{\star}:\text{\ensuremath{\mathbb{R}}}^{2}\times[-1,1]\to\text{\ensuremath{\mathbb{R}}}^{2}$,
which shifts the parametrization of $\mathcal{M}_{\lambda}$, such
that $\boldsymbol{W}(\boldsymbol{y},\lambda)=\boldsymbol{W}_{\mathit{fix}}(\boldsymbol{y}+\boldsymbol{y}^{\star}(\boldsymbol{y},\lambda),\lambda)$,
where $\boldsymbol{W}_{\mathit{fix}}$ is just one parametrization
of $\mathcal{M}_{\lambda}$. $\boldsymbol{y}^{\star}$ is an unknown
and will be calculated in step 4.
\item Introduce an artificial parameter $\varepsilon$, that slows down
the dynamics on $\mathcal{M}_{\lambda}$ to standstill at $\varepsilon=0$
and has no effect at $\varepsilon=1$. Then for $\varepsilon=0$ calculate
the invariant normal bundle of $\mathcal{M}_{\lambda}$, which is
formed by the subspace orthogonal to the kernel of the adjoint $\boldsymbol{A}_{0}^{\star}(\boldsymbol{y},\lambda)$
at each point on $\mathcal{M}_{\lambda}$.
\item Choose a coordinate shift $\boldsymbol{y}^{\star}$ so that $D_{2}\boldsymbol{W}$
falls into the invariant normal bundle of $\mathcal{M}_{\lambda}$
at $\varepsilon=0$, that is, $\boldsymbol{W}$ satisfies assumption
\nameref{A4BAR-Hyperbolicity}. This now fully specifies the immersion
$\boldsymbol{W}$.
\item Obtain the skeleton model by substituting the immersion $\boldsymbol{W}$
into (\ref{eq:NLS fullVF}).
\item Calculate the normal discontinuity gap from the dynamics in the invariant
normal bundle of $\mathcal{M}_{\lambda}$. Also determine $\sigma$,
the rate of convergence of the trajectory in the normal bundle with
initial condition $D_{2}\boldsymbol{W}$.
\end{enumerate}
\begin{proposition}
Following the six steps above yields the reduced order model of equations
(\ref{eq:NLS fullVF}) and (\ref{eq:NLS h Function}) in the form
of
\[
\begin{pmatrix}\dot{y}_{1}\\
\dot{y}_{2}\\
\varepsilon\dot{\kappa}
\end{pmatrix}=\begin{pmatrix}\begin{array}{l}
y_{2}+y_{2}^{\star}(y_{1},\lambda)\\
-\left(c^{2}(y_{1},\lambda)\left(\pi^{2}y_{1}-2\gamma(\lambda)\sin\pi\xi^{\star}\right)+2\beta\pi\left(y_{2}+y_{2}^{\star}(y_{1},\lambda)\right)\right)-D_{1}y_{2}^{\star}(y_{1},\lambda)\left(y_{2}+y_{2}^{\star}(y_{1},\lambda)\right)\\
\sigma\kappa+d^{+}(y_{1},\lambda)\dot{\lambda}
\end{array}\end{pmatrix}
\]
with switching function
\[
h_{\varepsilon}(y_{1},y_{2},\kappa,\lambda)=v_{0}-\left(y_{2}+y_{2}^{\star}(y_{1},\lambda)\right)\sin\pi\xi^{\star}-\varepsilon\kappa.
\]
The normal discontinuity gap is
\[
d^{\pm}(y_{1},\lambda)=\gamma'(\lambda)\frac{c^{2}(y_{1},\lambda)\cos^{-1}\frac{\beta}{c(y_{1},\lambda)}}{\pi\sqrt{c^{2}(y_{1},\lambda)-\beta^{2}}}\;\text{and}\;d^{+}(y_{1},\lambda)=d^{\pm}(y_{1},\lambda)-D_{2}y_{2}^{\star}(y_{1},\lambda)\sin\pi\xi^{\star}.
\]
The coordinate shift on the manifold in the velocity coordinate is
any function that satisfies the differential equation
\[
D_{2}y_{2}^{\star}(y_{1},\lambda)=\frac{4\Gamma\beta\gamma^{\prime}(\lambda)\gamma(\lambda)}{c^{2}(y_{1},\lambda)}\left(y_{1}-2\gamma(\lambda)\frac{\sin\pi\xi^{\star}}{\pi^{2}}\right)\sum_{k=2}^{\infty}\frac{\sin^{2}k\pi\xi^{\star}}{k^{3}\pi}.
\]
The instantaneous square of the wave speed at the contact point is
\begin{multline}
c^{2}(y_{1},\lambda)=1+\frac{\Gamma}{2}\left(\gamma^{2}(\lambda)\xi^{\star}(1-\xi^{\star})+\gamma(\lambda)\left(y_{1}-2\gamma(\lambda)\frac{\sin\pi\xi^{\star}}{\pi^{2}}\right)\sin\pi\xi^{\star}\right.,\\
\left.+\frac{\pi^{2}}{2}\left(y_{1}-2\gamma(\lambda)\frac{\sin\pi\xi^{\star}}{\pi^{2}}\right)^{2}\right)\label{eq:NLS-FinalWaveSpeed}
\end{multline}
where
\begin{align}
\gamma(\lambda) & =\frac{\left(b+9\Gamma^{2}\lambda\mu(1-\xi^{\star})^{2}\xi^{\star2}\right)^{2/3}-2\sqrt[3]{3}\Gamma(1-\xi^{\star})\xi^{\star}}{3^{2/3}\Gamma(1-\xi^{\star})\xi^{\star}\sqrt[3]{b+9\Gamma^{2}\lambda\mu(1-\xi^{\star})^{2}\xi^{\star2}}}\;\text{with}\label{eq:NLS-Gamma-Lambda}\\
b & =\sqrt{3}\sqrt{\Gamma^{3}(1-\xi^{\star})^{3}\xi^{\star3}\left(27\Gamma\lambda^{2}\mu^{2}(1-\xi^{\star})\xi^{\star}+8\right)}.
\end{align}
\end{proposition}
\begin{svmultproof}
These results are proven in lemmas \ref{lem:NLS-Immersion}, \ref{lem:NLS-Frechet},
\ref{lem:NLS-CoordinateShift}, \ref{lem:NLS-SkeletonModel} and \ref{lem:NLS-DiscontinuityGap}.
\end{svmultproof}

\subsection{The invariant manifold and its parametrization}

We identify the invariant manifold $\mathcal{M}_{\lambda}$ with the
spectral submanifold \cite{Haller2016} of the string's equilibrium
corresponding to its first natural frequency. When $\lambda$ is constant
the string has an equilibrium. We choose the smoothest invariant manifold
about the equilibrium corresponding to the first natural frequency
of the string, which is a unique two-dimensional linear subspace.
We note that the theory of Cabr\'e et al. \cite{CabreLlave2003}
does not apply, because damping makes backward-time solutions non-unique.
\begin{lemma}
\label{lem:NLS-Immersion}The immersion of the invariant manifold
$\mathcal{M}_{\lambda}$ about the equilibrium, as specified in steps
1 and 2 of the model reduction process is 

\begin{equation}
\boldsymbol{W}(\boldsymbol{y},\lambda)=\begin{pmatrix}\begin{array}{l}
\gamma(\lambda)\left(\xi(1-\xi^{\star})-H(\xi-\xi^{\star})(\xi-\xi^{\star})\right)+(y_{1}+y_{1}^{\star}(\boldsymbol{y},\lambda))\sin\pi\xi\\
(y_{2}+y_{2}^{\star}(\boldsymbol{y},\lambda))\sin\pi\xi
\end{array}\end{pmatrix},\label{eq:NLS-FullImmersion}
\end{equation}
where $\gamma(\lambda)$ is given by equation (\ref{eq:NLS-Gamma-Lambda}).
The coordinate shift $\boldsymbol{y}^{\star}=\left(y_{1}^{\star},y_{2}^{\star}\right)^{T}$
is not yet known. 
\end{lemma}
\begin{figure}
\begin{centering}
\includegraphics[width=0.7\linewidth]{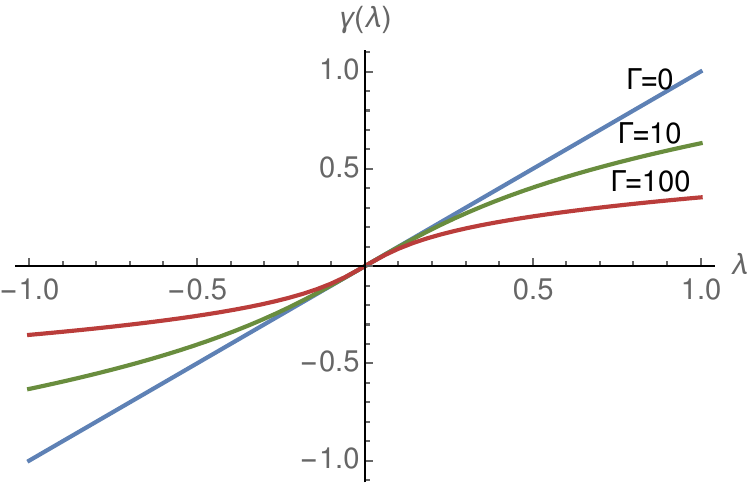}
\par\end{centering}
\caption{\label{fig:scGammaPlot}The graph of $\gamma(\lambda)$ for various
values of $\Gamma$. The other parameters are $\xi^{\star}=\sqrt{2}/2$
and $\mu=1$.}
\end{figure}

\begin{svmultproof}
We choose the representation of the invariant manifold as
\begin{equation}
\boldsymbol{x}=\boldsymbol{W}(\boldsymbol{y},\lambda)=\boldsymbol{x}^{\star}(\lambda)+\boldsymbol{W}_{1}\left(\boldsymbol{y}+\boldsymbol{y}^{\star}(\boldsymbol{y},\lambda)\right),\label{eq:NLS-immersion}
\end{equation}
where
\[
\boldsymbol{W}_{1}=\begin{pmatrix}\sin\pi\xi & 0\\
0 & \sin\pi\xi
\end{pmatrix}.
\]
A substitution of $\boldsymbol{W}$ into the invariance equation (\ref{eq:ASS-Invariance})
shows that $\boldsymbol{W}$ is indeed an immersion of an invariant
manifold and corresponds to the first natural frequency. Because $\boldsymbol{W}$
is linear in $\boldsymbol{y}$, $\mathcal{M}_{\lambda}$ is also the
smoothest invariant manifold.

The equilibrium $\boldsymbol{x}^{\star}(\lambda)$ is calculated by
setting the time-derivative to zero in equation (\ref{eq:NLS-u-undamp}),
which yields
\begin{equation}
\left(1+\frac{\Gamma}{2}\int_{0}^{1}\boldsymbol{x}_{1}'^{2}\mathrm{d}\xi\right)\boldsymbol{x}_{1}''=0,\;\left(1+\frac{\Gamma}{2}\int_{0}^{1}\boldsymbol{x}_{1}'^{2}\mathrm{d}\xi\right)\left(\boldsymbol{x}_{1}'(\xi^{\star}-)-\boldsymbol{x}_{1}'(\xi^{\star}+)\right)=\lambda\mu.\label{eq:NLSEquilibr}
\end{equation}
In equation (\ref{eq:NLSEquilibr}), $\frac{\Gamma}{2}\int_{0}^{1}\boldsymbol{x}_{1}'^{2}\mathrm{d}\xi$
is independent of $\xi$, therefore integrating (\ref{eq:NLSEquilibr})
twice and applying the boundary conditions we get 
\begin{equation}
\boldsymbol{x}_{1}^{\star}=\lambda\mu\left(1+\frac{\Gamma}{2}\int_{0}^{1}\left(\boldsymbol{x}_{1}^{\star}\right)^{\prime2}\mathrm{d}\xi\right)^{-1}\left(\xi(1-\xi^{\star})-H(\xi-\xi^{\star})(\xi-\xi^{\star})\right),\label{eq:NLS-x1star-implicit}
\end{equation}
which still needs to be solved for $\boldsymbol{x}_{1}^{\star}$.
We define 
\begin{equation}
\gamma=\lambda\mu\left(1+\frac{\Gamma}{2}\int_{0}^{1}\left(\boldsymbol{x}_{1}^{\star}\right)^{\prime2}\mathrm{d}\xi\right)^{-1},\label{eq:NLS-gamma-def}
\end{equation}
which yields 
\begin{equation}
\boldsymbol{x}_{1}^{\star}=\gamma(\lambda)\left(\xi(1-\xi^{\star})-H(\xi-\xi^{\star})(\xi-\xi^{\star})\right).\label{eq:NLS-x1star-gamma}
\end{equation}
Physically $\gamma(\lambda)\xi^{\star}(1-\xi^{\star})$ is the displacement
of the string at the contact point at the equilibrium. To find the
equation for $\gamma$ we substitute (\ref{eq:NLS-x1star-gamma})
into (\ref{eq:NLS-gamma-def}). We then evaluate the integral in (\ref{eq:NLS-gamma-def}),
that is,
\[
\int_{0}^{1}\left[\frac{\mathrm{d}}{\mathrm{d}\xi}\left(\xi(1-\xi^{\star})-H(\xi-\xi^{\star})(\xi-\xi^{\star})\right)\right]^{2}\mathrm{d}\xi=\xi^{\star}(1-\xi^{\star}),
\]
so that equation (\ref{eq:NLS-gamma-def}) becomes
\begin{equation}
\left(1+\frac{\Gamma}{2}\gamma^{2}\xi^{\star}(1-\xi^{\star})\right)\gamma=\lambda\mu.\label{eq:NLS gammaEquation}
\end{equation}
Equation (\ref{eq:NLS gammaEquation}) can be solved for $\gamma$
with a single real solution, which is given by equation (\ref{eq:NLS-Gamma-Lambda}),
that makes the equilibrium fully specified. Figure \ref{fig:scGammaPlot}
shows the values of $\gamma$ for various levels of nonlinearity.
The nonlinearity is hardening, because the string deforms less than
it would under the same force with the linear model. Substituting
the equilibrium into (\ref{eq:NLS-immersion}) yields equation (\ref{eq:NLS-FullImmersion}).
\end{svmultproof}

\subsection{Linearized dynamics about the invariant manifold}

The linearized dynamics about $\mathcal{M}_{\lambda}$ is characterized
by the Frechet derivative $\boldsymbol{A}_{1}(\boldsymbol{y},\lambda)$
of equation (\ref{eq:NLS fullVF}), which is calculated here. 
\begin{lemma}
\label{lem:NLS-Frechet}The Frechet derivative of $\boldsymbol{F}$
evaluated on $\mathcal{M}_{\lambda}$ is

\begin{equation}
\boldsymbol{A}_{1}(\boldsymbol{y},\lambda)\boldsymbol{z}=\begin{pmatrix}\begin{array}{l}
\boldsymbol{z}_{2}\\
\Gamma\int_{0}^{1}\overline{\boldsymbol{x}}_{1}^{\prime\prime}\boldsymbol{z}_{1}\mathrm{d}\xi\mathfrak{D}^{2}\overline{\boldsymbol{x}}_{1}-c^{2}(\boldsymbol{y},\lambda)\mathfrak{D}^{2}\boldsymbol{z}_{1}-2\beta\mathfrak{D}\boldsymbol{z}_{2}
\end{array}\end{pmatrix},\label{eq:NLS FrechetDeri}
\end{equation}
where
\begin{equation}
\int_{0}^{1}\overline{\boldsymbol{x}}_{1}^{\prime\prime}\boldsymbol{z}_{1}\mathrm{d}\xi=-\gamma(\lambda)\boldsymbol{z}_{1}(\xi^{\star})-\pi^{2}\left(y_{1}+y_{1}^{\star}(\boldsymbol{y},\lambda)\right)\int_{0}^{1}\boldsymbol{z}_{1}\sin\pi\xi\mathrm{d}\xi,\label{eq:NLS-coupling}
\end{equation}
and the instantaneous square of the wave speed on $\mathcal{M}_{\lambda}$
at the contact point is
\begin{multline}
c^{2}(\boldsymbol{y},\lambda)=1+\frac{\Gamma}{2}\left(\gamma^{2}(\lambda)\xi^{\star}(1-\xi^{\star})+\gamma(\lambda)\left(y_{1}+y_{1}^{\star}(\boldsymbol{y},\lambda)\right)\sin\pi\xi^{\star}\right.\\
\left.+\frac{\pi^{2}}{2}\left(y_{1}+y_{1}^{\star}(\boldsymbol{y},\lambda)\right)^{2}\right)\label{eq:NLS WaveSpeed}
\end{multline}
and 
\begin{equation}
\mathfrak{D}^{2}\overline{\boldsymbol{x}}_{1}=\pi^{2}\left(y_{1}+y_{1}^{\star}(\boldsymbol{y},\lambda)\right)\sin\pi\xi.\label{eq:NLS-Dsqx1bar}
\end{equation}
The domain of definition of $\boldsymbol{A}_{1}(\boldsymbol{y},\lambda)$
is 
\begin{equation}
\boldsymbol{\mathcal{D}}=\left\{ \left(\boldsymbol{x}_{1},\boldsymbol{x}_{2}\right)\in\boldsymbol{X}:\boldsymbol{x}_{1}^{\prime\prime},\boldsymbol{x}_{2}^{\prime\prime}\in L^{\infty}([0,1],\mathbb{R}),\boldsymbol{x}_{2}\in\mathrm{Lip}([0,1],\mathbb{R})\right\} \label{eq:NLS-DomainOfDefinition}
\end{equation}
and 
\begin{equation}
\boldsymbol{Z}=\overline{\boldsymbol{\mathcal{D}}}=C^{1}([0,1],\mathbb{R})\times C^{0}([0,1],\mathbb{R}).\label{eq:NLS-DomainClosure}
\end{equation}
\end{lemma}
\begin{svmultproof}
The only term in equation (\ref{eq:NLS fullVF}) not already linear
is 
\begin{equation}
\left(1+\frac{\Gamma}{2}\int_{0}^{1}\boldsymbol{x}_{1}'^{2}\mathrm{d}\xi\right)\mathfrak{D}^{2}\boldsymbol{x}_{1},\label{eq:NLS NonlinPart}
\end{equation}
which we now linearize about a general point $(\overline{\boldsymbol{x}}_{1},\overline{\boldsymbol{x}}_{2})\in\mathcal{M}_{\lambda}$.
The expression (\ref{eq:NLS NonlinPart}) is a product, hence we use
the product rule when differentiating it with respect to $\boldsymbol{x}_{1}$.
First we linearize (\ref{eq:NLS NonlinPart}) about $\overline{\boldsymbol{x}}_{1}$
and get

\[
\left[D_{\boldsymbol{x}_{1}}\left(1+\frac{\Gamma}{2}\int_{0}^{1}\boldsymbol{x}_{1}'^{2}\mathrm{d}\xi\right)\right]\boldsymbol{z}_{1}=\Gamma\int_{0}^{1}\overline{\boldsymbol{x}}_{1}'\boldsymbol{z}_{1}'\mathrm{d}\xi=-\Gamma\int_{0}^{1}\overline{\boldsymbol{x}}_{1}^{\prime\prime}\boldsymbol{z}_{1}\mathrm{d}\xi,
\]
where we have used that $\boldsymbol{z}_{1}$ must vanish at the boundaries
$\xi=0,1$. Therefore the first order Taylor expansion of (\ref{eq:NLS NonlinPart})
is
\begin{equation}
\left(1+\frac{\Gamma}{2}\int_{0}^{1}\boldsymbol{x}_{1}'^{2}\mathrm{d}\xi\right)\mathfrak{D}^{2}\boldsymbol{x}_{1}=-\Gamma\int_{0}^{1}\overline{\boldsymbol{x}}_{1}^{\prime\prime}\boldsymbol{z}_{1}\mathrm{d}\xi\mathfrak{D}^{2}\bar{\boldsymbol{x}}_{1}+\left(1+\frac{\Gamma}{2}\int_{0}^{1}\bar{\boldsymbol{x}}_{1}'^{2}\mathrm{d}\xi\right)\mathfrak{D}^{2}\boldsymbol{z}_{1}+\mathcal{O}(\boldsymbol{z}^{2}).\label{eq:NLS NonlinPart Expansion}
\end{equation}
The value $\overline{\boldsymbol{x}}_{1}$ is the first component
of the immersion of $\mathcal{M}_{\lambda}$,
\[
\overline{\boldsymbol{x}}_{1}(\xi)=\left[\boldsymbol{W}(\boldsymbol{y},\lambda)\right]_{1}=\gamma(\lambda)\left(\xi(1-\xi^{\star})-H(\xi-\xi^{\star})(\xi-\xi^{\star})\right)+\left(y_{1}+y_{1}^{\star}(\boldsymbol{y},\lambda)\right)\sin\pi\xi
\]
and when applying $\mathfrak{D}^{2}$, we get (\ref{eq:NLS-Dsqx1bar}).
The remaining term in the Taylor approximation (\ref{eq:NLS NonlinPart Expansion})
is
\[
\int_{0}^{1}\overline{\boldsymbol{x}}_{1}^{\prime2}(\xi)\mathrm{d}\xi=\gamma^{2}(\lambda)\xi^{\star}(1-\xi^{\star})+\gamma(\lambda)\left(y_{1}+y_{1}^{\star}(\boldsymbol{y},\lambda)\right)\sin\pi\xi^{\star}+\frac{\pi^{2}}{2}\left(y_{1}+y_{1}^{\star}(\boldsymbol{y},\lambda)\right)^{2}.
\]
We also define the square of the instantaneous wave-speed on $\mathcal{M}_{\lambda}$
and at the contact point as 
\[
c^{2}(\boldsymbol{y},\lambda)=1+\frac{\Gamma}{2}\int_{0}^{1}\boldsymbol{x}_{1}'^{2}\mathrm{d}\xi,
\]
which, after substitution becomes (\ref{eq:NLS WaveSpeed}). Gathering
all linear terms the Frechet derivative on the invariant manifold
becomes equation (\ref{eq:NLS FrechetDeri}).

The domain of definition of $\boldsymbol{A}_{1}(\boldsymbol{y},\lambda)$
must now include that $\boldsymbol{x}_{1}^{\prime}(\xi^{\star}-)-\boldsymbol{x}_{1}^{\prime}(\xi^{\star}+)=0$,
because the equilibrium is already included in the definition of $\boldsymbol{W}$.
However, $\boldsymbol{x}_{1}^{\prime}(\xi^{\star}-)-\boldsymbol{x}_{1}^{\prime}(\xi^{\star}+)=0$
is already satisfied if $\boldsymbol{x}_{1}^{\prime\prime}\in L^{\infty}([0,1],\mathbb{R})$,
therefore we arrive at (\ref{eq:NLS-DomainOfDefinition}). To determine
the closure $\overline{\boldsymbol{\mathcal{D}}},$ we first note
that $\boldsymbol{x}_{1}^{\prime}$ must be Lipschitz continuous and
therefore $\boldsymbol{x}_{1}$ is Lipschitz continuously differentiable.
The closure of such functions in the Lipschitz norm are the continuously
differentiable functions $C^{1}([0,1],\mathbb{R})$. Since $\mathfrak{D}$
is the square root of $\mathfrak{D}^{2}$ they have the same domain
of definition, therefore $\boldsymbol{x}_{2}$ is Lipschitz continuously
differentiable in $L^{\infty}$. The closure for this set in the $L^{\infty}$
norm is the set of continuous functions $C^{0}([0,1],\mathbb{R})$.
Summing up this argument we have found (\ref{eq:NLS-DomainClosure}).
Lipschitz functions are not dense in $C^{1}$ and continuous functions
are not dense in the set of bounded functions either, therefore $\boldsymbol{Z}\neq\boldsymbol{X}$.
\end{svmultproof}

\subsection{Invariant normal bundle and time-scale separation}

There is no small parameter in equation (\ref{eq:NLS fullVF}) that
controls the spectral gap between the tangential and normal dynamics
about $\mathcal{M}_{\lambda}$. We therefore introduce such a scaling
by artificially constructing $\boldsymbol{A}_{\varepsilon}(\boldsymbol{y},\lambda)$
such that for $\varepsilon=1$ we recover the original dynamics and
for $\varepsilon=0$ the tangential dynamics becomes $\dot{\boldsymbol{y}}=\boldsymbol{0}$
when time is re-scaled. This allows us to calculate the invariant
normal bundle of $\mathcal{M}_{\lambda}$ at $\varepsilon=0$ and
determine the parametrization $\mathcal{M}_{\lambda}$ (i.e., the
unknown coordinate shift $\boldsymbol{y}^{\star}(\boldsymbol{y},\lambda)$
in the immersion of $\mathcal{M}_{\lambda}$) such that $D_{2}\boldsymbol{W}(\boldsymbol{y},\lambda)$
is strictly in the invariant normal bundle of $\mathcal{M}_{\lambda}$.
\begin{lemma}
\label{lem:NLS-CoordinateShift}Applying steps 3 and 4 of the model
reduction procedure we find that the coordinate shift $\boldsymbol{y}^{\star}=\left(y_{1}^{\star},y_{2}^{\star}\right)^{T}$
becomes
\begin{equation}
y_{1}^{\star}(\lambda)=-2\gamma(\lambda)\frac{\sin\pi\xi^{\star}}{\pi^{2}}\label{eq:NLS y1Star}
\end{equation}
and
\begin{equation}
D_{2}y_{2}^{\star}(y_{1},\lambda)=\frac{4\Gamma\beta\gamma^{\prime}(\lambda)\gamma(\lambda)}{c^{2}(y_{1},\lambda)}\left(y_{1}-2\gamma(\lambda)\frac{\sin\pi\xi^{\star}}{\pi^{2}}\right)\sum_{k=2}^{\infty}\frac{\sin^{2}k\pi\xi^{\star}}{k^{3}\pi},\label{eq:NLS y2StarPrime}
\end{equation}
whose solution is 
\begin{multline}
y_{2}^{\star}(y_{1},\lambda)=\frac{16\beta\sum_{k=2}^{\infty}\frac{\sin^{2}k\pi\xi^{\star}}{k^{3}\pi^{3}}}{(1-\xi^{\star})\xi^{\star}}\times\\
\times\left(\frac{\left(b\left(d^{2}-2e\right)-ad\right)\tan^{-1}\left(\frac{2\gamma(\lambda)+d}{\sqrt{4e-d^{2}}}\right)}{\sqrt{4e-d^{2}}}+\frac{1}{2}(a-bd)\log(e+\gamma(\lambda)(\gamma(\lambda)+d))+b\gamma(\lambda)\right),\label{eq:NLS-y2Star}
\end{multline}
where
\[
a=-\pi^{2}y_{1},\,b=2\sin(\pi\xi^{\star}),\,d=\frac{y_{1}\sin(\pi\xi^{\star})}{(\xi^{\star}-1)\xi^{\star}},\,e=\frac{\pi^{2}\Gamma y_{1}^{2}+4}{2\Gamma\xi^{\star}-2\Gamma\xi^{\star}{}^{2}}.
\]
Using the coordinate shift (\ref{eq:NLS y1Star}), the square of the
instantaneous wave speed (\ref{eq:NLS WaveSpeed}) becomes (\ref{eq:NLS-FinalWaveSpeed}).
\end{lemma}
\begin{svmultproof}
Even though the mode shapes of the nonlinear string are the orthogonal
harmonic functions $\sin k\pi\xi$, the contact force $\lambda$ at
$\xi^{\star}$ makes these modes intricately coupled. This coupling
is represented by the term $\gamma(\lambda)\boldsymbol{z}_{1}(\xi^{\star})$
in equation (\ref{eq:NLS-coupling}). Nevertheless, we project $\boldsymbol{A}_{1}(\boldsymbol{y},\lambda)$
into two subspaces using the projection operators $\boldsymbol{P}:\boldsymbol{X}\to T_{\boldsymbol{y}}\mathcal{M}_{\lambda}$,
\begin{equation}
\boldsymbol{P}\boldsymbol{x}=\begin{pmatrix}2\sin\pi\xi\int_{0}^{1}\boldsymbol{x}_{1}(\Xi)\sin\pi\Xi\mathrm{d}\Xi\\
2\sin\pi\xi\int_{0}^{1}\boldsymbol{x}_{2}(\Xi)\sin\pi\Xi\mathrm{d}\Xi
\end{pmatrix}\label{eq:NLS P definition}
\end{equation}
and $\boldsymbol{Q}=\boldsymbol{I}-\boldsymbol{P}$. We calculate
the projected operators
\[
\boldsymbol{B}_{1}=\boldsymbol{P}\boldsymbol{A}_{1}(\boldsymbol{y},\lambda)\boldsymbol{P},\;\boldsymbol{B}_{12}=\boldsymbol{P}\boldsymbol{A}_{1}(\boldsymbol{y},\lambda)\boldsymbol{Q},\;\boldsymbol{B}_{21}=\boldsymbol{Q}\boldsymbol{A}_{1}(\boldsymbol{y},\lambda)\boldsymbol{P},\;\boldsymbol{B}_{2}=\boldsymbol{Q}\boldsymbol{A}_{1}(\boldsymbol{y},\lambda)\boldsymbol{Q},
\]
where we can show that $\boldsymbol{B}_{21}=\boldsymbol{0}$. Introducing
time-scale separation is just a multiplication of matrix $\boldsymbol{B}_{1}$
by $\varepsilon$, which represents the rescaled linearized dynamics
in the (invariant) tangent bundle of $\mathcal{M}_{\lambda}$. In
the new coordinate system the scaled linear operator becomes
\begin{equation}
\boldsymbol{A}_{\varepsilon}(\boldsymbol{y},\lambda)=\begin{pmatrix}\varepsilon\boldsymbol{B}_{1} & \boldsymbol{B}_{12}\\
\boldsymbol{0} & \boldsymbol{B}_{2}
\end{pmatrix}.\label{eq:NLS-A0-partition}
\end{equation}
Given the form of $\boldsymbol{A}_{\varepsilon}$, the bundle projections
assume the form
\[
\Pi^{s}=\begin{pmatrix}\boldsymbol{0} & \boldsymbol{C}\\
\boldsymbol{0} & \boldsymbol{I}
\end{pmatrix},\;\Pi^{c}=\boldsymbol{I}-\Pi^{s}=\begin{pmatrix}\boldsymbol{I} & -\boldsymbol{C}\\
\boldsymbol{0} & \boldsymbol{0}
\end{pmatrix},
\]
where $\boldsymbol{C}$ is an unknown operator. Expanding the bundle
invariance equation (\ref{eq:ASS-StableBundle}) yields
\begin{equation}
\begin{pmatrix}\boldsymbol{C}\boldsymbol{w}\\
\boldsymbol{w}
\end{pmatrix}=\begin{pmatrix}\varepsilon\boldsymbol{B}_{1} & \boldsymbol{B}_{12}\\
\boldsymbol{0} & \boldsymbol{B}_{2}
\end{pmatrix}\begin{pmatrix}\boldsymbol{C}\boldsymbol{v}\\
\boldsymbol{v}
\end{pmatrix},\label{eq:NLS-NormalBund-InvCond}
\end{equation}
which must hold for all $\boldsymbol{v}$, $\boldsymbol{P}\boldsymbol{v}=\boldsymbol{0}$.
Further expanding (\ref{eq:NLS-NormalBund-InvCond}) we get an equation
for $\boldsymbol{C}$ in the form of
\[
\boldsymbol{C}\boldsymbol{B}_{2}-\varepsilon\boldsymbol{B}_{1}\boldsymbol{C}=\boldsymbol{B}_{12}.
\]
The solution is $\boldsymbol{C}=\boldsymbol{B}_{12}\boldsymbol{B}_{2}^{-1}$
at the critical parameter value $\varepsilon=0$. We can now introduce
another coordinate system in which
\[
\boldsymbol{x}_{\ell}=\sum_{k=1}^{\infty}a_{\ell k}\sin k\pi\xi,\;\ell=1,2.
\]
We denote this transformation by $\boldsymbol{x}=\boldsymbol{T}\boldsymbol{a}$,
where $\boldsymbol{a}=\left(\boldsymbol{a}_{1},\boldsymbol{a}_{2}\right)^{T}$
with $\boldsymbol{a}_{\ell}=(a_{\ell1},a_{\ell2},\ldots)^{T}$, $\ell=1,2$.
The projections can be written as 
\[
\boldsymbol{T}^{-1}\boldsymbol{P}\boldsymbol{T}\boldsymbol{a}=(a_{11},a_{12})\;\text{and}\;\boldsymbol{T}^{-1}\boldsymbol{Q}\boldsymbol{T}\boldsymbol{a}=\left((a_{12},a_{13},\ldots),(a_{22},a_{23},\ldots)\right)^{T}.
\]
In this new coordinate system we have the operators
\[
\boldsymbol{B}_{12}\boldsymbol{T}\boldsymbol{a}=\begin{pmatrix}\begin{array}{l}
0\\
-\Gamma\left(\gamma(\lambda)\sum_{k=2}^{\infty}a_{1,k}\sin k\pi\xi^{\star}\right)\pi^{2}\left(y_{1}+y_{1}^{\star}(\boldsymbol{y},\lambda)\right)\sin\pi\xi
\end{array}\end{pmatrix},
\]
\[
\boldsymbol{B}_{2}\boldsymbol{T}\boldsymbol{a}=\begin{pmatrix}\begin{array}{l}
\sum_{k=2}^{\infty}a_{2,k}\sin k\pi\xi\\
-c^{2}(\boldsymbol{y},\lambda)\sum_{k=2}^{\infty}a_{1,k}k^{2}\pi^{2}\sin k\pi\xi-2\beta\sum_{k=2}^{\infty}a_{2,k}k\pi\sin k\pi\xi
\end{array}\end{pmatrix},
\]
\[
\boldsymbol{T}^{-1}\boldsymbol{B}_{2}\boldsymbol{T}=\begin{pmatrix}\boldsymbol{0} & \boldsymbol{I}\\
-c^{2}(\boldsymbol{y},\lambda)\pi^{2}\mathrm{diag}_{k=2}^{\infty}k^{2} & -2\beta\pi\mathrm{diag}_{k=2}^{\infty}k
\end{pmatrix},
\]
where $\mathrm{diag}_{k=2}^{\infty}k^{2}$ means an infinite diagonal
matrix with elements $k^{2}$ in the diagonal. The inverse $\boldsymbol{B}_{2}^{-1}$
is represented by 
\[
\boldsymbol{T}^{-1}\boldsymbol{B}_{2}^{-1}\boldsymbol{T}=\begin{pmatrix}-2\beta c^{-2}(\boldsymbol{y},\lambda)\pi^{-1}\mathrm{diag}_{k=2}^{\infty}k^{-1} & -c^{-2}(\boldsymbol{y},\lambda)\pi^{-2}\mathrm{diag}_{k=2}^{\infty}k^{-2}\\
\boldsymbol{I} & \boldsymbol{0}
\end{pmatrix}.
\]
We can now calculate $\boldsymbol{C}$ or its representation $\boldsymbol{C}\boldsymbol{T}\boldsymbol{a}=\boldsymbol{B}_{12}\boldsymbol{B}_{2}^{-1}\boldsymbol{T}\boldsymbol{a}$,
which becomes
\[
\boldsymbol{C}\boldsymbol{T}\boldsymbol{a}=\begin{pmatrix}\begin{array}{l}
0\\
-\Gamma\frac{\gamma(\lambda)\pi^{2}\left(y_{1}+y_{1}^{\star}(\boldsymbol{y},\lambda)\right)\sin\pi\xi}{c^{2}(\boldsymbol{y},\lambda)}\sum_{k=2}^{\infty}\left(-2\beta\frac{a_{1k}}{k\pi}-\frac{a_{2k}}{k\pi}\right)\sin k\pi\xi^{\star}
\end{array}\end{pmatrix}.
\]
The derivative of the immersion $\boldsymbol{W}$, when Fourier expanded
is
\begin{equation}
D_{2}\boldsymbol{W}(\boldsymbol{y},\lambda)=\begin{pmatrix}\begin{array}{l}
{\displaystyle 2\sum_{k=1}^{\infty}\gamma^{\prime}(\lambda)\frac{\sin k\pi\xi^{\star}}{k^{2}\pi^{2}}\sin k\pi\xi+D_{2}y_{1}^{\star}(\boldsymbol{y},\lambda)\sin\pi\xi}\\
{\displaystyle D_{2}y_{2}^{\star}(\boldsymbol{y},\lambda)\sin\pi\xi}
\end{array}\end{pmatrix},\label{eq:NLS-D2W_Four_Unres}
\end{equation}
hence the coordinates of $\boldsymbol{T}^{-1}\boldsymbol{P}D_{2}\boldsymbol{W}(\boldsymbol{y},\lambda)$
and $\boldsymbol{T}^{-1}\boldsymbol{Q}D_{2}\boldsymbol{W}(\boldsymbol{y},\lambda)$
are 
\[
a_{11}=2\gamma^{\prime}(\lambda)\frac{\sin\pi\xi^{\star}}{\pi^{2}}+D_{2}y_{1}^{\star}(\boldsymbol{y},\lambda),\;a_{21}=D_{2}y_{2}^{\star}(\boldsymbol{y},\lambda)
\]
and
\[
a_{1k}=2\gamma^{\prime}(\lambda)\frac{\sin k\pi\xi^{\star}}{k^{2}\pi^{2}},\;a_{2k}=0,\;k\ge2,
\]
respectively. The constraint \nameref{A4BAR-Hyperbolicity}, i.e.,
$\Pi^{c}D_{2}\boldsymbol{W}(\boldsymbol{y},\lambda)=\boldsymbol{0}$
gives
\begin{multline}
\left(\boldsymbol{I},-\boldsymbol{C}\right)\begin{pmatrix}\boldsymbol{P}D_{2}\boldsymbol{W}(\boldsymbol{y},\lambda)\\
\boldsymbol{Q}D_{2}\boldsymbol{W}(\boldsymbol{y},\lambda)
\end{pmatrix}=\begin{pmatrix}\begin{array}{l}
2\gamma^{\prime}(\lambda)\frac{\sin\pi\xi^{\star}}{\pi^{2}}+D_{2}y_{1}^{\star}(\boldsymbol{y},\lambda)\\
D_{2}y_{2}^{\star}(\boldsymbol{y},\lambda)
\end{array}\end{pmatrix}\\
+\begin{pmatrix}\begin{array}{l}
0\\
-\left(\frac{4\Gamma\beta\gamma^{\prime}(\lambda)\gamma(\lambda)}{c^{2}(\boldsymbol{y},\lambda)}\sum_{k=2}^{\infty}\frac{\sin^{2}k\pi\xi^{\star}}{k^{3}\pi}\right)\left(y_{1}+y_{1}^{\star}(\boldsymbol{y},\lambda)\right)
\end{array}\end{pmatrix}=\begin{pmatrix}0\\
0
\end{pmatrix},\label{eq:NLS-Qequ}
\end{multline}
which is an equation for $y_{1}^{\star}(\lambda)$ and $y_{2}^{\star}(y_{1},\lambda)$.
Equation (\ref{eq:NLS-Qequ}) is solved for $D_{2}y_{1}^{\star}$
and $D_{2}y_{2}^{\star}$, which are then integrated over $\lambda$.
The constant of integration for $y_{1}^{\star}$ is such that $y_{1}^{\star}(\boldsymbol{y},0)=0$,
which yields (\ref{eq:NLS y1Star}). However we notice that the constant
of integration does not play a role, so we present the simplest formula
for $y_{2}^{\star}(y_{1},\lambda)$, whose derivative is $D_{2}y_{2}^{\star}(y_{1},\lambda)$
without paying attention to the constant of integration. The result
of this integration is (\ref{eq:NLS-y2Star}). Evaluating the square
of the wave speed with this coordinate shift yields (\ref{eq:NLS-FinalWaveSpeed}).
\end{svmultproof}

\begin{remark}
\label{rem:NLS-Q0-Gamma0}For $\Gamma=0$ we have $\boldsymbol{C}=\boldsymbol{0}$
and also $D_{2}y_{2}^{\star}=0$. This implies that for the linear
string the bundle projection is simply $\boldsymbol{Q}$. The normal
bundle is independent of $\varepsilon$ and there is no need to introduce
time-scale separation. Instead of $\varepsilon$, $\Gamma$ can be
used to track the deformation of the invariant normal bundle, which
persists for a sufficiently small $\Gamma>0$ due to the properties
of $\mathcal{M}_{\lambda}$ \cite{Bates1998}.
\end{remark}

\subsection{The vector field $\boldsymbol{f}(\boldsymbol{y},\lambda)$ on the
invariant manifold}
\begin{lemma}
\label{lem:NLS-SkeletonModel}The skeleton model of equation (\ref{eq:NLS fullVF})
on the invariant manifold specified by (\ref{eq:NLS-FullImmersion})
and with coordinate shifts (\ref{eq:NLS y1Star}) and (\ref{eq:NLS-y2Star})
can be written as
\begin{equation}
\boldsymbol{f}(\boldsymbol{y},\lambda)=\begin{pmatrix}\begin{array}{l}
y_{2}+y_{2}^{\star}(y_{1},\lambda)\\
-c^{2}(\boldsymbol{y},\lambda)\left(\pi^{2}y_{1}-2\gamma(\lambda)\sin\pi\xi^{\star}\right)-\left(2\beta\pi+D_{1}y_{2}^{\star}(y_{1},\lambda)\right)\left(y_{2}+y_{2}^{\star}(y_{1},\lambda)\right)
\end{array}\end{pmatrix}.\label{eq:NLS-skeletonVF}
\end{equation}
After substituting the immersion (\ref{eq:NLS-FullImmersion}), the
switching function (\ref{eq:NLS h Function}) becomes
\begin{equation}
h_{0}(\boldsymbol{y},\lambda)=h(\boldsymbol{W}(\boldsymbol{y},\lambda))=v_{0}-\left(y_{2}+y_{2}^{\star}(y_{1},\lambda)\right)\sin\pi\xi^{\star}.\label{eq:NLS-skeleton-SWfun}
\end{equation}
\end{lemma}
\begin{svmultproof}
The dynamics on the invariant manifold $\mathcal{M}_{\lambda}$ is
given by the invariance condition
\[
D_{1}\boldsymbol{W}(\boldsymbol{y},\lambda)\boldsymbol{f}(\boldsymbol{y},\lambda)=\boldsymbol{F}(\boldsymbol{W}(\boldsymbol{y},\lambda),\lambda).
\]
This is an equation in the tangent bundle of $\mathcal{M}_{\lambda}$,
therefore it makes sense to project it using $\boldsymbol{P}$, as
defined by (\ref{eq:NLS P definition}), to find $\boldsymbol{f}$.
We first calculate that
\begin{equation}
\boldsymbol{P}D_{1}\boldsymbol{W}(\boldsymbol{y},\lambda)=\begin{pmatrix}1 & 0\\
D_{1}y_{2}^{\star}(y_{1},\lambda) & 1
\end{pmatrix}.\label{eq:NLS-P.D1W}
\end{equation}
By inverting the matrix (\ref{eq:NLS-P.D1W}) we find that the reduced
vector field is 
\begin{equation}
\boldsymbol{f}(\boldsymbol{y},\lambda)=\begin{pmatrix}1 & 0\\
-D_{1}y_{2}^{\star}(y_{1},\lambda) & 1
\end{pmatrix}\boldsymbol{P}\boldsymbol{F}(\boldsymbol{W}(\boldsymbol{y},\lambda),\lambda).\label{eq:NLS-ProjectionApplied}
\end{equation}
Next we substitute the immersion (\ref{eq:NLS-immersion}) so that
the vector field (\ref{eq:NLS fullVF}) on the manifold becomes
\begin{equation}
\boldsymbol{F}(\boldsymbol{W}(\boldsymbol{y},\lambda),\lambda)=\begin{pmatrix}\begin{array}{l}
(y_{2}+y_{2}^{\star}(y_{1},\lambda))\sin\pi\xi\\
-\left(y_{1}-2\gamma(\lambda)\frac{\sin\pi\xi^{\star}}{\pi^{2}}\right)c^{2}(\boldsymbol{y},\lambda)\pi^{2}\sin\pi\xi-2\beta\left(y_{2}+y_{2}^{\star}(y_{1},\lambda)\right)\pi\sin\pi\xi
\end{array}\end{pmatrix}.\label{eq:NLS-ImmersionSubstituted}
\end{equation}
Substituting (\ref{eq:NLS-ImmersionSubstituted}) into (\ref{eq:NLS-ProjectionApplied})
yields the reduced vector field (\ref{eq:NLS-skeletonVF}). The switching
function, defined by equation (\ref{eq:NLS h Function}) after substituting
the immersion becomes (\ref{eq:NLS-skeleton-SWfun}).
\end{svmultproof}

\subsection{\label{subsec:NLS-discGap_sigma}The normal discontinuity gap $d^{\pm}$
and decay rate $\sigma$}

The normal discontinuity gap $d^{\pm}$ measures the discontinuity
of the correction about the invariant manifold with initial conditions
$D_{2}\boldsymbol{W}(\boldsymbol{y},\lambda)$ at $t=0$ and determines
the uniqueness of solutions according to theorem \ref{thm:Cont-Uniqueness}.
We calculate $d^{\pm}$ for the $\varepsilon\to0$ limit, when the
dynamics in the normal bundle of $\mathcal{M}_{\lambda}$ becomes
autonomous. Therefore it is sufficient to evaluate the limit $\lim_{t\downarrow0}Dh\cdot\mathrm{e}^{\boldsymbol{A}_{0}(\boldsymbol{y},\lambda)t}D_{2}\boldsymbol{W}(\boldsymbol{y},\lambda)$. 
\begin{lemma}
\label{lem:NLS-DiscontinuityGap}The normal discontinuity gap in the
limit $\varepsilon\to0$ is
\begin{equation}
d^{\pm}(y_{1},\lambda)=\gamma'(\lambda)\frac{c^{2}(y_{1},\lambda)\cos^{-1}\frac{\beta}{c(y_{1},\lambda)}}{\pi\sqrt{c^{2}(y_{1},\lambda)-\beta^{2}}}.\label{eq:NLS-dPlusMinus}
\end{equation}
The rate of decay as defined by (\ref{eq:MR-dichotomy}) is
\begin{equation}
\sigma=-2\pi\beta.\label{eq:NLS-sigma-value}
\end{equation}
\end{lemma}
\begin{svmultproof}
The calculation is carried out using Fourier series, hence we write
the series expansion
\[
D_{2}\boldsymbol{W}(\boldsymbol{y},\lambda)=\begin{pmatrix}\begin{array}{l}
2\gamma^{\prime}(\lambda)\sum_{k=2}^{\infty}\frac{\sin k\pi\xi^{\star}}{k^{2}\pi^{2}}\sin k\pi\xi\\
D_{2}y_{2}^{\star}(y_{1},\lambda)\sin\pi\xi
\end{array}\end{pmatrix},
\]
which is calculated from (\ref{eq:NLS-D2W_Four_Unres}) by substituting
(\ref{eq:NLS y1Star}). Since $D_{2}\boldsymbol{W}(\boldsymbol{y},\lambda)$
is in the invariant normal bundle of the critical manifold, it is
sufficient to restrict the dynamics there. We use the decomposition
of $\mathrm{e}^{\boldsymbol{A}_{0}(\boldsymbol{y},\lambda)t}$ as
given by (\ref{eq:NLS-A0-partition}) to arrive at the expression
\[
\mathrm{e}^{\boldsymbol{A}_{0}(\boldsymbol{y},\lambda)t}D_{2}\boldsymbol{W}(\boldsymbol{y},\lambda)=\mathrm{e}^{\boldsymbol{B}_{2}t}\boldsymbol{Q}D_{2}\boldsymbol{W}(\boldsymbol{y},\lambda)+\boldsymbol{P}D_{2}\boldsymbol{W}(\boldsymbol{y},\lambda)+\boldsymbol{B}_{12}\int_{0}^{t}\mathrm{e}^{\boldsymbol{B}_{2}\tau}\boldsymbol{Q}D_{2}\boldsymbol{W}(\boldsymbol{y},\lambda)\mathrm{d}\tau,
\]
where 
\begin{align*}
\boldsymbol{P}D_{2}\boldsymbol{W}(\boldsymbol{y},\lambda) & =\begin{pmatrix}\begin{array}{l}
0\\
D_{2}y_{2}^{\star}(y_{1}\lambda)\sin\pi\xi
\end{array}\end{pmatrix},\\
\boldsymbol{Q}D_{2}\boldsymbol{W}(\boldsymbol{y},\lambda) & =\begin{pmatrix}\begin{array}{l}
2\gamma^{\prime}(\lambda)\sum_{k=2}^{\infty}\frac{\sin k\pi\xi^{\star}}{k^{2}\pi^{2}}\sin k\pi\xi\\
0
\end{array}\end{pmatrix}.
\end{align*}
The relevant component of the  solution is 
\begin{equation}
\boldsymbol{z}_{2}(t)=D_{2}y_{2}^{\star}(y_{1},\lambda)\sin\pi\xi-2\gamma^{\prime}(\lambda)\sum_{k=2}^{\infty}\mathrm{e}^{-\pi k\beta t}\frac{\pi c^{2}\sin\left(\pi kt\sqrt{c^{2}-\beta^{2}}\right)}{\sqrt{c^{2}-\beta^{2}}}\frac{\sin k\pi\xi^{\star}}{k^{2}\pi^{2}}\sin k\pi\xi.\label{eq:NLS-z2Sol}
\end{equation}
The limit $d^{+}=-\lim_{t\downarrow0}\boldsymbol{z}_{2}(t)\vert_{\xi=\xi^{\star}}$
is calculated as
\begin{equation}
d^{+}(y_{1},\lambda)=\lim_{t\downarrow0}Dh\cdot\boldsymbol{z}(t)=\gamma'(\lambda)\frac{c^{2}\cos^{-1}\nicefrac{\beta}{c}}{\pi\sqrt{c^{2}-\beta^{2}}}-D_{2}y_{2}^{\star}(y_{1},\lambda)\sin\pi\xi^{\star}\label{eq:NLS-dPlus}
\end{equation}
and therefore we have shown (\ref{eq:NLS-dPlusMinus}). The calculation
of (\ref{eq:NLS-dPlus}) involves lengthy algebraic manipulations,
converting the product of exponentials and trigonometric functions
in (\ref{eq:NLS-z2Sol}) into sums of pure exponential expressions,
which yields a sum of series with exponential terms. Each of the sub-series
converge to logarithms of exponential polynomials. It then turns out
that the result has discontinuities due to branch cuts of the complex
logarithm and the limit at the branch cut brings the result. The detailed
calculation (with slightly different notation) can be found in section
II of the Electronic Supplementary Material of \cite{SzalaiMZfriction}. 

The decay rate (\ref{eq:NLS-sigma-value}) is found by reading off
the smallest exponent from formula (\ref{eq:NLS-z2Sol}).
\end{svmultproof}

\begin{remark}
The normal discontinuity gap $d^{\pm}$ is a local property of the
string; it depends on material properties and the tension in the string.
However, $d^{\pm}$ is independent of the boundary conditions and
the position where the string is forced.
\end{remark}

\subsection{Spectrum of the normal dynamics on $\Sigma$}

We use theorem \ref{thm:Normal-hyperbolicity} to find out whether
there exists an attracting critical manifold. 
\begin{lemma}
\label{lem:NLS-CharFunction}The characteristic function determining
the stability of the critical manifold $\text{\ensuremath{\mathcal{M}}}_{\lambda}^{crit}$
within $\Sigma$ is given by
\begin{multline}
\Delta(s)=2\gamma^{\prime}(\lambda)\frac{1}{s}\left(\Gamma\gamma(\lambda)\pi^{2}\left(y_{1}+y_{1}^{\star}(\lambda)\right)\sin\pi\xi^{\star}\right)\times\\
\times\sum_{k=2}^{\infty}\left(\frac{\sin^{2}k\pi\xi^{\star}}{k^{2}\pi^{2}}-\frac{c^{2}(y_{1},\lambda)\sin^{2}k\pi\xi^{\star}}{s^{2}+2s\beta\pi k+c^{2}(y_{1},\lambda)\pi^{2}k^{2}}\right)\\
+2\gamma^{\prime}(\lambda)s\sum_{k=2}^{\infty}\frac{c^{2}(y_{1},\lambda)\sin^{2}k\pi\xi^{\star}}{s^{2}+2s\beta\pi k+c^{2}(y_{1},\lambda)\pi^{2}k^{2}}.\label{eq:NLS-chfun}
\end{multline}
\end{lemma}
\begin{svmultproof}
The characteristic function (\ref{eq:TimeScale-CharFunction}), whose
roots define stability, is formally written as
\[
\Delta(s)=Dh\cdot\left(s\left(s-\boldsymbol{A}_{0}(\boldsymbol{y},\lambda)\right)^{-1}D_{2}\boldsymbol{W}(\boldsymbol{y},\lambda)-D_{2}\boldsymbol{W}(\boldsymbol{y},\lambda)\right).
\]
It is possible to find a convergent series expansion of $\Delta(s)$
by using Fourier series. Let us, for the moment, define $\boldsymbol{x}=\left(s-\boldsymbol{A}_{0}(\boldsymbol{y},\lambda)\right)^{-1}D_{2}\boldsymbol{W}(\boldsymbol{y},\lambda)$,
which is obtained by solving
\begin{equation}
\left(s-\boldsymbol{A}_{0}(\boldsymbol{y},\lambda)\right)\boldsymbol{x}=D_{2}\boldsymbol{W}(\boldsymbol{y},\lambda)\label{eq:NLS-Delta-resolv}
\end{equation}
for $\boldsymbol{x}$. We separate the solution into the first Fourier
coefficient and the rest, such that 
\begin{equation}
\boldsymbol{x}=\left(\left(x_{11},\boldsymbol{x}_{21}\right),\left(x_{12},\boldsymbol{x}_{22}\right)\right),\,\boldsymbol{x}_{2\ell}=\left\{ x_{2\ell,k}\right\} _{k=2}^{\infty},\label{eq:NLS-x-FourierNotation}
\end{equation}
where $x_{11}$ and $x_{12}$ are the coefficients of $\sin\pi\xi$
and $x_{2\ell,k}$ are the coefficients of $\sin k\pi\xi$ in the
Fourier expansion of $\boldsymbol{x}$. Now expanding equation (\ref{eq:NLS-Delta-resolv})
and using the notation (\ref{eq:NLS-x-FourierNotation}) gives
\begin{equation}
s\begin{pmatrix}x_{11}\\
x_{12}
\end{pmatrix}-\begin{pmatrix}\begin{array}{l}
0\\
-\Gamma\left(\gamma(\lambda)\sum_{k=2}^{\infty}x_{21,k}\sin k\pi\xi^{\star}\right)\pi^{2}\left(y_{1}+y_{1}^{\star}(\lambda)\right)\sin\pi\xi
\end{array}\end{pmatrix}=\begin{pmatrix}0\\
D_{2}y_{2}^{\star}(y_{1},\lambda)
\end{pmatrix}\label{eq:NLS-x-Delta-resolv-P}
\end{equation}
and
\begin{equation}
s\begin{pmatrix}x_{21,k}\\
x_{22,k}
\end{pmatrix}-\begin{pmatrix}x_{22,k}\\
-c^{2}(y_{1},\lambda)\pi^{2}k^{2}x_{21,k}-2\beta\pi kx_{22,k}
\end{pmatrix}=\begin{pmatrix}2\gamma^{\prime}(\lambda)\frac{\sin k\pi\xi^{\star}}{k^{2}\pi^{2}}\\
0
\end{pmatrix}.\label{eq:NLS-x-Delta-resolv-Q}
\end{equation}
The solution to equation (\ref{eq:NLS-x-Delta-resolv-Q}) for the
$k\ge2$ Fourier coefficients is
\[
x_{22,k}=-2\gamma^{\prime}(\lambda)\frac{c^{2}(y_{1},\lambda)\sin k\pi\xi^{\star}}{s^{2}+2s\beta\pi k+c^{2}(y_{1},\lambda)\pi^{2}k^{2}}
\]
and
\[
x_{21,k}=\frac{1}{s}2\gamma^{\prime}(\lambda)\left(\frac{\sin k\pi\xi^{\star}}{k^{2}\pi^{2}}-\frac{c^{2}(y_{1},\lambda)\sin k\pi\xi^{\star}}{s^{2}+2s\beta\pi k+c^{2}(y_{1},\lambda)\pi^{2}k^{2}}\right).
\]
For the first Fourier coefficients the solution of equation (\ref{eq:NLS-x-Delta-resolv-P})
is $sx_{11}=0$ and
\begin{multline*}
sx_{12}=D_{2}y_{2}^{\star}(y_{1},\lambda)-\Gamma\frac{1}{s}2\gamma^{\prime}(\lambda)\gamma(\lambda)\pi^{2}\left(y_{1}+y_{1}^{\star}(\lambda)\right)\times\\
\times\sum_{k=2}^{\infty}\left(\frac{1}{k^{2}\pi^{2}}-\frac{c^{2}(y_{1},\lambda)}{s^{2}+2s\beta\pi k+c^{2}(y_{1},\lambda)\pi^{2}k^{2}}\right)\sin^{2}k\pi\xi^{\star}.
\end{multline*}
The series expansion of the characteristic function $\Delta(s)$ using
notation (\ref{eq:NLS-x-FourierNotation}) is 
\[
\Delta(s)=-sx_{12}\sin\pi\xi^{\star}-s\sum_{k=2}^{\infty}x_{22,k}\sin k\pi\xi^{\star}+D_{2}y_{2}^{\star}(y_{1},\lambda)\sin\pi\xi^{\star}
\]
and substituting system parameters yields (\ref{eq:NLS-chfun}).
\end{svmultproof}

\begin{figure}
\begin{centering}
\includegraphics[width=0.5\linewidth]{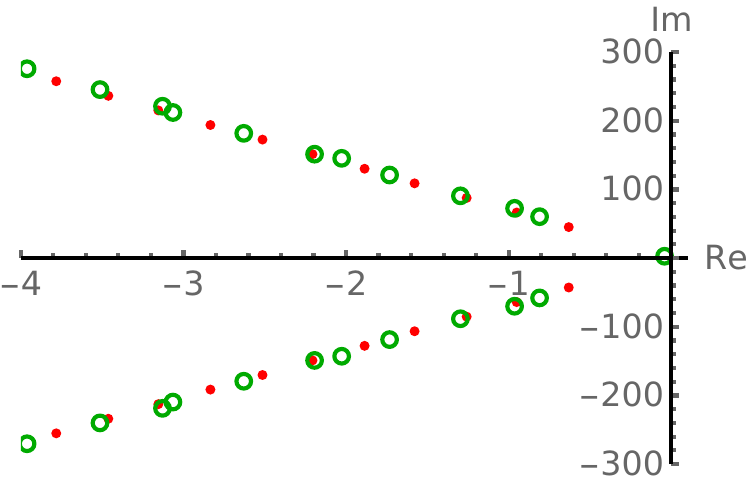}\includegraphics[width=0.5\linewidth]{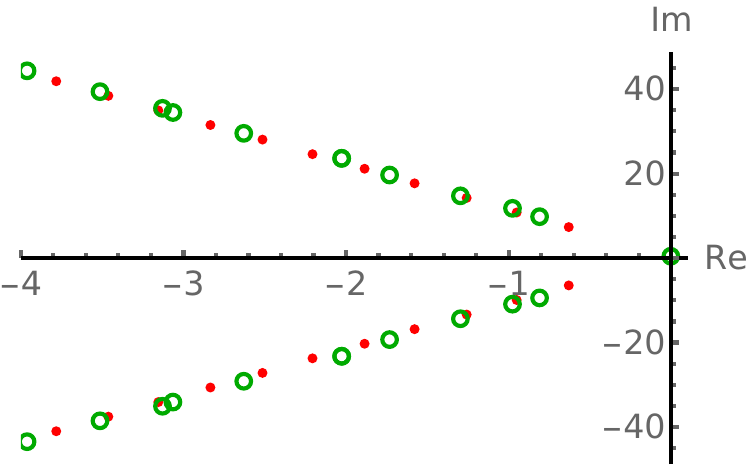}
\par\end{centering}
\caption{\label{fig:NLS-roots} Green circles show the roots of the characteristic
function (\ref{eq:NLS-chfun}) at $\lambda=1$ and $y_{1}=1$, with
parameters $\mu=1$, $\beta=0.1$ and $\xi^{\star}=\sqrt{2}/2$. The
red dots show the eigenvalues of $\boldsymbol{A}_{0}(\boldsymbol{y},\lambda)$
without the zeros of the tangent bundle for comparison. Left: The
nonlinear case with $\Gamma=20$; Right: linear case with $\Gamma=0$.
Note that the natural frequencies are much higher for the nonlinear
case, because the string has significantly more tension due to deformation.}
\end{figure}
 Figure \ref{fig:NLS-roots} shows the roots of (\ref{eq:NLS-chfun})
at an attracting point of the critical manifold. It can be seen that
there is a real root near zero, while other roots are well inside
the left complex half space. It seems that roots of $\Delta(s)$ are
perturbations of the eigenvalues of $\boldsymbol{A}_{0}(\boldsymbol{y},\lambda)$
apart from the rightmost root, that appears due to switching.

\begin{figure}
\begin{centering}
\includegraphics[width=0.7\linewidth]{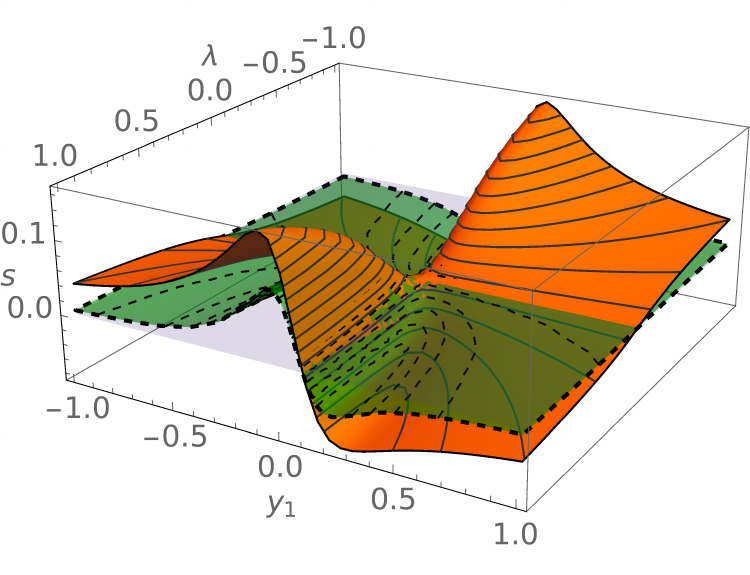}
\par\end{centering}
\caption{\label{fig:NLS-charRoot}The rightmost root of the characteristic
function (\ref{eq:NLS-chfun}) (orange, solid lines) and $3\times(-\sigma d^{-}/d^{\pm})$
in (\ref{eq:Conc-CritManifDefinition}) (green, dashed lines) determine
whether the critical manifold is attracting at any given value of
$\boldsymbol{y}$ and $\lambda$. Parameters are $\Gamma=10$, $\mu=1$,
$\beta=0.1$ and $\xi^{\star}=\sqrt{2}/2$.}
\end{figure}

The plot of this rightmost root is shown in figure \ref{fig:NLS-charRoot}
in orange (solid lines), which indicates that the critical manifold
is partly attracting (negative root), partly repelling (positive root).
This might be surprising because the system dissipates energy as a
whole. However, the constraint $h=0$ and nonlinearity couple the
dynamics in the tangent and normal bundles and energy is exchanged
between them causing instability. In contrast, there is no such coupling
in the linear string ($\Gamma=0$), the green root near the origin
in figure \ref{fig:NLS-roots} remains at the origin and therefore
the normal dynamics is neutrally stable.
\begin{remark}
\label{rem:NLS-EpsOneJustification}Continuing from remark \ref{rem:NLS-Q0-Gamma0}
we find that for $\Gamma=0$ and for all $\varepsilon\in[0,1]$ the
characteristic function (\ref{eq:NLS-chfun}) is valid due to $\boldsymbol{A}_{1}$
being constant. Let us denote the zero root of the characteristic
function (\ref{eq:NLS-chfun}) for $\Gamma=0$ by $\sigma^{0}$. The
invariant vector bundle corresponding to $\sigma^{0}$ is then isomorphic
to $\mathcal{M}_{\lambda}\times\mathbb{R}$. For $\Gamma>0$ the invariant
vector bundle of $\sigma^{0}$ is continuously perturbed. The perturbation
turns $\sigma^{0}$ into a small Lyapunov exponent of the now non-autonomous
dynamics within the invariant vector bundle. $\Gamma>0$ can be chosen
such that $\sigma_{s}<-\varepsilon\left|\sigma^{0}\right|$, that
is, the invariant vector bundle is an attracting normally hyperbolic
invariant manifold in $\left(\mathcal{M}_{\lambda}\times\boldsymbol{X}\right)\cap\Sigma$,
that is the phase space of the corrected model in $\Sigma$. The dynamics
in the invariant vector bundle is represented by the reduced order
model on $\mathcal{M}_{\lambda}\times\mathbb{R}$.
\end{remark}

\subsection{Equivalent reduced order model on $\Sigma_{\varepsilon}$}

We now investigate the dynamics of the string on $\Sigma_{\varepsilon}$.
The dynamics outside $\Sigma_{\varepsilon}$ is given by $\dot{\boldsymbol{y}}=\boldsymbol{f}(\boldsymbol{y},\lambda)$,
$\varepsilon\dot{\kappa}=\sigma\kappa$ and $\dot{\lambda}=0$. The
skeleton model on $\Sigma_{\varepsilon}$ is formally given by equation
(\ref{eq:RED-SwitchingEq}), while the reduced order model including
a qualitative approximation of the normal dynamics is given by (\ref{eq:Conc-Sigma-dynamics}).
The complication with equations (\ref{eq:RED-SwitchingEq}) and (\ref{eq:Conc-Sigma-dynamics})
is that they involve the lengthy term $y_{2}^{\star}(y_{1},\lambda)$
as shown by equation (\ref{eq:NLS-y2Star}). It is possible to eliminate
$y_{2}^{\star}(y_{1},\lambda)$ using the transformation
\begin{equation}
\overline{\boldsymbol{y}}=\left(y_{1},y_{2}+y_{2}^{\star}(y_{1},\lambda)\right)^{T}.\label{eq:NLS-SimplyTrafo}
\end{equation}
The vector field on the invariant manifold now involves $\dot{\lambda}$
in the form of $\dot{\overline{\boldsymbol{y}}}=\overline{\boldsymbol{f}}(\overline{\boldsymbol{y}},\lambda,\dot{\lambda})$,
where
\begin{equation}
\overline{\boldsymbol{f}}(\overline{\boldsymbol{y}},\lambda,\dot{\lambda})=\begin{pmatrix}\begin{array}{l}
\overline{y}_{2}\\
-\left(c^{2}(\overline{y}_{1},\lambda)\left(\pi^{2}\overline{y}_{1}-2\gamma(\lambda)\sin\pi\xi^{\star}\right)+2\beta\pi\overline{y}_{2}\right)+D_{2}y_{2}^{\star}(\overline{y}_{1},\lambda)\dot{\lambda}
\end{array}\end{pmatrix}.\label{eq:NLS-SimplifiedVF}
\end{equation}
The main advantage of this formulation is that the function defining
the switching manifold $\Sigma$ becomes simpler, namely 
\begin{align}
\overline{h}_{\varepsilon}(\overline{\boldsymbol{y}},\kappa,\lambda) & =v_{0}-\overline{y}_{2}\sin\pi\xi^{\star}+\kappa,\label{eq:NLS-Simplified-hEps}
\end{align}
which is independent of $\lambda$. The stick dynamics on $\Sigma$
has the same dependence on $\dot{\lambda}$ as before, because 
\begin{multline}
\frac{\mathrm{d}}{\mathrm{d}t}\overline{h}_{\varepsilon}(\overline{\boldsymbol{y}},\kappa,\lambda)=D_{1}h_{\varepsilon}(\overline{\boldsymbol{y}},\kappa,\lambda)\overline{\boldsymbol{f}}(\overline{\boldsymbol{y}},\lambda,\dot{\lambda})+\dot{\kappa}\\
=\left(\left(c^{2}(\overline{y}_{1},\lambda)\left(\pi^{2}\overline{y}_{1}-2\gamma(\lambda)\sin\pi\xi^{\star}\right)+2\beta\pi\overline{y}_{2}\right)-D_{2}y_{2}^{\star}(\overline{y}_{1},\lambda)\dot{\lambda}\right)\sin\pi\xi^{\star}+\dot{\kappa}.\label{eq:NLS-Simplified-hEpsDeri}
\end{multline}
Note that $d^{-}=D_{2}h_{0}(\overline{\boldsymbol{y}},\lambda)=-D_{2}y_{2}^{\star}(\overline{y}_{1},\lambda)$
appears in equation (\ref{eq:NLS-Simplified-hEpsDeri}), and remains
the coefficient of $\dot{\lambda}$. The last equation we need is
(\ref{eq:Conc-kappa-DE}) that describes $\kappa$. Note that the
transformation (\ref{eq:NLS-SimplyTrafo}) does not change the values
of $d^{+}$, $d^{-}$ and $d^{\pm}$ given by (\ref{eq:NLS-dPlus})
and (\ref{eq:NLS-dPlusMinus}), because they do not depend on $y_{2}$.
As a result we have 
\begin{equation}
\left.\begin{array}{rl}
\dot{\overline{y}_{1}} & =\overline{y}_{2}\\
\dot{\overline{y}_{2}} & =\dot{\kappa}/\sin\pi\xi^{\star}\\
\dot{\lambda} & =\frac{c^{2}(\overline{y}_{1},\lambda)\left(\pi^{2}\overline{y}_{1}-2\gamma(\lambda)\sin\pi\xi^{\star}\right)+2\beta\pi\overline{y}_{2}+\varepsilon^{-1}\sigma\kappa}{d^{\pm}(\overline{y}_{1},\lambda)}\\
\dot{\kappa} & =-\frac{d^{+}(\overline{y}_{1},\lambda)}{d^{\pm}(\overline{y}_{1},\lambda)}\left(c^{2}(\overline{y}_{1},\lambda)\left(\pi^{2}\overline{y}_{1}-2\gamma(\lambda)\sin\pi\xi^{\star}\right)+2\beta\pi\overline{y}_{2}\right)-\varepsilon^{-1}\frac{d^{-}(\overline{y}_{1},\lambda)}{d^{\pm}(\overline{y}_{1},\lambda)}\sigma\kappa
\end{array}\right\} .\label{eq:NLS-ReducedModel-atStick}
\end{equation}

\subsection{Dynamics of the skeleton model on $\Sigma_{0}$}

In this section we explore the dynamics of the skeleton model, which
is the same as the dynamics on the critical manifold, when $\varepsilon=0$
in equation (\ref{eq:Conc-Sigma-dynamics}). The dynamics on the critical
manifold can be found by setting $\overline{y}_{2}=v_{0}/\sin\pi\xi^{\star}$and
$\dot{\overline{y}_{2}}=0$ so that $h=0$ and $\dot{h}=0$, then
solving $\dot{\overline{\boldsymbol{y}}}=\overline{\boldsymbol{f}}(\overline{\boldsymbol{y}},\lambda,\dot{\lambda})$
as an algebraic equation with (\ref{eq:NLS-SimplifiedVF}) on the
right side for $\dot{\lambda}$, that is,
\begin{equation}
\left.\begin{array}{rl}
\dot{\overline{y}}_{1} & =v_{0}/\sin\pi\xi^{\star}\\
\dot{\lambda} & ={\displaystyle \frac{c^{2}(\overline{y}_{1},\lambda)\left(\pi^{2}\overline{y}_{1}-2\gamma(\lambda)\sin\pi\xi^{\star}\right)+2\beta\pi v_{0}/\sin\pi\xi^{\star}}{D_{2}y_{2}^{\star}(\overline{y}_{1},\lambda)}}
\end{array}\right\} .\label{eq:NLS-CritStick}
\end{equation}
As we noted in theorem \ref{thm:ExtUnique} in section \ref{sec:skeleton},
solutions of this model pass through the boundaries $\Sigma^{\pm}$
if $D_{2}y_{2}^{\star}(\overline{y}_{1},\lambda)\sin\pi\xi^{\star}=-d^{-}>0$.
To avoid any problem with having the wrong sign of $d^{-}$ we re-scale
time by $D_{2}y_{2}^{\star}(\overline{y}_{1},\lambda)$ and get
\begin{equation}
\left.\begin{array}{rl}
\dot{\overline{y}}_{1} & =D_{2}y_{2}^{\star}(\overline{y}_{1},\lambda)v_{0}/\sin\pi\xi^{\star}\\
\dot{\lambda} & ={\displaystyle c^{2}(\overline{y}_{1},\lambda)\left(\pi^{2}\overline{y}_{1}-2\gamma(\lambda)\sin\pi\xi^{\star}\right)+2\beta\pi v_{0}/\sin\pi\xi^{\star}}
\end{array}\right\} ,\label{eq:NLS-CritStick-DEN}
\end{equation}
whose forward-time solutions always pass through $\Sigma^{\pm}$.
This allows for a straightforward numerical solution.
\begin{figure}
\begin{centering}
\includegraphics[width=0.49\linewidth]{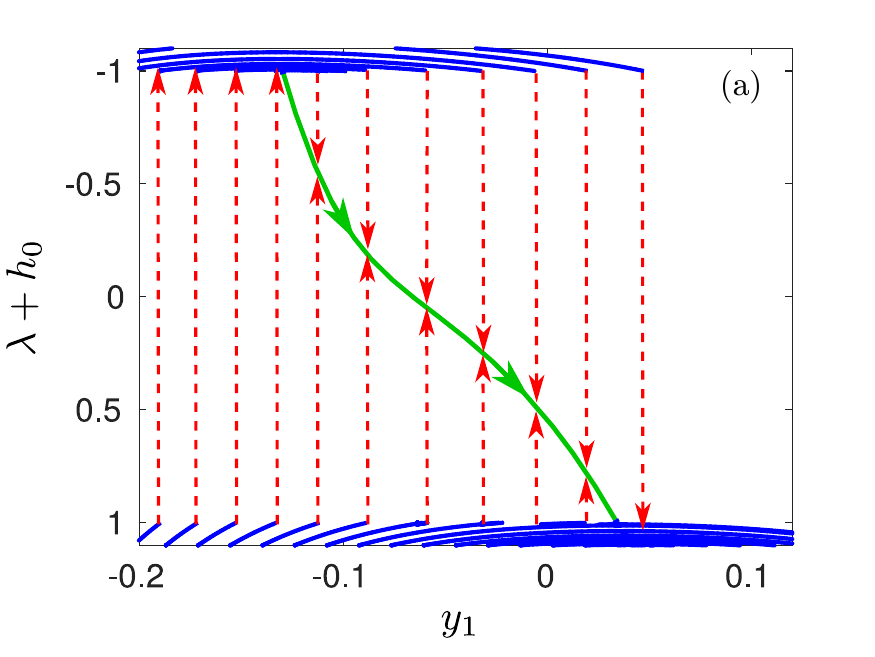}\includegraphics[width=0.49\linewidth]{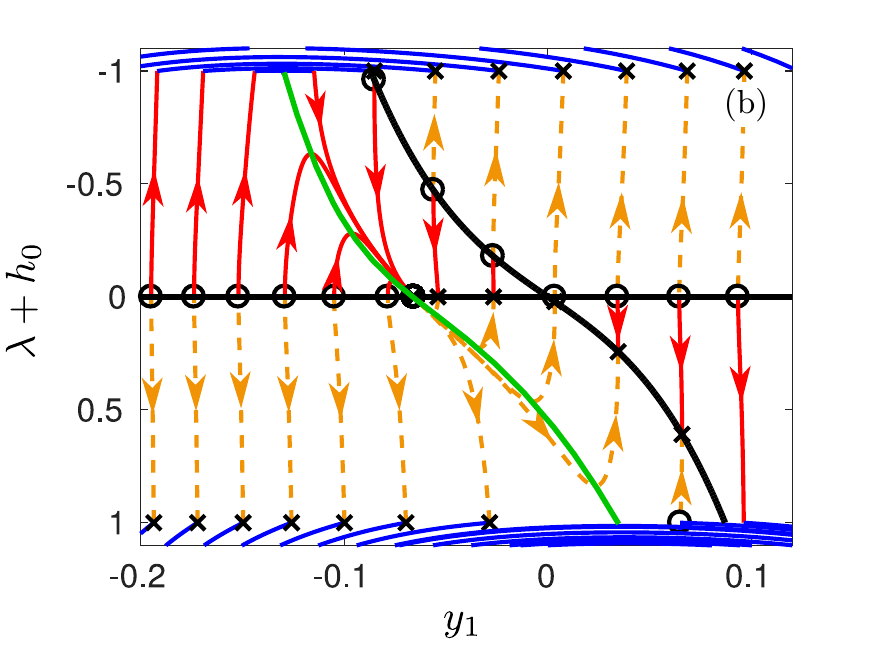}
\par\end{centering}
\caption{\label{fig:NLS-Crit-Phase}Phase portrait of the skeleton model. (a)
Utkin's closure, assuming $D_{2}y_{2}^{\star}(\overline{y}_{1},\lambda)=0$.
The dashed red lines show how the solution jumps between values of
$\lambda$. Blue lines are trajectories with $\dot{\lambda}=0$. (b)
The red lines are solutions with $D_{2}y_{2}^{\star}(\overline{y}_{1},\lambda)>0$,
so that the red and blue lines are a continuation of each other. The
dashed orange lines represent solutions with $D_{2}y_{2}^{\star}(\overline{y}_{1},\lambda)<0$.
At points marked by black crosses the solutions come together infinitely
fast, when $\lambda\protect\neq\pm1$. The points marked by circles
repel solutions infinitely fast, when $\lambda\protect\neq\pm1$.
The black lines correspond to $D_{2}y_{2}^{\star}(\overline{y}_{1},\lambda)=0$
and the thick green line within $\Sigma$ is the nullcline $\dot{\lambda}=0$.
The parameters are $\beta=0.1$, $\Gamma=20$, $\xi^{\star}=\sqrt{2}/2$.}
\end{figure}

Let us first recall, what Utkin's closure would produce if we disregarded
$D_{2}y_{2}^{\star}(\overline{y}_{1},\lambda)$. The solution would
be given by the equation
\begin{equation}
\left.\begin{array}{rl}
\dot{\overline{y}}_{1} & =v_{0}/\sin\pi\xi^{\star}\\
0 & =c^{2}(\overline{y}_{1},\lambda)\left(\pi^{2}\overline{y}_{1}-2\gamma(\lambda)\sin\pi\xi^{\star}\right)+2\beta\pi v_{0}/\sin\pi\xi^{\star}
\end{array}\right\} ,\label{eq:NLS-UtkinStick}
\end{equation}
which is partly algebraic. Figure \ref{fig:NLS-Crit-Phase}(a) shows
the phase portrait. The dashed orange lines correspond to $\lambda$
values jumping between either $\pm1$ or the solution of the algebraic
constraint in equation (\ref{eq:NLS-UtkinStick}). The continuous
green line represents $\lambda$ values that are admissible by the
constraint in equation (\ref{eq:NLS-UtkinStick}). Since $\dot{\overline{y}}_{1}$
is a positive constant, solutions can only move in one direction along
the green line. This is typical of friction oscillators and it is
consistent with rigid body dynamics, where forces are allowed to be
discontinuous.

A different picture emerges when $D_{2}y_{2}^{\star}(\overline{y}_{1},\lambda)$
is considered. Figure \ref{fig:NLS-Crit-Phase}(b) shows a typical
phase portrait of (\ref{eq:NLS-CritStick}). Parts of trajectories
are denoted by dashed lines when $D_{2}y_{2}^{\star}(\overline{y}_{1},\lambda)<0$.
The arrows indicate the correct forward direction of time. Black lines
indicate when $D_{2}y_{2}^{\star}(\overline{y}_{1},\lambda)=0$ and
therefore equation (\ref{eq:NLS-CritStick}) is singular and the direction
of time changes in equation (\ref{eq:NLS-CritStick-DEN}). The black
lines also form a set of nullclines of equation (\ref{eq:NLS-CritStick-DEN}),
because at these points $\dot{\overline{y}}_{1}=0$. Another nullcline
is shown in green, where $\dot{\lambda}=0$. At the intersection of
the green and black lines equation (\ref{eq:NLS-CritStick-DEN}) has
an equilibrium, which is a node. The weak stable manifold of this
equilibrium is close to the green nullcline of $\dot{\lambda}=0$. 

The phase portrait of figure \ref{fig:NLS-Crit-Phase}(b) is not typical
for a friction oscillator. Yet, the skeleton model is obtained through
a careful reduction to an invariant manifold, where we made sure any
perturbation due to the discontinuities would only affect the invariant
normal bundle. Applying Filippov's closure at the boundaries $\lambda=\pm1$
yields sliding solutions. However this implies that solutions coming
from either side of $\Sigma$ cannot enter $\Sigma$ while $\lambda$
stays at $\pm1$. Having a fixed value of $\lambda$ is not physical
in a friction oscillator. The singularities within $\Sigma$ are reached
infinitely fast. This resembles the dynamics of the van der Pol oscillator
at the fold point of its critical manifold \cite{Kanamaru:2007} or
in general the dynamics of singularly perturbed systems \cite{kuehn2015multiple}.
For example, the equilibrium of equation (\ref{eq:NLS-CritStick-DEN}),
formed by the intersection of two nullclines resembles folded-node
singularities \cite{Wechselberger2012,Kristiansen2017Flat}. Therefore,
our best chance to gain more insight is to consider the reduced order
model (\ref{eq:NLS-CritStick}) that includes a representation of
the dynamics in the invariant normal bundle of $\mathcal{M}_{\lambda}$
as described in section \ref{subsec:QualitativeModel}.

\subsection{Dynamics of the reduced order model on $\Sigma_{\varepsilon}$}

\begin{figure}
\begin{centering}
\includegraphics[width=0.9\linewidth]{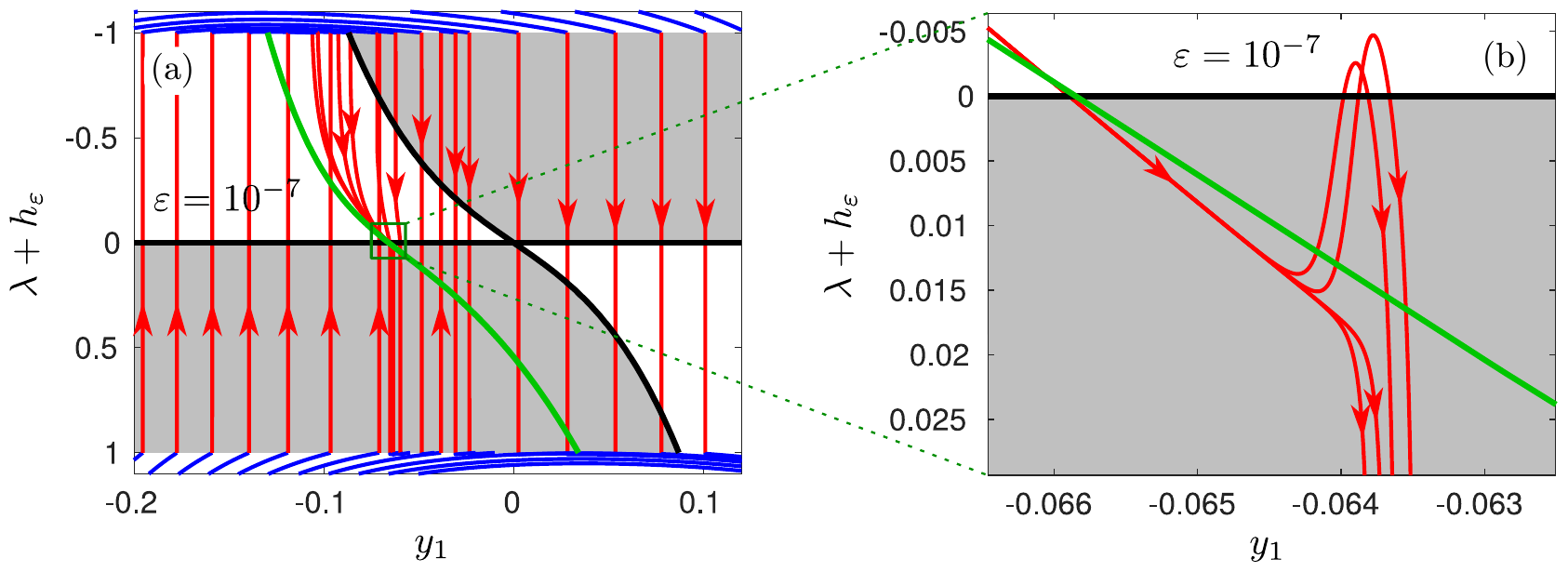}\\
\includegraphics[width=0.9\linewidth]{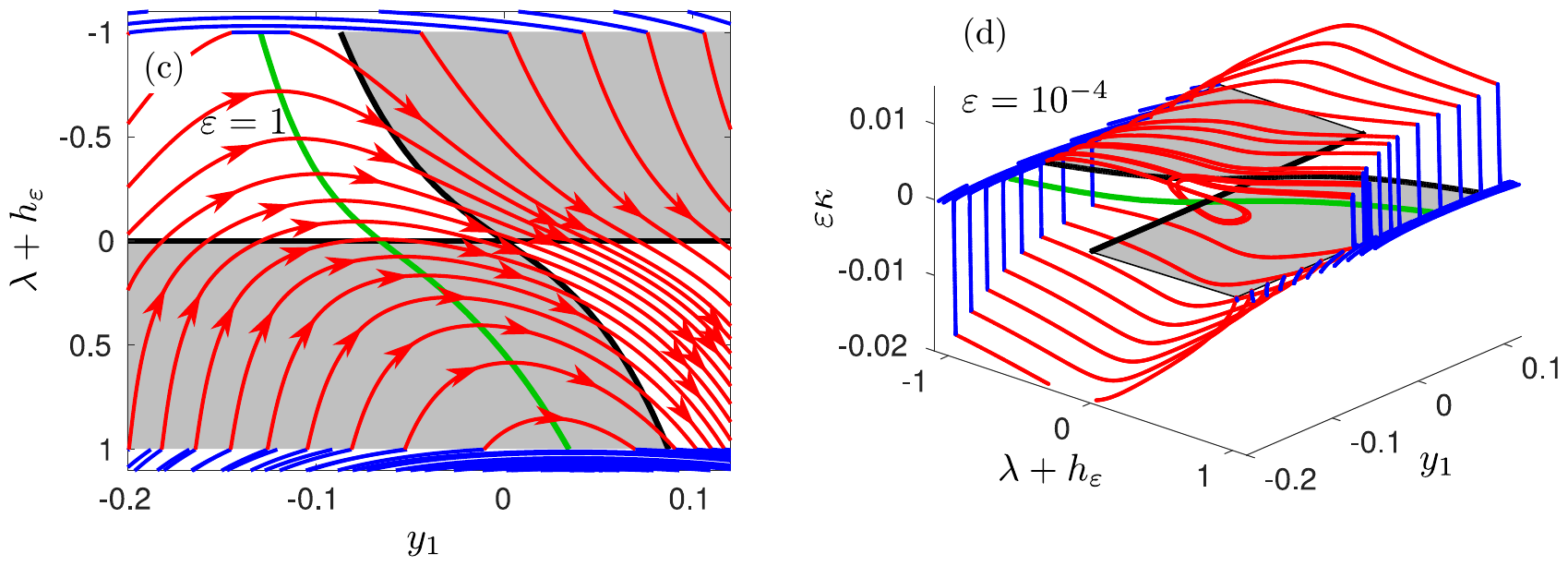}
\par\end{centering}
\caption{\label{fig:NLS-ConcSim}Projections of phase portraits of the reduced
order model (\ref{eq:NLS-ReducedModel-atStick}) in $\Sigma$ for
the nonlinear string. The repelling regions of the critical manifold
are shaded. The black lines correspond to the boundaries, where normal
hyperbolicity of the critical manifold is lost. The thick green line
in $\Sigma$ is the nullcline where $\dot{\overline{y}}_{2}=0$. Parameters
are $\beta=0.1$, $\Gamma=20$, $\xi^{\star}=\sqrt{2}/2$. a) slightly
perturbed dynamics on the critical manifold with $\varepsilon=10^{-7}$;
b) a detailed view of the dynamics about the stable node with $\varepsilon=10^{-7}$;
c) dynamics in $\Sigma_{\varepsilon}$ with the original time-scale
($\varepsilon=1$); d) three dimensional view of the dynamics in $\Sigma_{\varepsilon}$
with $\varepsilon=10^{-4}$, note the fast decay of $\kappa$ outside
of $\Sigma$.}
\end{figure}

In this section we investigate the reduced order model (\ref{eq:Conc-FullSystem}),
which is the extension of the skeleton model by a single variable
representing the dynamics in the normal bundle of $\mathcal{M}_{\lambda}$.
The dynamics on $\Sigma_{\varepsilon}$ is given by equation (\ref{eq:NLS-ReducedModel-atStick})
with parameters derived in section \ref{subsec:NLS-discGap_sigma}.
Proposition \ref{prop:same-stability} shows that the reduced order
model captures the stability of $\mathcal{M}_{\lambda}$ for $\varepsilon=0$
well. Figure \ref{fig:NLS-charRoot} confirms this: in the illustrated
part of the phase space the stability of the critical manifold of
the corrected model and the reduced order model is the same. The critical
manifold is repelling where $d^{-}>0$ as per proposition \ref{prop:Conc-Hyperbolicity}.
The skeleton model does not capture the repelled trajectories and
also displays singular dynamics for $d^{-}>0$ as shown in figure
\ref{fig:NLS-Crit-Phase}(b). Here we illustrate that the positive
value of $d^{\pm}$ for the reduced order model resolves the singularities
that occur in the skeleton model according to proposition \ref{prop:Conc-Uniqueness}.

We first choose a small parameter value $\varepsilon=10^{-7}$ to
show the qualitative differences between equation (\ref{eq:NLS-CritStick})
and equation (\ref{eq:NLS-ReducedModel-atStick}). Figure \ref{fig:NLS-ConcSim}(a)
shows the two-dimensional projection of the phase portrait ignoring
variable $\kappa$. When trajectories start in the shaded part with
$\lambda=1$, where the critical manifold is repelling, they quickly
pass to $\lambda=-1$ without much change in $y_{1}$, while $\kappa$
exponentially explodes. Trajectories starting with $\lambda=-1$,
in the region where the critical manifold is attracting, follow the
manifold while being attracted to the stable node of (\ref{eq:NLS-CritStick})
at the intersection of the green and black lines. At the node, the
stability of the critical manifold changes and trajectories are again
repelled with growing magnitude of $\kappa$. This is illustrated
in figure \ref{fig:NLS-ConcSim}(b). After passing the node, trajectories
tend to either $\lambda=\pm1$. It is then likely that trajectories
will start a violent oscillation between $\lambda=\pm1$, because
they interact with the two repelling parts of the critical manifold.
This dynamics has some resemblance to figure \ref{fig:NLS-Crit-Phase}(b)
except that there is no need to re-scale time, since there is no division
by $d^{-}$.

Increasing $\varepsilon$ leads to less violent oscillations between
$\lambda=\pm1$, which eventually continues without reaching $\lambda=\pm1$.
Such a case is shown in figure \ref{fig:NLS-ConcSim}(d), where the
oscillation is reduced to a single loop about the line where the critical
manifold becomes repelling. For $\varepsilon=1$ the dynamics becomes
relatively slow for all variables and resembles that of typical friction
oscillators with well defined stick and slip phases. This phase portrait
is shown in figure \ref{fig:NLS-ConcSim}(c). For $\varepsilon$ sufficiently
large the time scale of the normal dynamics ($\kappa$ variable) becomes
much longer than the dynamics of the rest of the variables and therefore
during a stick phase $\kappa$ does not change much, which also means
that the instability of the critical manifold loses its influence
on the dynamics. Indeed, the leading characteristic root of $\Delta(s)$
is a small perturbation of the zero root, hence it is easily dominated
by other time-scales. In fact by removing nonlinearity ($\Gamma=0$)
this root remains zero, hence $\kappa$ simply becomes an integral
of other quantities without a dynamics of it own. In our example at
$\varepsilon=1$, $\kappa$ is almost without its own dynamics. The
justification why $\varepsilon$ can be increased to $\varepsilon=1$
can be found in remark \ref{rem:NLS-EpsOneJustification}.

The conclusion from the analysis is that simply applying reduction
to an invariant manifold is not sufficient, one needs to take into
account at least a qualitative approximation of the normal dynamics.
This is because the skeleton model (\ref{eq:RED-skeleton}) over-emphasizes
instabilities and turns them into singularities. The main component
that makes the reduced order model (\ref{eq:Conc-FullSystem}) well
behaved is that $d^{\pm}$ is positive in all parts of the phase space.
For the nonlinear string example, $d^{\pm}$ is the velocity jump
of the contact point due to a unit jump in $\lambda$, i.e., the contact
force. Therefore in light of Newton's second law it is understandable
why $d^{\pm}>0$. In case we had found $d^{\pm}=0$ the reduced order
model (\ref{eq:Conc-FullSystem}), including an approximation of the
normal dynamics about $\mathcal{M}_{\lambda}$, would not be necessary,
the skeleton model would be sufficient.

\section{Conclusion}

In this paper we have investigated PWS systems on Banach spaces with
non-dense domain of definition. Specific application areas that satisfy
this assumption are the elastodynamics equations \cite{Kausel2006},
delay equations \cite{Diekmann1995} or age-dependent population dynamics
\cite{metz1986dynamics}. Such systems are different from other classes
of PWS systems, because they can have unique solutions under general
conditions. We were also able to construct a finite-dimensional reduced
order model that inherits key properties of an infinite dimensional
model. Non-dense domain of definition can arise if the phase space
is non-reflexive, for example when the phase space consists of continuous,
bounded or Lipschitz continuous functions. In some cases boundary
conditions can make the domain non-dense \cite{Neubrander1988}. 

The key quantity that decides uniqueness of solutions is the normal
discontinuity gap, which is due to discontinuous trajectories that
systems with non-dense domains have. For the linear and nonlinear
string the normal discontinuity gap represents the velocity jump of
a contact point in response to a unit jump in force. The presence
of the normal discontinuity gap allows the dynamics inside the switching
manifold to become smooth. As a result, two new discontinuity boundaries
arise, where trajectories can enter or leave the switching manifold.
If the normal discontinuity gap is positive, trajectories cross the
new discontinuity boundaries under general conditions.

Despite uniqueness of solutions, invariant manifolds that extend over
the switching manifold do not exist. We have assumed the existence
of an invariant manifold when the switching parameter of the vector
field is constant. This invariant manifold does not persist when the
switching parameter varies, but we have found that pieces of this
manifold do persist, while discontinuities between the persisting
pieces develop along the two new discontinuity boundaries. We have
also shown that switching can make the invariant manifold repelling.
However in the example of the nonlinear string the invariant manifold
is repelling only in a single direction, which can be captured by
a scalar variable. We have constructed a reduced order model that
captures this instability. It remains to be shown under what conditions
there is a spectral gap between the reduced model and the rest of
the dynamics, so that the reduced order model captures all the essential
dynamics. We have only shown that the invariant manifold becomes repelling
within the reduced order model and within the infinite dimensional
system under the same conditions through a real root (see proposition
\ref{prop:same-stability}).

While the theory presented is incomplete, we hope that the results
in this paper will find applications in simulating PWS continuum systems.
Using the reduced order model instead of the skeleton model eliminates
singularities and allows for a unique solution. This allows well conditioned
numerical schemes that lead to robust solutions unlike what is currently
possible \cite{KANE19991}. While it is not proven that the reduced
order model fully captures all dynamics, we expect that this will
be shown in the future either in general or under further conditions.

We have demonstrated the model reduction procedure on a bowed nonlinear
string model. In this example we have found that the skeleton model
has nonphysical singularities, where the friction force between the
bow and the string remains at its maximal limit. The skeleton model
also has a singularity that resembles a folded node of singularly
perturbed systems \cite{Wechselberger2005,Kristiansen2017Flat}. After
correcting the skeleton model with the dynamics that arises in the
eliminated parts of the system due to switching, the pictures becomes
clearer. It turns out that the correction is a largely decaying motion
with the possibility of an instability along a one dimensional subspace.
When this possible instability is taken into account, the model becomes
free of singularities and the dynamics resembles what a friction oscillator
would exhibit when the friction force is regularized \cite{StomayorTeixeira}.
\begin{acknowledgements}
The author would like to thank Alan R.~Champneys and S.~John Hogan
for feedback on the manuscript. He would also like to thank the anonymous
reviewers who have helped with the clarity of the text and the accuracy
of calculations. The author is especially thankful to Galit Szalai,
who has proofread the final draft.
\end{acknowledgements}

\appendix

\section{\label{sec:Trivial Pert Solution}Solution of the correction term
in the introductory example}

This appendix details the solution of equation

\begin{equation}
\ddot{w}(\xi,t)=w^{\prime\prime}(\xi,t)-\ddot{\lambda}\left(u_{0}(\xi)-y^{\star}\frac{\sin\pi\xi}{\sin\pi\xi^{\star}}\right),\label{eq:Trivial Pert Wave-1}
\end{equation}
with initial condition
\begin{equation}
w(\xi,0)=\dot{w}(\xi,0)=0.\label{eq:Trivial Pert Wave IC-1}
\end{equation}
and boundary conditions $w(0,t)=w(1,t)=0$. The solution of (\ref{eq:Trivial Pert Wave-1})
is then substituted into the switching function 
\begin{equation}
h=v_{0}-\dot{y}-\dot{\lambda}\left(u_{0}(\xi^{\star})-y^{\star}\right)-\dot{w}(\xi^{\star},t),\label{eq:Trivial SWcond Pert-1}
\end{equation}
which replaces $h$ in equation (\ref{eq:Trivial 1DOF}) of section
2. We also assume that the solution starts with $\dot{\lambda}(0)=\ddot{\lambda}(0)=0$,
which occurs for example, when $h\vert_{t=0}\neq0$. 

The solution of equations (\ref{eq:Trivial Pert Wave-1}) and (\ref{eq:Trivial Pert Wave IC-1})
can be expressed using the variation-of-constants formula. Assume
that $w_{h}(\xi,t)$ is the solution of the homogeneous equation $\ddot{w}_{h}(\xi,t)=w_{h}^{\prime\prime}(\xi,t)$
with initial and boundary conditions, as in
\begin{equation}
w_{h}(\xi,0)=0,\;\dot{w}_{h}(\xi,0)=y^{\star}\frac{\sin\pi\xi}{\sin\pi\xi^{\star}}-u_{0}(\xi),\;w_{h}(0,t)=w_{h}(1,t)=0,\label{eq:Trivial HomWave IC}
\end{equation}
then the solution of (\ref{eq:Trivial Pert Wave-1}) for the velocity
with zero initial condition (\ref{eq:Trivial Pert Wave IC-1}) is
\begin{equation}
\dot{w}(\xi,t)=\int_{0}^{t}\dot{w}_{h}(\xi,t-\tau)\ddot{\lambda}(\tau)\mathrm{d}\tau.\label{eq:Trivial Convolution}
\end{equation}
Note that we express the velocity here, because that is what appears
in the switching function (\ref{eq:Trivial SWcond Pert-1}). Using
integration by parts twice transforms equation (\ref{eq:Trivial Convolution})
into
\begin{align}
\dot{w}(\xi,t) & =\dot{w}_{h}(\xi,0)\dot{\lambda}(t)-\dot{w}_{h}(\xi,t)\dot{\lambda}(0)+\ddot{w}_{h}(\xi,0)\lambda(t)-\ddot{w}_{h}(\xi,t)\lambda(0)\nonumber \\
 & \quad+\int_{0}^{t}\dddot{w}_{h}(\xi,t-\tau)\lambda(\tau)\mathrm{d}\tau.\label{eq:Trivial ByParts Integrated}
\end{align}
We now evaluate equation (\ref{eq:Trivial ByParts Integrated}) at
$\xi=\xi^{\star}$ and notice that a number of terms vanish. From
equation (\ref{eq:Trivial Pert Wave-1}), $\ddot{\lambda}(0)=0$ and
$w_{h}(\xi,0)=0$ it follows that $\ddot{w}_{h}(\xi^{\star},0)=0$,
and from equation (\ref{eq:Trivial HomWave IC}) we obtain $\dot{w}_{h}(\xi^{\star},0)=y^{\star}-u_{0}(\xi^{\star})$.
Further, using $\dot{\lambda}(0)=0$ brings (\ref{eq:Trivial ByParts Integrated})
into

\[
\dot{w}(\xi^{\star},t)=\left(y^{\star}-u_{0}(\xi^{\star})\right)\dot{\lambda}(t)-\ddot{w}_{h}(\xi^{\star},t)\lambda(0)+\int_{0}^{t}\dddot{w}_{h}(\xi^{\star},t-\tau)\lambda(\tau)\mathrm{d}\tau.
\]
And therefore the switching function (\ref{eq:Trivial SWcond Pert-1})
becomes
\begin{equation}
h=v_{0}-\dot{y}+\ddot{w}_{h}(\xi^{\star},t)\lambda(0)-\int_{0}^{t}\dddot{w}_{h}(\xi^{\star},t-\tau)\lambda(\tau)\mathrm{d}\tau.\label{eq:Trivial SWcond Convolution}
\end{equation}

The homogeneous solution $\dot{w}_{h}(\xi,t)$ is found using d'Alembert's
method, which states that there are functions $f$ and $g$ such that
$w_{h}(\xi,t)=f(\xi+t)+g(\xi-t)$. Given the initial conditions (\ref{eq:Trivial HomWave IC}),
we have $f(\xi)+g(\xi)=0$ and
\[
\dot{f}(\xi)-\dot{g}(\xi)=y^{\star}\frac{\sin\pi\xi}{\sin\pi\xi^{\star}}-u_{0}(\xi).
\]
It follows that $g(\xi)=-f(\xi)$ and therefore $2\dot{f}(\xi)=y^{\star}\frac{\sin\pi\xi}{\sin\pi\xi^{\star}}-u_{0}(\xi)$.
Next we define $\varphi(\xi)=\ddot{f}(\xi)$, and we get the solution
(for the acceleration) in the form
\begin{align*}
\ddot{w}_{h}(\xi,t) & =\varphi(\xi+t)-\varphi(\xi-t),
\end{align*}
where
\[
\varphi(\xi)=\frac{1}{2}\left(\pi y^{\star}\frac{\cos\pi\xi}{\sin\pi\xi^{\star}}+H(\xi-\xi^{\star})-\left(1-\xi^{\star}\right)\right),\;\xi\in[0,1].
\]
Now we evaluate the boundary conditions, that is $\ddot{w}_{h}(0,t)=\ddot{w}_{h}(1,t)=0$,
which gives $\varphi(\xi)=\varphi(-\xi)$ and $\varphi(1+\xi)=\varphi(1-\xi)$
and by recursion yields
\begin{equation}
\varphi(\xi)=\varphi(2k+\xi),\;\varphi(\xi)=\varphi(2k-\xi)\label{eq:Trivial Phi BC Recursion}
\end{equation}
for $k\in\mathbb{Z}$. Using the initial condition (\ref{eq:Trivial HomWave IC})
and equation (\ref{eq:Trivial Phi BC Recursion}) we find that 
\[
\varphi(\xi)=\frac{1}{2}\left(\pi y^{\star}\frac{\cos\pi\xi}{\sin\pi\xi^{\star}}+H(\left(\xi\mod2\right)-\xi^{\star})+H\left(2-\xi^{\star}-\left(\xi\mod2\right)\right)-2+\xi^{\star}\right),\;\xi\in\mathbb{R}.
\]
We can now expand that 
\begin{align}
\ddot{w}_{h}(\xi^{\star},t) & =-\pi y^{\star}\sin\pi t+\frac{1}{2}H\left(2-2\xi^{\star}-\left(t\mod2\right)\right)-\frac{1}{2}H\left(\left(t\mod2\right)-2\xi^{\star}\right).\label{eq:Trivial wh ddot}
\end{align}
Then the third derivative that appears in the convolution can be written
as
\[
\dddot{w}_{h}(\xi^{\star},t)=-\pi^{2}y^{\star}\cos\pi t+\sum_{k=-\infty}^{\infty}\left(\delta(t-2k)-\frac{1}{2}\delta(t+2\xi^{\star}-2k)-\frac{1}{2}\delta(t-2\xi^{\star}-2k)\right),
\]
where $\delta$ is the Dirac-delta distribution. Due to the convolution
integral (\ref{eq:Trivial Convolution}), we are only taking into
account past values of $\lambda$, which yields 
\begin{equation}
h=v_{0}-\dot{y}-\beta(t),\label{eq:Trivial SW Final-1}
\end{equation}
where 
\begin{multline*}
\beta(t)=\frac{1}{2}\lambda(t)+\sum_{k=1}^{2k<t}\lambda(t-2k)-\frac{1}{2}\sum_{k=0}^{2k<t+2\xi^{\star}}\lambda(t+2\xi^{\star}-2k)\\
-\frac{1}{2}\sum_{k=1}^{2k<t-2\xi^{\star}}\lambda(t-2\xi^{\star}-2k)+\pi^{2}y^{\star}\int_{0}^{t}\cos\pi(t-\tau)\lambda(\tau)\mathrm{d}\tau-\ddot{w}_{h}(\xi^{\star},t)\lambda(0)
\end{multline*}
We can also transform the last remaining integral into a differential
equation by introducing 
\[
\kappa=\pi^{2}y^{\star}\int_{0}^{t}\cos\pi(t-\tau)\lambda(\tau)\mathrm{d}\tau,
\]
which then gives the initial value problem
\begin{equation}
\ddot{\kappa}=\pi^{2}\left(y^{\star}\dot{\lambda}-\kappa\right),\;\kappa(0)=0,\,\dot{\kappa}(0)=\pi^{2}y^{\star}\lambda(0).\label{eq:Trivial Kappa wIC}
\end{equation}
Note that the harmonic term in (\ref{eq:Trivial wh ddot}) is canceled
by the homogeneous solution of (\ref{eq:Trivial Kappa wIC}), hence

\begin{multline}
\beta(t)=\frac{1}{2}\lambda(t)+\sum_{k=1}^{2k<t}\lambda(t-2k)-\frac{1}{2}\sum_{k=0}^{2k<t-2\xi^{\star}}\lambda(t+2\xi^{\star}-2k)-\frac{1}{2}\sum_{k=1}^{2k<t+2\xi^{\star}}\lambda(t-2\xi^{\star}-2k)+\kappa\\
\quad-\frac{1}{2}\left(H\left(2-2\xi^{\star}-\left(t\mod2\right)\right)-H\left(\left(t\mod2\right)-2\xi^{\star}\right)\right)\lambda(0)\label{eq:Trivial BetaFin}
\end{multline}
and the initial condition of (\ref{eq:Trivial Kappa wIC}) vanishes
\begin{equation}
\ddot{\kappa}=\pi^{2}\left(y^{\star}\dot{\lambda}-\kappa\right),\;\kappa(0)=0,\,\dot{\kappa}(0)=0.\label{eq:Trivial Kappa noIC}
\end{equation}

The switching function (\ref{eq:Trivial SW Final-1}) with (\ref{eq:Trivial BetaFin})
and (\ref{eq:Trivial Kappa noIC}) takes into account the full perturbation
(\ref{eq:Trivial Pert Wave-1}) exactly. If we are seeking to solve
for a finite time interval, infinitely long delays in (\ref{eq:Trivial BetaFin})
can be neglected. In case of very short simulation on the interval
$0\le t<\min\left(2\xi^{\star},2-2\xi^{\star}\right)$ it is sufficient
to use
\[
\beta(t)=\frac{1}{2}\lambda(t)+\kappa(t)-\frac{1}{2}\lambda(0),
\]
which then yields 
\[
h=v_{0}-\dot{y}-\frac{1}{2}\lambda-\kappa+\frac{1}{2}\lambda(0),
\]
which is the result we sought.

\section{\label{sec:Proof-Skeleton-Unique}Proof of theorem \ref{thm:ExtUnique}}

The following proof of theorem \ref{thm:ExtUnique} investigates whether
a trajectory approaching $\Sigma_{0}^{\pm}$ can be continued uniquely
after reaching $\Sigma_{0}^{\pm}$ in the two cases set out by the
theorem.
\begin{svmultproof}
Both of the equations (\ref{eq:RED-SmoothEq}) and (\ref{eq:RED-SwitchingEq})
that govern the dynamics on the two sides of $\Sigma_{0}^{\pm}$ already
have unique solutions. We need to exclude the possibility that a trajectory
can be continued by both equations (\ref{eq:RED-SmoothEq}) and (\ref{eq:RED-SwitchingEq})
simultaneously and also exclude the existence of a sliding trajectory
on $\Sigma_{0}^{\pm}$. We denote the solution of (\ref{eq:RED-SmoothEq})
by $(\boldsymbol{\eta}(t),\lambda^{\star})$, and the solution of
(\ref{eq:RED-SwitchingEq}) by $(\boldsymbol{\sigma}(t),\lambda(t))$
either of which can form $\mathcal{T}$. 

First we prove case 1. The speed of solutions relative to $\Sigma_{0}^{\pm}$
on the two sides of $\Sigma_{0}^{\pm}$ are given by $\frac{\mathrm{d}}{\mathrm{d}t}h_{0}(\boldsymbol{\eta}(t),\lambda^{\star})$
and $\dot{\lambda}$, respectively. Trajectories cross $\Sigma_{0}^{\pm}$
if the signs of these two quantities are the same. We calculate that
\begin{equation}
\frac{\mathrm{d}}{\mathrm{d}t}h_{0}(\boldsymbol{\eta}(t),\lambda^{\star})=D_{1}h_{0}(\boldsymbol{\eta}(t),\lambda^{\star})\boldsymbol{f}(\boldsymbol{\eta}(t),\lambda^{\star})\label{eq:hDeriSmooth}
\end{equation}
and rearrange equation (\ref{eq:RED-hDerivative}) into
\begin{equation}
-D_{2}h_{0}(\boldsymbol{\sigma}(t),\lambda(t))\dot{\lambda}(t)=D_{1}h_{0}(\boldsymbol{\sigma}(t),\lambda(t))\boldsymbol{f}(\boldsymbol{\sigma}(t),\lambda(t)).\label{eq:lambdaDeriSwitch}
\end{equation}
At $t=0$ the right sides of (\ref{eq:hDeriSmooth}) and (\ref{eq:lambdaDeriSwitch})
are equal. Using assumption (\ref{eq:UniquenessCondition}) we infer
that $\frac{\mathrm{d}}{\mathrm{d}t}h_{0}(\boldsymbol{\eta}(t),\lambda^{\star})$
and $\dot{\lambda}$ have the same sign at $t=0$, hence $\mathcal{T}$
has a unique continuation transversely through $\Sigma_{0}^{\pm}$.

We now show case 2 of the theorem. Assume that $\mathcal{T}$ is tangent
to $\Sigma_{0}^{\pm}$ to order $\ell-1$. This means that either
\begin{equation}
\frac{\mathrm{d}^{k}}{\mathrm{d}t^{k}}\lambda(t)\vert_{t=0}=0,\;0<k<\ell\label{eq:ExtUni-SwitchTangency}
\end{equation}
if $\mathcal{T}$ is a trajectory of (\ref{eq:RED-SwitchingEq}) or
\begin{equation}
\frac{\mathrm{d}^{k}}{\mathrm{d}t^{k}}h_{0}(\boldsymbol{\eta}(t),\lambda^{\star})\vert_{t=0}=0,\;0<k<\ell\label{eq:ExtUni-SmoothTangency}
\end{equation}
if $\mathcal{T}$ is a trajectory of (\ref{eq:RED-SmoothEq}). 

We first assume that $\mathcal{T}$ is a trajectory of (\ref{eq:RED-SwitchingEq})
and (\ref{eq:ExtUni-SwitchTangency}) holds. Let us consider
\begin{equation}
\frac{\mathrm{d}^{k}}{\mathrm{d}t^{k}}h_{0}(\boldsymbol{\sigma}(t),\lambda(t))\vert_{t=0}=\left.\sum_{j=0}^{k}\binom{k}{j}\frac{\partial^{k}}{\partial\tau^{j}\partial\vartheta^{k-j}}h_{0}(\boldsymbol{\sigma}(\tau),\lambda(\vartheta))\right|_{\tau=\vartheta=0}=0,\label{eq:ExtUni-full-kth-deri}
\end{equation}
which is the constraint that keeps the trajectory on $\Sigma$ (cf.
(\ref{eq:RED-hDerivative})). Any derivative of order $j$ with respect
$\vartheta$ in formula (\ref{eq:ExtUni-full-kth-deri}) includes
a $\frac{\mathrm{d}^{j}}{\mathrm{d}t^{j}}\lambda(t)$ factor. Using
(\ref{eq:ExtUni-SwitchTangency}) we can simplify (\ref{eq:ExtUni-full-kth-deri})
to
\begin{equation}
\frac{\mathrm{d}^{k}}{\mathrm{d}t^{k}}h_{0}(\boldsymbol{\sigma}(t),\lambda(t))\vert_{t=0}=\frac{\mathrm{d}^{k}}{\mathrm{d}t^{k}}h_{0}(\boldsymbol{\sigma}(t),\lambda^{\star})\vert_{t=0}=0,\;0<k<\ell\label{eq:ExtUni-k-th-deri}
\end{equation}
and

\begin{equation}
\frac{\mathrm{d}^{\ell}}{\mathrm{d}t^{\ell}}h_{0}(\boldsymbol{\sigma}(t),\lambda(t))\vert_{t=0}=\frac{\mathrm{d}^{\ell}}{\mathrm{d}t^{\ell}}h_{0}(\boldsymbol{\sigma}(t),\lambda^{\star})\vert_{t=0}+D_{2}h_{0}(\boldsymbol{y}^{\star},\lambda^{\star})\frac{\mathrm{d}^{\ell}}{\mathrm{d}t^{\ell}}\lambda(t)\vert_{t=0}=0.\label{eq:ExtUni-ell-th-deri}
\end{equation}

We now show that 
\begin{equation}
\frac{\mathrm{d}^{k}}{\mathrm{d}t^{k}}h_{0}(\boldsymbol{\sigma}(t),\lambda^{\star})\vert_{t=0}=\frac{\mathrm{d}^{k}}{\mathrm{d}t^{k}}h_{0}(\boldsymbol{\eta}(t),\lambda^{\star}),\;0<k\le\ell.\label{eq:ExtUni-dummy-1}
\end{equation}
The left side of (\ref{eq:ExtUni-dummy-1}) is an algebraic expression
of $D_{1}^{j}h_{0}(\boldsymbol{y}^{\star},\lambda^{\star})$ and $D^{j}\boldsymbol{\sigma}(0)$,
$0<j\le k$. Also, $D^{j}\boldsymbol{\sigma}(0)$ can be written as
\[
D^{j}\boldsymbol{\sigma}(\tau)=\left.\sum_{l=0}^{j-1}\binom{j-1}{l}\frac{\partial^{j-1}}{\partial\tau^{l}\partial\vartheta^{j-l-1}}\boldsymbol{f}(\boldsymbol{\sigma}(\tau),\lambda(\vartheta))\right|_{\tau=\vartheta=0},
\]
which depends on derivatives of $\lambda$ up to order $j-1$ that
are all zero according to (\ref{eq:ExtUni-SwitchTangency}). Therefore
none of $D^{j}\boldsymbol{\sigma}(0)$, $0<j\le k$ depends on the
non-zero $\frac{\mathrm{d}^{k}}{\mathrm{d}t^{k}}\lambda(0)$. The
same holds true for the right side of (\ref{eq:ExtUni-dummy-1}) where
$\lambda$ is assumed to be constant, which proves (\ref{eq:ExtUni-dummy-1}).
Substituting (\ref{eq:ExtUni-dummy-1}) into (\ref{eq:ExtUni-k-th-deri})
implies that (\ref{eq:ExtUni-SmoothTangency}) follows from (\ref{eq:ExtUni-SwitchTangency}).
Further, substituting (\ref{eq:ExtUni-dummy-1}) into (\ref{eq:ExtUni-ell-th-deri})
we find that 
\begin{equation}
\frac{\mathrm{d}^{\ell}}{\mathrm{d}t^{\ell}}h_{0}(\boldsymbol{\sigma}(t),\lambda(t))\vert_{t=0}=\frac{\mathrm{d}^{\ell}}{\mathrm{d}t^{\ell}}h_{0}(\boldsymbol{\eta}(t),\lambda^{\star})\vert_{t=0}+D_{2}h_{0}(\boldsymbol{y}^{\star},\lambda^{\star})\frac{\mathrm{d}^{\ell}}{\mathrm{d}t^{\ell}}\lambda(t)\vert_{t=0}=0.\label{eq:ExtUni-ell-th-deri-mixed}
\end{equation}
So far we have shown that if vector field (\ref{eq:RED-SwitchingEq})
is tangent of order $\ell-1$ to $\Sigma_{0}^{\pm}$ at $(\boldsymbol{y}^{\star},\lambda^{\star})$
then so is (\ref{eq:RED-SmoothEq}) and the orientation of the tangencies
are the same. We now show that this sufficient condition is also necessary. 

Using equation (\ref{eq:ExtUni-ell-th-deri-mixed}) for $\ell=1$
does not require assumption (\ref{eq:ExtUni-SwitchTangency}). It
directly follows from equation (\ref{eq:ExtUni-ell-th-deri-mixed})
that 
\[
\frac{\mathrm{d}}{\mathrm{d}t}h_{0}(\boldsymbol{\eta}(t),\lambda^{\star})\vert_{t=0}=0\;\implies\;\frac{\mathrm{d}}{\mathrm{d}t}\lambda(t)=0.
\]
Now knowing that $\frac{\mathrm{d}}{\mathrm{d}t}\lambda(t)=0$ we
can apply (\ref{eq:ExtUni-ell-th-deri-mixed}) for $\ell=2$ and conclude
that 
\[
\frac{\mathrm{d}^{2}}{\mathrm{d}t^{2}}h_{0}(\boldsymbol{\eta}(t),\lambda^{\star})\vert_{t=0}=0\;\implies\;\frac{\mathrm{d}^{2}}{\mathrm{d}t^{2}}\lambda(t)=0.
\]
Repeating this procedure a sufficient number of times shows that (\ref{eq:ExtUni-SmoothTangency})
implies (\ref{eq:ExtUni-SwitchTangency}).

In summary, (\ref{eq:ExtUni-SmoothTangency}) holds if and only if
(\ref{eq:ExtUni-SwitchTangency}) holds and we have the equality
\begin{equation}
\frac{\mathrm{d}^{\ell}}{\mathrm{d}t^{\ell}}h_{0}(\boldsymbol{\eta}(t),\lambda^{\star})\vert_{t=0}=-D_{2}h_{0}(\boldsymbol{y}^{\star},\lambda^{\star})\frac{\mathrm{d}^{\ell}}{\mathrm{d}t^{\ell}}\lambda(t)\vert_{t=0}.\label{eq:ExtUni-final-deri}
\end{equation}

Using assumption (\ref{eq:UniquenessCondition}) and equation (\ref{eq:ExtUni-final-deri})
we conclude that the order and orientation of the tangency is the
same on both side of $\Sigma_{0}^{\pm}$. If $\ell$ is odd, trajectory
$\mathcal{T}$ passes through $\Sigma_{0}^{\pm}$ at the tangency.
If $\ell$ is even, the trajectory continues on the same side of $\Sigma_{0}^{\pm}$
and there is no joining trajectory from the other side of $\Sigma_{0}^{\pm}$.
Conditions (\ref{eq:UniquenessCondition}) and (\ref{eq:ExtUni-TangencyCondition})
also imply that either case 1 or 2 holds for points on $\Sigma_{0}^{\pm}$
in a sufficiently small open neighborhood of $(\boldsymbol{y}^{\star},\lambda^{\star})$.
This excludes cases where trajectories are forced onto $\Sigma_{0}^{\pm}$
for a non-zero length of time and implies that there cannot be trajectories
joining $(\boldsymbol{y}^{\star},\lambda^{\star})$ from within $\Sigma_{0}^{\pm}$.
Therefore, there is a unique continuation of $\mathcal{T}$ for $t>0$
(or $t<0$) sufficiently small.
\end{svmultproof}

\section{\label{sec:UniquenessProof}Proofs of lemma \ref{lem:INF-discontinuity-gap}
and theorem \ref{thm:Cont-Uniqueness}}

There are three steps to the proof of theorem \ref{thm:Cont-Uniqueness}.
First we derive a differential equation for $\lambda$ from the algebraic
constraint 
\[
h(\boldsymbol{W}(\boldsymbol{y},\lambda)+\boldsymbol{z})=0
\]
 by differentiation, which follows from lemma \ref{lem:INF-discontinuity-gap}.
At this step we claim continuity of $\dot{\lambda}$. By investigating
the resulting equation we move to step two and establish that $\dot{\lambda}$
is indeed continuous in $\Sigma$. In the final step we show that
trajectories transverse to $\Sigma^{\pm}$ cross $\Sigma^{\pm}$ when
$d^{\pm}(\boldsymbol{y},\boldsymbol{z},\lambda)>0$, which concludes
the proof.

For convenience we copy here lemma \ref{lem:INF-discontinuity-gap}.
\begin{lemma}[Lemma \ref{lem:INF-discontinuity-gap}]
Assume \nameref{A-EvolutionOperator} and that $\lambda$ is continuously
differentiable and $\boldsymbol{y}$, $\boldsymbol{z}$ satisfy the
differential equations
\begin{equation}
\left.\begin{array}{rl}
\dot{\boldsymbol{y}} & =\boldsymbol{f}(\boldsymbol{y},\lambda)\\
\dot{\boldsymbol{z}} & =\boldsymbol{A}_{1}(\boldsymbol{y},\lambda)\boldsymbol{z}-D_{2}\boldsymbol{W}(\boldsymbol{y},\lambda)\dot{\lambda}
\end{array}\right\} \label{eq:INF-lin-equation}
\end{equation}
on the interval $t\in[s,s+\epsilon)$, $\epsilon>0$ with an initial
condition $\boldsymbol{y}(s)\in G$, $\boldsymbol{z}(s)\in\boldsymbol{\mathcal{D}}$.
Then the right-side derivative of $h$ as a function of time is calculated
as
\begin{equation}
\frac{\mathrm{d}}{\mathrm{d}t^{+}}h(\boldsymbol{W}(\boldsymbol{y},\lambda)+\boldsymbol{z})=Dh(\boldsymbol{W}(\boldsymbol{y},\lambda)+\boldsymbol{z})\cdot D_{1}\boldsymbol{W}(\boldsymbol{y},\lambda)\boldsymbol{f}(\boldsymbol{y},\lambda)-d^{\pm}(\boldsymbol{y},\boldsymbol{z},\lambda)\dot{\lambda}+\boldsymbol{A}_{1}(\boldsymbol{y},\lambda)\boldsymbol{z},\label{eq:INF-hDerivative-1}
\end{equation}
where
\begin{equation}
d^{\pm}(\boldsymbol{y},\boldsymbol{z},\lambda)=\lim_{\delta\downarrow0}Dh(\boldsymbol{W}(\boldsymbol{y},\lambda)+\boldsymbol{z})\cdot\left(\boldsymbol{K}(t+\delta,t)-D_{2}\boldsymbol{W}(\boldsymbol{y},\lambda)\right).\label{eq:dPlusMinusDefinition-1}
\end{equation}
\end{lemma}
\begin{svmultproof}
Consider equation (\ref{eq:INF-lin-equation}) with solution $\boldsymbol{y}$
and $\boldsymbol{z}$ on the interval $[s,s+\epsilon)$. We start
with the expression 
\begin{equation}
\left.\frac{\mathrm{d}}{\mathrm{d}t^{+}}h(\boldsymbol{W}(\boldsymbol{y},\lambda)+\boldsymbol{z})\right|_{t=s}=Dh(\boldsymbol{W}(\boldsymbol{y},\lambda)+\boldsymbol{z})\cdot\left(D_{1}\boldsymbol{W}(\boldsymbol{y},\lambda)\boldsymbol{f}(\boldsymbol{y},\lambda)+D_{2}\boldsymbol{W}(\boldsymbol{y},\lambda)\dot{\lambda}+\left.\frac{\mathrm{d}}{\mathrm{d}t^{+}}\boldsymbol{z}\right|_{t=s}\right).\label{eq:INF-A-formal-h-deri}
\end{equation}
and show that it can be transformed into (\ref{eq:INF-hDerivative-1}).
Let us define $\overline{\boldsymbol{x}}=\boldsymbol{W}(\boldsymbol{y},\lambda)+\boldsymbol{z}\vert_{t=s}$.
The only unresolved term 
\[
Dh(\overline{\boldsymbol{x}})\cdot\left.\frac{\mathrm{d}}{\mathrm{d}t^{+}}\boldsymbol{z}\right|_{t=s}
\]
is obtained by taking the derivative of
\begin{equation}
Dh(\overline{\boldsymbol{x}})\cdot\boldsymbol{z}(t)=Dh(\overline{\boldsymbol{x}})\cdot\left(\boldsymbol{U}(t,s)\boldsymbol{z}(s)-\int_{s}^{t}\boldsymbol{K}(t,\tau)\dot{\lambda}(\tau)\mathrm{d}\tau\right)\label{eq:INF-Projected-z(t)}
\end{equation}
while assuming constant $\overline{\boldsymbol{x}}$. The integral
in (\ref{eq:INF-Projected-z(t)}) is well defined in the Riemann sense,
because of assumption \nameref{A-PertTrajectories} and because $\dot{\lambda}$
is continuous. To simplify notation we define 
\begin{equation}
\eta(t,s)=Dh(\overline{\boldsymbol{x}})\cdot\boldsymbol{K}(t,s),\label{eq:INF-A-eta-definition}
\end{equation}
so that 
\begin{equation}
\frac{\mathrm{d}}{\mathrm{d}t^{+}}Dh(\overline{\boldsymbol{x}})\cdot\boldsymbol{z}(t)\vert_{t=s}=Dh(\overline{\boldsymbol{x}})\cdot\boldsymbol{A}_{1}(\boldsymbol{y},\lambda)\boldsymbol{z}(s)-\left.\frac{\mathrm{d}}{\mathrm{d}t^{+}}\left(\int_{s}^{t}\eta(t,\tau)\dot{\lambda}(\tau)\mathrm{d}\tau\right)\right|_{t=s}.\label{eq:INF-A-deri-zt}
\end{equation}
Assumptions \nameref{A-PertTrajectories} and \nameref{A-EvolutionOperator}
imply that $\eta(t,s)$ is continuous for $t>s$, but also allow a
discontinuity at $t=s$, which needs to be taken into account. Differentiating
the convolution in (\ref{eq:INF-A-deri-zt}) yields 
\begin{equation}
\frac{\mathrm{d}}{\mathrm{dt}}\int_{s}^{t}\eta(t,\tau)\dot{\lambda}(\tau)\mathrm{d}\tau=\eta(t,t)\dot{\lambda}(t)+\int_{s}^{t}D_{1}\eta(t,\tau)\dot{\lambda}(\tau)\mathrm{d}\tau.\label{eq:INF-A-deri-convolution}
\end{equation}
The integral on the right side of (\ref{eq:INF-A-deri-convolution})
is approximated by a Riemann sum
\[
\int_{s}^{t}D_{1}\eta(t,\tau)\dot{\lambda}(\tau)\mathrm{d}\tau\approx\sum_{k=0}^{m-1}\delta D_{1}\eta(t,s+k\delta+c_{k})\dot{\lambda}(s+k\delta+c_{k})+\mathcal{O}(\delta m^{-1}),
\]
where $\delta=\left(t-s\right)/m$, $m>1$ is an integer and $c_{k}\in(0,\delta)$.
We use finite differences to approximate the derivative
\[
D_{1}\eta(t,s+k\delta+c_{k})\approx\delta^{-1}\left(\eta(t+c_{k},s+k\delta+c_{k})-\eta(t-\delta+c_{k},s+k\delta+c_{k})\right).
\]
The scheme of the finite difference is such that for $k=m-1$ the
second and first argument of $\eta$ are equal in one of the terms,
i.e., $t-\delta+c_{k}=s+k\delta+c_{k}$, which takes into account
the discontinuity of $\eta$. If this discontinuity is not taken into
account, the integral in the limit $t\downarrow s$ would vanish.
In summary we have the integral
\[
\int_{s}^{t}D_{1}\eta(t,\tau)\dot{\lambda}(\tau)\mathrm{d}\tau\approx\sum_{k=0}^{m-1}\left(\eta(t+c_{k},s+k\delta+c_{k})-\eta(t-\delta+c_{k},s+k\delta+c_{k})\right)\dot{\lambda}(s+k\delta+c_{k}).
\]
Taking the limit $t\downarrow s$, is the same as $\delta\downarrow0$,
hence we calculate that 
\[
\lim_{\delta\downarrow0}\eta(t+c_{k},s+k\delta+c_{k})-\eta(t-\delta+c_{k},s+k\delta+c_{k})=\begin{cases}
0 & k<m-1\\
\lim_{t\downarrow s}\eta(t,s)-\eta(s,s) & k=m-1
\end{cases},
\]
which implies that
\[
\lim_{t\downarrow s}\int_{s}^{t}D_{1}\eta(t,\tau)\dot{\lambda}(\tau)\mathrm{d}\tau=\left(\lim_{t\downarrow s}\eta(t,s)-\eta(s,s)\right)\dot{\lambda}(s).
\]
Using the definition of $\eta$ and formula (\ref{eq:INF-A-deri-convolution})
gives
\begin{align*}
\left.\frac{\mathrm{d}}{\mathrm{d}t^{+}}\left(\int_{s}^{t}\eta(t,\tau)\dot{\lambda}(\tau)\mathrm{d}\tau\right)\right|_{t=s} & =\cancel{\eta(s,s)}\dot{\lambda}(s)+\left(\lim_{t\downarrow s}\eta(t,s)-\cancel{\eta(s,s)}\right)\dot{\lambda}(s)\\
 & =\lim_{t\downarrow s}Dh(\overline{\boldsymbol{x}})\cdot\boldsymbol{K}(t,s)\dot{\lambda}(s),
\end{align*}
which in turn is put into (\ref{eq:INF-A-formal-h-deri})
\begin{align*}
\frac{\mathrm{d}}{\mathrm{d}t^{+}}h(\boldsymbol{W}(\boldsymbol{y},\lambda)+\boldsymbol{z}) & =Dh(\boldsymbol{W}(\boldsymbol{y},\lambda)+\boldsymbol{z})\cdot\left(D_{1}\boldsymbol{W}(\boldsymbol{y},\lambda)\boldsymbol{f}(\boldsymbol{y},\lambda)+D_{2}\boldsymbol{W}(\boldsymbol{y},\lambda)\dot{\lambda}+\boldsymbol{A}_{1}(\boldsymbol{y},\lambda)\boldsymbol{z}(s)\right)\\
 & \quad-\lim_{t\downarrow s}Dh(\overline{\boldsymbol{x}})\cdot\boldsymbol{K}(t,s)\dot{\lambda}.
\end{align*}
This proves lemma \ref{lem:INF-discontinuity-gap}.
\end{svmultproof}

We now prove theorem \ref{thm:Cont-Uniqueness}. 
\begin{svmultproof}
Proof of theorem \ref{thm:Cont-Uniqueness}. We want to derive $\dot{\lambda}$
from the constraint $h(\overline{\boldsymbol{x}})=0$, where $\overline{\boldsymbol{x}}=\boldsymbol{W}(\boldsymbol{y},\lambda)+\boldsymbol{z}$.
We take the time derivative $\frac{\mathrm{d}}{\mathrm{d}t^{+}}h(\boldsymbol{W}(\boldsymbol{y},\lambda)+\boldsymbol{z})=0$
and use the expression (\ref{eq:INF-hDerivative-1}) to solve for
\begin{equation}
\dot{\lambda}=\frac{1}{d^{\pm}(\boldsymbol{y},\boldsymbol{z},\lambda)}Dh(\overline{\boldsymbol{x}})\cdot\left(D_{1}\boldsymbol{W}(\boldsymbol{y},\lambda)\boldsymbol{f}(\boldsymbol{y},\lambda)+\boldsymbol{A}_{1}(\boldsymbol{y},\lambda)\boldsymbol{z}(t)\right).\label{eq:INF-A-lambdaDeri}
\end{equation}
This concludes the first step of the proof. We now demonstrate that
$\dot{\lambda}$ is indeed continuous. We only need to recall that
the formal solution for any history of $\lambda$ is
\[
\boldsymbol{z}(t)=\boldsymbol{U}(t,s)\boldsymbol{z}(s)-\int_{s}^{t}\boldsymbol{K}(t,\tau)\dot{\lambda}(\tau)\mathrm{d}\tau,
\]
which allows us to write that
\begin{equation}
Dh(\overline{\boldsymbol{x}})\cdot\boldsymbol{A}_{1}(\boldsymbol{y},\lambda)\boldsymbol{z}(t)=\frac{\mathrm{d}}{\mathrm{d}t}Dh(\overline{\boldsymbol{x}})\cdot\boldsymbol{U}(t,s)\boldsymbol{z}(s)-\int_{s}^{t}\frac{\mathrm{d}}{\mathrm{d}t}Dh(\overline{\boldsymbol{x}})\cdot\boldsymbol{K}(t,\tau)\dot{\lambda}(\tau)\mathrm{d}\tau.\label{eq:INF-A-Dh-zDeri}
\end{equation}
Under the assumptions of theorem \ref{thm:Cont-Uniqueness}, the expression
(\ref{eq:INF-A-Dh-zDeri}) is continuous in $t$ and so is (\ref{eq:INF-A-lambdaDeri}),
which concludes the second part of the proof.

Finally, we need to show that if $Dh(\overline{\boldsymbol{x}})\cdot\left(D_{1}\boldsymbol{W}(\boldsymbol{y},\lambda)\boldsymbol{f}(\boldsymbol{y},\lambda)+\boldsymbol{A}_{1}(\boldsymbol{y},\lambda)\boldsymbol{z}(t)\right)\neq0$,
trajectories cross $\Sigma^{\pm}$. This renders trajectories unique
as they pass through $\Sigma^{\pm}$. Indeed, if $\frac{\mathrm{d}}{\mathrm{d}t^{+}}h$
outside of $\Sigma$ and $\dot{\lambda}$ inside of $\Sigma$, but
right on the boundary $\Sigma^{\pm}$, have the same sign, trajectories
cross $\Sigma^{\pm}$. We only need to evaluate (\ref{eq:INF-hDerivative})
with $\dot{\lambda}=0$ and compare that to formula (\ref{eq:INF-A-lambdaDeri}).
The two values differ by the factor $d^{\pm}(\boldsymbol{y},\boldsymbol{z},\lambda)$,
hence if $d^{\pm}(\boldsymbol{y},\boldsymbol{z},\lambda)>0$ and $Dh(\overline{\boldsymbol{x}})\cdot\left(D_{1}\boldsymbol{W}(\boldsymbol{y},\lambda)\boldsymbol{f}(\boldsymbol{y},\lambda)+\boldsymbol{A}_{1}(\boldsymbol{y},\lambda)\boldsymbol{z}(t)\right)\neq0$,
trajectories cross $\Sigma^{\pm}$ and solutions are unique.
\end{svmultproof}

\end{document}